\documentclass [11pt ]{article}

\usepackage{geometry}
\usepackage[english]{babel}
\usepackage[utf8]{inputenc}
\usepackage{fontenc}
\usepackage{newlfont}
\usepackage{indentfirst}
\usepackage{fancyhdr}
\usepackage{textcomp}
\usepackage{graphicx}
\usepackage{amsmath}
\usepackage{amssymb}
\usepackage{amsthm}
\usepackage{mathrsfs}
\usepackage{float}
\usepackage{multicol}
\usepackage[dvips]{hyperref}
\theoremstyle{plain}
\newtheorem{Teo}{Theorem}[section]
\newtheorem{Prop}[Teo]{Proposition} 
 
\newtheorem{Lem}[Teo]{Lemma}
\theoremstyle{definition}
\newtheorem{Def}{Definition}

\theoremstyle{remark}
\newtheorem{Rem}{Remark}
\newtheorem{Oss}{Observation}

\numberwithin{Teo}{section}
\numberwithin{Rem}{section}
\numberwithin{Def}{section}
\numberwithin{Not}{section}
\numberwithin{Oss}{section}

\title{On a quasi-additivity formula for the capacity on Ahlfors-regular spaces}
\author{Cavina Michelangelo}
\begin{document}
\maketitle
\begin{abstract}
In this work we prove formulas of quasi-additivity for the capacity associated to kernels of radial type in the setting of the boundary of a tree structure and in the setting of compact Ahlfors-regular spaces. We also define a notion of harmonic extension, to one additional variable, of a function defined over a compact Ahlfors-regular space, and we prove a result of tangential convergence of the harmonic extension to the values at the boundary.\\
\end{abstract}

\tableofcontents

\section*{Introduction}\addcontentsline{toc}{section}{Introduction}

Potential theory was born out of the theory of electrostatic potential, from the work of C.F. Gauss. The first notion of capacity dates back to the 1830's and it is the notion of electrostatic capacitance of a compact set $K \subseteq R^3$. We consider a distribution of charge $\mu$ over $\mathbb R^3$, which is a positive measure on $\mathbb R^3$. The amount of charge on a conductor $K \subseteq \mathbb R^3$ is equal to $\mu(K)$. Given a distribution of charge $\mu$, the electrostatic potential associated to $\mu$ at the point $y \in \mathbb R^3$ is defined by 
\begin{equation}
V_{\mu}(y):=\int_{\mathbb R^3} \frac{d\mu(x)}{|x-y|},
\end{equation}
and the energy associated to $\mu$ is defined by 
\begin{equation}
\mathcal E(\mu):=\int_{\mathbb R^3}\int_{\mathbb R^3}\frac{1}{|x-y|}d\mu(x)d\mu(y).
\end{equation}
A charge distribution $\tilde \mu$ is called an equilibrium charge distribution for a compact set $K\subseteq \mathbb R^3$ if $\tilde \mu$ is supported on $K$ and if the potential $V_{\tilde \mu}(y)$ is equal to a constant value $V_{\tilde \mu}$ for all $y \in K$, except for a ``small" set of exceptional points. An equilibrium charge distribution $\tilde \mu$ minimizes the value of the energy $\mathcal E(\mu)$ for charge distributions $\mu$ such that $\mu(K)=\tilde \mu(K)$. Using these notions we define the electrostatic capacitance of a compact set $K \subseteq \mathbb R^3$ to be ratio between the amount of charge in $K$ and the value of the electrostatic potential associated to an equilibrium charge distribution, i.e.
\begin{equation}
C_{\text{E.S.}}(K):= \frac{\tilde \mu(K)}{V_{\tilde \mu}}.
\end{equation}

Potential theory has always been strictly connected to the theory of Hilbert spaces and harmonic analysis. In the 1830's Gauss proved the existence of equilibrium potentials by minimizing a quadratic integral, the energy (see \cite{Gauss}). The same result was proved with modern mathematical rigor by O. Frostman in the 1930's (see \cite{Frostman}). This theme kept on growing during the 1940's, and it was made expecially clear in the work of H. Cartan (see \cite{Cartan}).\\
In the following decades a mathematical concept of potential theory, disconnected from the theory of the electrostatics, was developed. G. Choquet (see \cite{Choquet}), introduced a definition of capacity and capacitability in 1950, and in the 1970's Frostman developed the first rigorous mathematical definition of potential theory. In the early 1960's Maz'ya developed a non linear potential theory which is connected to the theory of function spaces.\\
In 1990 D.R. Adams and L.I. Hedberg developed an axiomatic definition of potential theory for metric measure spaces (see \cite{AE}) which allows to prove many results in a very general setting. In this work, we refer to \cite{AE} for the definitions and main theorems about potential theory. In the 1990's J. Heinonen and P. Koskela developed a potential theory on metric spaces based on the notions of rectifiable curves and of ``weak gradients" (see \cite{HK}).\\

In this work we study problems from the theory of the potential on metric measure spaces, in the case of tree boundaries and in the case of compact Ahlfors-regular spaces.\\
To define a notion of potential theory on a metric measure space $(X,d,m)$ we consider a kernel
\begin{equation}
K: X \times X \longrightarrow \mathbb R_0^+
\end{equation}
satisfying proper conditions (see \cite{AH}), and we define the $L^p$ capacity (associated to $K$) of a compact set $E\subseteq X$ in the following way:
\begin{equation}\label{L^p capacity definition in intro}
C_{K,p}(E):= \inf \left\{\| f\|^p_{L^p(X,m)} \; \big| \; f \in L^p(X,m), \; K*f(x)\geq 1 \; \forall x \in E\right \},
\end{equation}
where $K*f$ denotes the potential of $f$, and it is defined by
\begin{equation}
K*f(x):= \int_X K(x,y)f(y)dm(y).
\end{equation}

A well-known result from potential theory (see \cite{AH}) is that the capacity $C_{K,p}(E)$ of a set $E \subseteq X$ is a sub additive map
\begin{equation*}
C_{K,p}: \{\text{Borel subsets of $X$}\} \longrightarrow [0,+\infty],
\end{equation*}
hence we have
\begin{equation}
 C_{K,p}\left(\bigcup_{j \in \mathbb N}E_j\right) \leq \sum_{j \in \mathbb N} C_{K,p}(E_j) 
\end{equation}
for all $E_j \subseteq X$.\\
In general the capacity is very far from superadditive, an example of this property comes from the theory of electrostatic potential. Suppose $B(0,r) \subseteq R^3$ is the metric ball centered at 0 of radius $r>0$. When we consider the electrostatic capacitance $C_{\text{E.S.}}(B(0,r))$ we have 
\begin{equation}
C_{\text{E.S.}}(B(0,r))=C_{\text{E.S.}}(\mathring{B(0,r)})=C_{\text{E.S.}}(\partial B(0,r)).
\end{equation}
If we consider an infinite family $\{r_j\}$ such that $\frac r 2 <r_j<r$ and $r_{j_1}\neq r_{j_2}$ for all $j_1 \neq j_2$, then we have $\partial B(0,r_{j_1}) \cap \partial B(0,r_{j_2})=\emptyset$ for all $j_1 \neq j_2$, and $\partial B(0,r_{j_1}) \subseteq B(0,r)$ for all $j$, and 
\begin{equation}
\begin{aligned}
C_{\text{E.S.}}(B(0,r))<+\infty,\\
\sum_j C_{\text{E.S.}}(\partial B(0,r_j))=+\infty,
\end{aligned}
\end{equation}
which is a counter example to the quasi-additivity of the electrostatic capacitance.\\
A common problem in harmonic analysis is to find notions of ``properly separated" families of sets $\{E_j\}_{j \in \mathbb N}$ implying that
\begin{equation}
\sum_{j \in \mathbb N} C_{K,p}(E_j)\leq C \cdot C_{K,p}\left(\bigcup_{j \in \mathbb N}E_j\right),
\end{equation}
for an appropriate constant $C\geq 1$.\\
Results of this kind have been developed through the years, including results using a notion separation of the sets $E_j$ based on Whitney decompositions (see \cite{A1}, \cite{AE}).\\

In the first section of this work we consider the quasi-additivity formula in the article \cite{AB} from H. Aikawa and A.A. Borichev, and we extend it to the setting of the theory of the potential on tree boundaries. Then, we use the results from the theory of dyadic potential (see \cite{ARSW1}) to find a quasi-additivity formula for compact Ahlfors-regular spaces.\\
Aikawa and Borichev consider the setting of $\mathbb R^n$ and define a theory of the potential based on radial non-increasing kernels. Given a non-increasing function $\phi(r)>0$ for $r>0$ such that
\begin{equation}
\lim_{r\rightarrow 0} \phi(r)=+\infty, \quad \lim_{r \rightarrow +\infty} \phi(r)=0,
\end{equation}
and such that
\begin{equation}
\int_{0}^R \phi(r)r^{n-1}dr <+\infty,
\end{equation}
they consider the radial kernel
\begin{equation}\label{Aikawa Borichev kernels definition intro}
K(x,y):=K(\phi(\lvert x-y \rvert))
\end{equation}
and consider the $L^p$-capacity $C_{K,p}(E)$ on $\mathbb R^n$ associated to the kernel $K$ defined by (\ref{L^p capacity definition in intro}), i.e
\begin{equation}
C_{K,p}(E):= \inf \left\{\| f\|^p_{L^p(\mathbb R^n)} \; \big| \; f \in L^p(\mathbb R^n), \; K*f(x)\geq 1 \; \forall x \in E\right \},
\end{equation}
where the potential $K*f$ is defined by
\begin{equation}
K*f(x):= \int_{\mathbb R^n} K(x,y)f(y)dy.
\end{equation}

Aikawa and Borichev define the following ``Aikawa-Borichev radius", which gives a notion of separation for sets which yields a quasi-additivity formula for the capacity $C_{K,p}$.
\begin{Def}[Aikawa-Borichev radius]\label{Aikawa Borichev radius definition in intro}
Let $x \in \mathbb R^n$, $r>0$. We define $\eta_p(r)$ to be the radius $R>0$ such that
\begin{equation}
 m(B(x,R)) = C_{K,p}(B(x,r)),
\end{equation}
and we also define
\begin{equation}
\eta^*_p(r):= \max\{2r,\eta_p(r)\}.
\end{equation}
\end{Def}
Aikawa and Borichev use this definition to formulate the following result about a quasi-additivity formula for the capacity $C_{K,p}$.
\begin{Teo}[Aikawa-Borichev quasi-additivity in $\mathbb R^n$, see Theorem 5, \cite{AB}]\label{theorem 5 intro}
There exists a universal constant $A$ depending only on $K$, $n$ and $p$ such that: for any countable family $B(x_j,r_j) \subseteq \mathbb R^n$ such that the family $\{B(x_j,\eta^*_p(r_j))\}$ is pairwise disjoint, for all compact sets $E_j$ such that $E_j \subseteq B(x_j,r_j)$, we have
\begin{equation}
\sum_j C_{K,p}(E_j) \leq A \cdot C_{K,p}\bigg(\bigcup_j E_j\bigg).
\end{equation}
\end{Teo}

In our work we extend Theorem \ref{theorem 5 intro} to the setting of tree boundaries for the capacity associated to integrable radial kernels. Let $(\partial T,\rho,m)$ be a metric measure space, where $\partial T$ is the boundary of a rooted tree $T$, $\rho$ is the natural (ultrametric) distance on $\partial T$ and $m$ is a $\sigma$-finite Borel measure. We recall that the distance $\rho$ satisfies
\begin{equation}
\rho(x,y)=\delta^{n(x,y)},
\end{equation}
for a specific $0<\delta<1$, where $n(x,y) \in \mathbb N$ is the number of edges shared by the geodesics connecting $x$ and $y$ to the root $o$. In the following, $\delta$ denotes such value.\\
We consider (with a slight abuse of notation) a radial kernel
\begin{equation*}
K:\partial T \times \partial T \longrightarrow \mathbb R_0^+
\end{equation*}
\begin{equation*}
K(x,y)=K(\rho(x,y)),
\end{equation*}
such that the following integrability condition holds:
\begin{equation}
\max\left(\sup_{x \in \partial T}\int_{\partial T}K(x,y)dm(y),\sup_{y \in \partial T}\int_{\partial T}K(x,y)dm(x)\right)<+\infty.
\end{equation}

According to the definition (\ref{L^p capacity definition in intro}) we define the $L^p$-capacity of a compact set $E\subseteq \partial T$ associated to the kernel $K$ by
\begin{equation}
C_{K,p}(E):=\inf \bigg\{ \|f\|^p_{L^p(X,m)} \; \bigg | \; K * f(x)\geq 1 \; \; \; \forall x \in E\bigg\},
\end{equation}
where
\[ K * f(x) := \int_X K(x,y) f(y) dm(y).\]
Following the definition (\ref{Aikawa Borichev radius definition in intro}) we define the Aikawa-Borichev radius in the setting of tree boundaries by considering
\begin{equation}
\eta_p(x,r):= \inf \bigg \{ \delta^{n-\frac 1 2} \in \mathbb R \; \bigg | \; n \in \mathbb N, \; m(B_{\rho}(x,\delta^{n-\frac 1 2}))\geq C_{K,p}(B_{\rho}(x,r))\bigg \}
\end{equation}
for $x \in \partial T$ and $r>0$, where $B_{\rho}(x,r)$ denotes the metric ball of center $x$ and radius $r$ with respect to the distance $\rho$, and we define
\begin{equation}
\eta^*_p(x,r):= \max\{r,\eta_p(x,r)\}.
\end{equation}
From this definition we prove the following quasi-additivity formula for the capacity $C_{K,p}$, which extends Theorem \ref{theorem 5 intro} to the setting of tree boundaries.
\begin{Teo}[Quasi-additivity of the capacity for tree boundaries]\label{quasi additivity tree boundaries theorem intro}
There exists a universal constant $A$ depending only on $K$, $\partial T$ and $p$ such that: for any countable family $B_{\rho}(x_j,r_j) \subseteq \partial T$ such that the family $\{B_{\rho}(x_j,\eta^*_p(x_j,r_j))\}$ is pairwise disjoint, for all compact sets $E_j$ such that $E_j \subseteq B_{\rho}(x_j,r_j)$, we have
\begin{equation}
\sum_j C_{K,p}(E_j) \leq A \cdot C_{K,p}\bigg(\bigcup_j E_j\bigg).
\end{equation}
\end{Teo}
We observe that Theorem \ref{quasi additivity tree boundaries theorem intro} holds in a greater level of generality compared to Theorem \ref{theorem 5 intro}, because it does not require monotonicity of the radial kernel $K(x,y)=K(\rho(x,y))$.\\

In the second section of our work we focus on potential theory on compact Ahlfors-regular spaces, and we study the quasi-additivity of the Riesz capacity.\\
Let $(X,d,m)$ be a metric measure space. Let $Q>0$. We say that $X$ is a $Q$-regular Ahlfors space if there exist constants $0<C_1\leq C_2<+\infty$ such that
\begin{equation}
C_1 r^Q \leq m(B(x,r)) \leq C_2 r^Q
\end{equation}
for all $x \in X$, for all $0<r<\text{diam}(X)$. The parameter $Q$ represents the dimension of the space $X$.\\
Given a compact Ahlfors $Q$-regular space $(X,m,d)$ and $1<p<+\infty$ we consider the Riesz kernel
\begin{equation}\label{Riesz kernel Ahlfors regular intro definition}
\begin{aligned}
K_{X,s}:X\times X &\longrightarrow \mathbb R\\
(x,y) &\longmapsto \frac{1}{d(x,y)^{Q \cdot s}},
\end{aligned}
\end{equation}
where $\frac{1}{p'}\leq s <1$.

Following (\ref{L^p capacity definition in intro}) we define the $L^p$-capacity of a compact set $E\subseteq X$ associated to the kernel $K_{X,s}$ by
\begin{equation}\label{Riesz capacity Ahlfors regular intro definition}
C_{K_{X,s},p}(E):= \inf \left\{\| f\|^p_{L^p(X,m)} \; \big| \; f \in L^p(X,m), \; K_{X,s}*f(x)\geq 1 \; \forall x \in E\right \},
\end{equation}
where the Riesz potential $K_{X,s}*f$ is defined by
\begin{equation}
K_{X,s}*f(x):= \int_X K_{X,s}(x,y)f(y)dm(y).
\end{equation}
Following the previous approach we define the Aikawa-Borichev radius for the setting of compact Ahlfors-regular spaces by setting
\begin{equation}
\eta_{X,p}(x,r):= \inf \left\{ R>0 \; | \; m(B_d(x,R))\geq C_{K_{X,s},p}(B_d(x,r))\right\},
\end{equation}
for $x \in X$, $r>0$, and by defining
\begin{equation}
\eta^*_{X,p}(x,r):=\max \{ r,\eta_{X,p}(x,r)\}.
\end{equation}
The second main result of our work is the following quasi-additivity formula for the capacity $ C_{K_{X,s},p}$, which extends Theorem \ref{theorem 5 intro} to the setting of compact Ahlfors-regular spaces.
\begin{Teo}[Quasi-additivity of the capacity for compact Ahlfors-regular spaces]\label{quasi additivity Ahlfors regular theorem intro}
There exists a universal constant $A$ depending only on $K$, $X$, $s$ and $p$ such that: for any countable family $B_{d}(x_j,r_j) \subseteq X$ such that the family $\{B_{d}(x_j,\eta^*_{X,p}(x_j,r_j))\}$ is pairwise disjoint, for all compact sets $E_j$ such that $E_j \subseteq B_{d}(x_j,r_j)$, we have
\begin{equation}
\sum_j C_{K_{X,s},p}(E_j) \leq A \cdot C_{K_{X,s},p}\bigg(\bigcup_j E_j\bigg).
\end{equation}
\end{Teo}

In the third and fourth sections of our work we define a harmonic-type extension to the space $X \times (0,+\infty)$ of the Riesz potential on a compact Ahlfors-regular space $X$, and we study its behaviour at the boundary $X \times \{0\}$.\\
Studying the behaviour of harmonic functions at the boundary is a common problem in harmonic analysis. In the article \cite{AB}, Aikawa and Borichev use theorem \ref{theorem 5 intro} to prove a result about the convergence of the Harmonic extension of a potential which extends a classical result by A. Nagel, W. Rudin and J.H. Shapiro, see \cite{NRS}.\\
Let $1<p<+\infty$, let $0<\alpha\leq n$. Let $f\in L^p(\mathbb R^n)$. Let us consider the Bessel potential of $f$
\begin{equation}
B_{\alpha}*f(x):= \mathfrak F^{-1}\left((1+4  \pi^2 \lvert\xi\rvert^2 )^{-\frac{\alpha}{2}} \mathfrak F(f)(\xi)\right)(x),
\end{equation}
where $\mathfrak F(f)$ denotes the Fourier transform of $f$, and let us consider the harmonic extension of the Bessel potential of $f$, which is defined over $\mathbb R^{n+1}_+:= \mathbb R^n \times (0,+\infty)$ by
\begin{equation*}
PI(B_{\alpha}*f):\mathbb R^{n+1}_+\longrightarrow \mathbb R
\end{equation*}
\begin{equation*}
PI(B_{\alpha}*f)(x_0,y_0):=\int_{\mathbb R^n} \frac{\Gamma(\frac{n+1}{2})}{\pi^{\frac{n+1}{2}}} \frac{y_0}{\left(y_0^2+|x_0-x|^2\right)^{\frac{n+1}{2}}}B_{\alpha}*f(x)dx.
\end{equation*}
Nagel, Rudin and Shapiro prove the following result about the boundary behaviour of $PI(B_s*f)$.
\begin{Teo}[Convergence of $PI(B_s*f)$ at the boundary $\mathbb R^n \times \{0\}$, see \cite{NRS}, Theorem 5.5.]\label{Nagel Rudin Shapiro theorem in intro}
Let $1<p<+\infty$, $\frac 1 p + \frac{1}{p'}=1$, let $\frac{1}{p'}\leq s <1$. Let $f \in L^p(\mathbb R^n)$. Then, except eventually on a set of $x_0$'s of measure 0, $PI(B_{(1-s)n}*f)(x,y)$ converges to $B_{(1-s)n}*f(x_0)$ as $(x,y) \rightarrow (x_0,0)$ within the regions defined by:
\begin{enumerate}
\item \[y>c | x-x_0|^\frac{1}{1-p+s\cdot p} \quad \text{if } \frac{1}{p'}<s<1,\]
\item \[y>\exp\left(-c | x-x_0|^{-n(p'-1)}\right) \quad \text{if } s=\frac{1}{p'}.\]
\end{enumerate}
\end{Teo}
Aikawa and Borichev prove the following result, which extends  Theorem \ref{Nagel Rudin Shapiro theorem in intro} to the case of potentials $K*f$ associated to the radial kernels $K$ defined in (\ref{Aikawa Borichev kernels definition intro}).
\begin{Teo}[Tangential convergence of $PI(K*f)$ at the boundary $\mathbb R^n \times \{0\}$, see \cite{AB}, Theorem 11]\label{Aikawa Borichev theorem 11}
Let $1<p<+\infty$, $\frac 1 p + \frac{1}{p'}=1$, let $\frac{1}{p'}\leq s<1$. Let $f \in L^p(\mathbb R^n)$. Let $C_{K,p}$ be the capacity associated to the kernel $K$ defined by (\ref{Aikawa Borichev kernels definition intro}). Then, for almost all $x \in \mathbb R^n$ there exists a region $\Omega_{K,p,x} \subseteq \mathbb R^{n+1}_+$ such that $x$ belongs to the boundary of $\Omega_{K,p,x}$, and such that
\begin{equation}
\lim_{\underset{P \in_{\Omega_{K,p,x}}}{P \rightarrow(x,0)}}PI(K*f)(P)=K*f(x),
\end{equation}
and the region $\Omega_{K,p,x}$ has an order of contact of tagential type at the boundary point $(x,0) \in \partial \mathbb R^{n+1}_+$.
\end{Teo}
If we consider the case of the Bessel potential, where $K=B_{(1-s)n}$, we observe that (see \cite{ARSW1}, Corollary 28) it follows that we can choose a region $\Omega_{B_{(1-s)n},p,x}$ such that:
\begin{enumerate}
\item Case $\frac{1}{p'}<s<1$: the region $\Omega_{B_{(1-s)n},p,x}$  has an order of contact of polynomial type at the boundary point $(x,0) \in \partial \mathbb R^{n+1}_+$.
\item Case $s=\frac{1}{p'}$: the region $\Omega_{B_{(1-s)n},p,x}$  has an order of contact of exponential type at the boundary point $(x,0) \in \partial \mathbb R^{n+1}_+$.
\end{enumerate}
So Theorem \ref{Aikawa Borichev theorem 11} is an extension of Theorem \ref{Nagel Rudin Shapiro theorem in intro}.\\

In our work we extend Theorem \ref{Aikawa Borichev theorem 11} to the setting of compact Ahlfors-regular spaces for Riesz potentials. We consider a compact Ahlfors $Q$-regular space $(X,d,m)$, $Q>0$, and we define a notion of ``harmonic-type" extension in the following way:
\begin{Def}[Poisson Integral in $X \times (0,+\infty)$]
Let $f \in L^p(X)$, let $(x,y) \in X \times(0,+\infty)$. We define the Poisson Integral
\begin{equation}
PI(f)(x,y):= \int_{X} C(x,y) \cdot \frac{1}{y^Q} \sum_{k=0}^{+\infty} \frac{\chi_{B_d(x,2^k y)}(z)}{2^{(Q+1)k}}f(z)dm(z).
\end{equation}
Here $C(x,y)$ is the constant that normalizes the Poisson kernel, and it can be proved that
\begin{equation*}
C(x,y)\approx 1.
\end{equation*}
\end{Def}
In subsections 3.2 and 3.3 we show that this notion of harmonic extension shares some of the properties of the classical harmonic extension in $\mathbb R^{n+1}_+$, including a Harnack-type inequality.\\
With this notion of harmonic extension we prove the final result in our work (see Theorem \ref{Theorem 11} and Observation \ref{final observation}).

\begin{Teo}[Tangential convergence of $PI(K_{X,s}*f)$ at the boundary $X\times \{0\}$]\label{tangential convergence Ahlfors regular intro}
Let $(X,m,d)$ be a compact Ahlfors $Q$-regular space, $Q>0$. Let $1<p<+\infty$, $\frac 1 p + \frac{1}{p'}=1$, let $\frac{1}{p'}\leq s <1$. Let $K_{X,s}$ be the Riesz kernel defined by (\ref{Riesz kernel Ahlfors regular intro definition}). Let $f \in L^p(\mathbb R^n)$. Then, except eventually on a set of $x_0$'s of measure 0, $PI(K_{X,s}*f)(x,y)$ converges to $K_{X,s}*f(x_0)$ as $(x,y) \rightarrow (x_0,0)$ within the regions defined by:
\begin{enumerate}
\item \[y>c \cdot d( x,x_0)^\frac{1}{1-p+s\cdot p} \quad \text{if } \frac{1}{p'}<s<1,\]
\item \[y>\exp\left(-c \cdot d( x,x_0)^{-Q}\right) \quad \text{if } s=\frac{1}{p'}.\]
\end{enumerate}
\end{Teo}

This theorem generalizes the result by Nagel, Rudin and Shapiro to the setting of compact Ahlfors-regular spaces, and constitutes a starting point for a potential theory on Ahlfors-regular spaces analogous to the classical potential theory on $\mathbb R^n$. We believe that Theorem \ref{tangential convergence Ahlfors regular intro} can be generalized to non-compact Ahlfors-regular spaces.\\

Our work is heavily inspired by \cite{AB}. We reference \cite{AE} and \cite{AH} for the general notions and facts about the theory of the potential. We reference \cite{ARSW1} for the theory of the potential on tree boundaries.\\

\textbf{Keywords}--- Ahlfors-regular metric spaces, tree boundaries, potential theory, non-linear capacity, quasi-additivity, harmonic extension, boundary behaviour.\\
\textbf{Mathematical Subject Classification 2020}--- 31E05, 31C45, 31B05, 31B25.

\section*{Notations} \addcontentsline{toc}{section}{Notations}
Let $a,b \in \mathbb R$. We write $a \lesssim b$ (respectively $a \gtrsim b$) if and only if there exists a constant $0<C<+\infty$ such that $a \leq C \cdot b$ (respectively $a \geq C \cdot b$). Here the constant $C$ does not depend on any of the parameters of the problem.\\
We write $a \lesssim_{(p_1,p_2,\dots,p_n)} b$ (respectively $a \gtrsim_{(p_1,p_2,\dots,p_n)} b$) if and only if there exists a constant $C=C(p_1,p_2,\dots,p_n)$, $0<C<+\infty$, such that $a \leq C \cdot b$ (respectively $a \geq C \cdot b$). Here the constant $C$ depends on the paramenters $p_1,p_2,\dots,p_n$. \\
We write $a \approx b$ if and only if both $a\lesssim b$ and $a \gtrsim b$ hold.\\
We write $a \approx_{(p_1,p_2,\dots,p_n)} b$ if and only if both $a\lesssim_{(p_1,p_2,\dots,p_n)} b$ and $a \gtrsim_{(p_1,p_2,\dots,p_n)} b$ hold.\\
Let $(X,d)$ be a metric space. Let $x \in X$, $r\geq0$. We denote by $B_d(x,r)$ the metric ball of radius $r$ and center $x$, i.e.
\begin{equation}
B_d(x,r):=\{y \in X \; | \; d(x,y)<r\}.
\end{equation}
\section{Quasi-additivity on tree boundaries}
In this section we prove a quasi-additivity formula for capacities associated to radial kernels in the setting of the tree boundaries.\\

\subsection{Setting of the problem}
Let $T$ be a tree. Suppose every node in $T$ has at least 2 children. In this section $X:= \partial T$ will denote the boundary of $T$. $X$ is a metric space, where the metric on $X$ is given by the natural (ultrametric) distance
\[\rho(x,y)=\delta^{\text{count}(x \wedge y,o)},\]
where $\delta<1$ is a fixed constant, and $\text{count}(x \wedge y,o)$ is the number of edges shared by the geodesics connecting $x$ and $y$ to the origin $o$.\\
Let $m$ be a $\sigma$-finite Borel measure on $X$.\\
Let $\mathcal D$ denote the set
\begin{equation}
\mathcal D:=\{0\} \cup \{\delta^{n} \; | \; n \in \mathbb N\}.
\end{equation}
Let $K:\mathcal D \rightarrow \mathbb R^+$ be a function. Suppose $K$ is lower semi-continuos in 0, i.e.
\begin{equation}
\liminf_{n\rightarrow +\infty} K(\delta^{n}) \geq K(0).
\end{equation}
We define, with a small abuse of notation, the kernel $K(x,y):=K(\rho(x,y))$ for $x,y \in X$. It follows that the function
\[x \longmapsto K(x,y_0)\]
is lower semi-continuous for every choice of $y_0 \in X$.\\
Suppose the kernel $K$ satisfies the following conditions:
\begin{equation}\label{function K integrability}
\sup_{x \in X}\int_{X} K(x,y)dm(y)<+\infty, \quad \quad  \sup_{y \in X}\int_{X} K(x,y)dm(x)<+\infty.
\end{equation}
Let us denote
\begin{equation}
\|K\|_1:= \max\bigg\{\sup_{x \in X}\int_{X} K(x,y)dm(y),  \sup_{y \in X}\int_{X} K(x,y)dm(x)\bigg \}<+\infty.
\end{equation}
\begin{Def}
Let $1<p<+\infty$. The capacity of a compact set $E\subseteq X$ is
\begin{equation}
C_{K,p}(E):=\inf \bigg\{ \|f\|^p_{L^p(X,m)} \; \bigg | \; K * f(x)\geq 1 \; \; \; \forall x \in E\bigg\},
\end{equation}
where
\[ K * f(x) := \int_X K(x,y) f(y) dm(y).\]
\end{Def}
\begin{Def}
Let $x \in X$, $r>0$. We define, when it exists, the radius
\begin{equation}
\eta_p(x,r):= \inf \bigg \{ \delta^{n-\frac 1 2} \in \mathbb R \; \bigg | \; n \in \mathbb N, \; m(B_{\rho}(x,\delta^{n-\frac 1 2}))\geq C_{K,p}(B_{\rho}(x,r))\bigg \}.
\end{equation}
We also define
\begin{equation}\label{bigger radius in tree boundary definition}
\eta^*_p(x,r):= \max\{r,\eta_p(x,r)\}.
\end{equation}
\end{Def}
It follows that $B_{\rho}(x,r)\subseteq B_{\rho}(x,\eta^*_p(x,r))$.
\begin{Oss}\label{observation existence eta star}
The radius $\eta_p(x,r)$ does not exist when $X$ is compact for $x \in X$ and $r>0$ such that
\begin{equation}
C_{K,p}(B_{\rho}(x,r))>m(X).
\end{equation}
However, all the propositions and theorems using $\eta^*_p$ can be proved by separating the cases where $\eta^*_p$ is not defined, using other properties, like the compactness of $X$. This follows from the properties of the Riesz capacity of a ball in an Ahlfors-regular space (see Proposition \ref{capacity of ball general formula proposition}). We will always assume that the radius $\eta^*_p$ exists in the following proofs. 
\end{Oss} 
For the proof of the following theorem we require a Young's inequality for tree boundaries, which is stated in the following lemma.
\begin{Lem}[Young's inequality for tree boundaries]\label{Young inequality lemma}
  Let $K = K (x, y)$ be a kernel on a metric measure space $(X,\rho,m)$, and $1 \leq p
  \leq +\infty$. Then,
  \[\| K * f \|_{L^{p}(X,m)}:= \left[ \int_X \left( \int_X K (x, y) f (y) dm(y) \right)^p dm(x) \right]^{\frac 1 p}
     \leq \|K\|_1 \| f\|_{L^{p}(X,m)} , \]
where
\[\|K\|_1:=\max \left\{ \sup_{x\in X} \int_X K (x, y) dm(y), \sup_{y\in X} \int_X K (x, y) dm(x)
     \right\}.\]
\end{Lem}

\proof
  We prove this result by interpolation. For $p = +\infty$ and $f \geq 0$, we have
  \begin{eqnarray*}
    \int_X K (x, y) f (y) dm(y) & \leq & \int_X K (x, y) dm(y) \| f\|_{L^{\infty}(X,m)}\leq\\
    & \leq & \left( \sup_{x\in X} \int_X K (x, y) d m(y) \right) \| f\|_{L^{\infty}(X,m)},
  \end{eqnarray*}
  hence 
\[\| K * f \|_{L^{\infty}(X,m)} \leq \left( \sup_{x\in X} \int_X K (x, y) dm(y) \right) \| f \|_{L^{\infty}(X,m)}.\]
For $p = 1$ we have
  \begin{eqnarray*}
    \int_X\int_X K (x, y) f (y) dm(y) dm(x) & = & \int_X K (x, y) dm(x) \int_X f (y) dm(y)\leq\\
    & \leq & \left( \sup_{y\in X} \int_X K (x, y) dm(x) \right) \int_X f (y) dm(y),
  \end{eqnarray*}
hence
\[\| K * f \|_{L^{1}(X,m)} \leq \left( \sup_{y\in X} \int K (x, y) dm(x)
  \right) \| f \|_{L^{1}(X,m)}.\]
  The result follows from Riesz-Thorin interpolation theorem.
\endproof
\subsection{Quasi additivity for tree boundaries}
We are now going to prove the first result in this work.
\begin{Teo}[Quasi-additivity for tree boundaries]\label{Quasi additivity for boundary of tree theorem}
Let $ J$ be a countable (or finite) set of indices. Let $\{ B_{\rho}(x_j, r_j)\}_{j \in J}$ be a family of metric balls in $X$ such that $\eta_p(x_j,r_j)$ exists for all $j \in J$. Suppose $E \subseteq X$ is a compact subset of $\bigcup_{j\in J}  B_{\rho}(x_j,r_j)$. Suppose $\{ B_{\rho}(x_j, \eta^*_p(x_j,r_j))\}_{j \in J}$ is disjoint.\\
Then
\begin{equation}\label{Thesis}
C_{K,p}(E) \leq \sum_{j \in J} C_{K,p}(E\cap B_{\rho}(x_j,r_j)) \leq A \cdot C_{K,p}(E),
\end{equation}
where $A=A(X,K,p)$, $1<A<+\infty$, is a constant depending only on $X$, $K$ and $p$.
\end{Teo}
For the constant $A$ in (\ref{Thesis}) we get the value $A=[ (2^{p'-1}+1) \|K\|_1^{p'} +2^{{p'}-1}]^\frac{1}{{p'}-1}$, where $\frac 1 p + \frac{1}{p'}=1$, which can be reformulated to
\begin{equation}
A\approx_{(p)} \max\{1, \|K\|_1^p\}.
\end{equation}
For the proof of Theorem \ref{Quasi additivity for boundary of tree theorem} we recall the dual definition of capacity (see \cite{AH}).
\begin{Teo}[Dual definition of capacity]
Let $\frac{1}{p}+\frac{1}{p'}=1$. Then
\begin{equation}
C_{K,p}(E)=\sup \bigg\{\|\mu\|^p \; \bigg | \; \mu \text{ is concentrated on }E; \; \|K*\mu\|_{L^{p'}(X,m)}\leq 1\bigg\},
\end{equation}
where
\[\|\mu\| := \int_X d\mu(x); \quad K*\mu(y):= \int_X K(x,y)d\mu(x)\]
\end{Teo}
\proof [Proof of Theorem \ref{Quasi additivity for boundary of tree theorem}]
If $C_{K,p}(E)=0$ the proof is trivial by the monotonicity of the capacity.\\
Suppose $C_{K,p}(E)>0$. Let $E_j:= E \cap B_{\rho}(x_j,r_j)$. Without loss of generality, we may assume $ C_{K,p}(E_j)>0$ for all $j$. Indeed, let $J_0:=\{j \in J \; | \; C_{K,p}(E_j)>0\}$. It follows that $C_{K,p}(E_j)=0$ for all $j \in J \backslash J_0$, so
\[\sum_{j \in J} C_{K,p}(E\cap B_{\rho}(x_j,r_j))=\sum_{j \in J_0} C_{K,p}(E\cap B_{\rho}(x_j,r_j)),\]
hence we will assume $ C_{K,p}(E_j)>0$ for all $j \in J$.\\
To prove the thesis (\ref{Thesis}) it is sufficient to prove that
\begin{equation}
\sum_{j \in J} C_{K,p}(E\cap B_{\rho}(x_j,r_j))\leq A  \cdot C_{K,p}(E).
\end{equation}
For simplicity we will assume $J \subseteq \mathbb N$. Let $E_j := E \cap B_{\rho}(x_j,r_j)$ for $j \in J$. Let  $\frac{1}{p}+\frac{1}{p'}=1$. By the dual definition of capacity for every $j \in J$ there exists an equilibrium measure $\mu^{E_j}$ such that
\begin{equation}\label{epsilon optimality}
\begin{cases}
\mu^{E_j} \text{ is concentrated on } E_j,\\
\|K*\mu^{E_j}\|_{L^{p'}(X,m)}=1,\\
\| \mu^{E_j}\|^p =C_{k,p}(E_j).
\end{cases}
\end{equation}
Let $\mu^*_j:=C_{k,p}(E_j)^{\frac{1}{p'}} \mu^{E_j}$, and $\mu^*:= \sum_j \mu^*_j$.\\
We get
\begin{equation}\label{star first condition}
\|\mu^*_j\|^p=C_{k,p}(E_j)^{\frac{p}{p'}} \|\mu^{E_j}\|^p.
\end{equation}
From (\ref{epsilon optimality}) and  (\ref{star first condition})  we get
\begin{equation}\label{third epsilon optimality}
\| \mu^*_j\|^p= C_{k,p}(E_j)^{1+ \frac{p}{p'}},
\end{equation}
which reformulates to
\begin{equation} \label{fourth epsilon optimality}
\| \mu^*_j\| =C_{k,p}(E_j).
\end{equation}
Now we are going to prove that
\begin{equation}
\|K*\mu^*\|_{L^{p'}(X,m)}^{p'} \leq C \cdot \underset {j \in J}{\sum} C_{K,p}(E_j),
\end{equation}
where $C=C(X,K,p)$ is a constant depending only on $X$, $K$ and $p$.\\
Let $j \in J$. We define the measure
\begin{equation}\label{definition of mu prime}
d\mu'_j(y):= \frac{\|\mu^*_j\|}{m(B_{\rho}(x_j,\eta^*_p(x_j,r_j)))} \chi_{B_{\rho}(x_j,\eta^*_p(x_j,r_j))}(y)d m(y); \quad \mu':=\underset{j \in J}{\sum}\mu'_j.
\end{equation}
By construction $\|\mu_j'\|=\|\mu_j^*\|$, so, using (\ref{fourth epsilon optimality}), we get
\begin{equation}\label{mu star bound}
\|\mu_j'\|\leq C_{K,p}(E_j).
\end{equation}
By the definition of the function $\eta^*_p$ we have
\begin{equation}\label{measure greater than capacity}
m(B_{\rho}(x_j,\eta^*_p(x_j,r_j)))\geq C_{K,p}(B_{\rho}(x_j,r_j))\geq C_{K,p}(B_{\rho}(x_j,r_j)\cap E)= C_{K,p}(E_j),
\end{equation}
so, using (\ref{definition of mu prime}) and (\ref{measure greater than capacity}), we get
\begin{align*}
d\mu'_j =&  \frac{\|\mu^*_j\|}{m(B_{\rho}(x_j,\eta^*_p(x_j,r_j)))} \chi_{B_{\rho}(x_j,\eta^*_p(x_j,r_j))}dm\leq\\
& \frac{C_{K,p}(E_j)}{C_{K,p}(E_j)} \chi_{B_{\rho}(x_j,\eta^*_p(x_j,r_j))}dm \leq\\
& \chi_{B_{\rho}(x_j,\eta^*_p(x_j,r_j))}dm.
\end{align*}
By construction $d\mu'$ has a density with respect to $dm$, so we may write $d\mu'=f \cdot dm$, where
\begin{equation}
f\leq \underset{j \in J}{\sum}\chi_{B_{\rho}(x_j,\eta^*_p(x_j,r_j))},
\end{equation}
and we get
\begin{equation}\label{f properties by construction}
\|\mu'\|=\|f\|_{L^1(X,m)}.
\end{equation}
By hypothesis the sets $B_{\rho}(x_j,\eta^*_p(x_j,r_j))$ are disjoint, so we have
\begin{equation}\label{f lower than 1}
f\leq 1.
\end{equation}
Using (\ref{function K integrability}), (\ref{mu star bound}), (\ref{f properties by construction}), (\ref{f lower than 1}) and Lemma \ref{Young inequality lemma} we get
\begin{align*}
\|K*\mu'\|^{p'}_{L^{p'}(X,m)}=&\int_{X}\bigg[\int_{X} K(\rho(x,y))d\mathcal \mu'(y)\bigg]^{p'} dm(x)=\\
&\int_{X}\bigg[\int_{X} K(\rho(x,y))f(y)dm(y)\bigg]^{p'} dm(x)\leq \\
&\|K\|^{p'}_1 \|f\|^{p'}_{L^{p'}(X,m)}\leq\|K\|^{p'}_1 \|f\|^{{p'}-1}_{L^{\infty}(X,m)} \|f\|_{L^1(X,m)}\leq\\
& \|K\|^{p'}_1 \cdot 1^{{p'}-1} \cdot  \|f\|_{L^1(X,m)}= \|K\|^{p'}_1 \|\mu'\|= \|K\|^{p'}_1\sum_j  \|\mu_j'\|  \leq \\
&\|K\|^{p'}_1 \sum_j C_{K,p}(E_j).
\end{align*}
So we proved the following estimate for the $q$-norm of the potential of the measure $\mu'$:
\begin{equation}\label{K bound}
\|K*\mu'\|^{p'}_{L^{p'}(X,m)}\leq \|K\|^{p'}_1 \sum_j C_{K,p}(E_j),
\end{equation}
where $\|K\|^{p'}_1$ is a constant depending only on $X$, $K$ and $p$, and, by (\ref{function K integrability}), $\|K\|^{p'}_1<+\infty$.\\

We are now going to prove an estimate for $K*\mu^*(\tilde x)$, for $\tilde x \in X$.\\
Let $j \in J$. Let $\tilde x \in X$ be a point such that $\tilde x \not \in B_{\rho}(x_j,\eta_p^*(x_j,r_j))$. The measure $\mu'_j$ is concentrated on the set $B_{\rho}(x_j,\eta_p^*(x_j,r_j))$, so we get
\begin{align*}
K*\mu'_j(\tilde x)=& \int_X K(\tilde x,y) d\mu'_j(y)=\int_{B_{\rho}(x_j,\eta_p^*(x_j,r_j))} K(\tilde x,y) d\mu'_j(y)\geq\\
&\bigg(\min_{y \in B_{\rho}(x_j,\eta_p^*(x_j,r_j))}K(\tilde x,y)\bigg) \mu'_j(B_{\rho}(x_j,\eta_p^*(x_j,r_j)))=\\
&\bigg(\min_{y \in B_{\rho}(x_j,\eta_p^*(x_j,r_j))}K(\tilde x,y)\bigg) \|\mu^*_j\| \frac{m( B_{\rho}(x_j,\eta_p^*(x_j,r_j)))}{m( B_{\rho}(x_j,\eta_p^*(x_j,r_j)))}=\\
&\bigg(\min_{y \in B_{\rho}(x_j,\eta_p^*(x_j,r_j))}K(\tilde x,y)\bigg) \|\mu^*_j\|.
\end{align*}
Moreover, $\mu^*_j$ is concentrated on $ E_j\subseteq B_{\rho}(x_j,\eta_p^*(x_j,r_j))$, so we get
\begin{align*}
K*\mu^*_j(\tilde x)=& \int_X K(\tilde x,y) d\mu^*_j(y)=\int_{B_{\rho}(x_j,\eta_p^*(x_j,r_j))} K(\tilde x,y) d\mu^*_j(y)\leq\\
&\bigg(\max_{y \in B_{\rho}(x_j,\eta_p^*(x_j,r_j))}K(\tilde x,y)\bigg) \mu^*_j(B_{\rho}(x_j,\eta_p^*(x_j,r_j)))=\\
&\bigg(\max_{y \in B_{\rho}(x_j,\eta_p^*(x_j,r_j))}K(\tilde x,y)\bigg)  \|\mu^*_j\|.
\end{align*}
So we proved that
\begin{equation}\label{max and min condition}
K*\mu'_j(\tilde x)\geq\bigg(\min_{y \in B_{\rho}(x_j,\eta_p^*(x_j,r_j))}K(\tilde x,y)\bigg) \|\mu^*_j\|; \quad K*\mu^*_j(\tilde x)\leq \bigg(\max_{y \in B_{\rho}(x_j,\eta_p^*(x_j,r_j))}K(\tilde x,y)\bigg)  \|\mu^*_j\|.
\end{equation}
However, $K(\tilde x,y)=K(\rho(\tilde x,y))$, and $\tilde x \not \in B_{\rho}(x_j,\eta_p^*(x_j,r_j))$, so, by definition of $\rho$, if follows that
\begin{equation}
\rho(\tilde x,y_1)=\rho(\tilde x,y_2) \quad \text{for all } y_1,y_2 \in B_{\rho}(x_j,\eta_p^*(x_j,r_j)).
\end{equation}
So the function
\[y \mapsto \rho(\tilde x,y)\]
is constant when $y \in B_{\rho}(x_j,\eta_p^*(x_j,r_j))$, which proves that
\begin{equation}\label{min equals max}
\min_{y \in B_{\rho}(x_j,\eta_p^*(x_j,r_j))}K(\tilde x,y)=\max_{y \in B_{\rho}(x_j,\eta_p^*(x_j,r_j))}K(\tilde x,y).
\end{equation}
So, using (\ref{max and min condition}) and  (\ref{min equals max}) , we get
\begin{equation}\label{K*mu star j estimate}
 K*\mu^*_j(\tilde x) \leq K*\mu'_j(\tilde x) \quad \text{for } \tilde x \not \in B_{\rho}(x_j,\eta_p^*(x_j,r_j)).
\end{equation}
If $\tilde x \not \in \bigcup_j B_{\rho}(x_j,\eta_p^*(x_j,r_j))$ then, by applying (\ref{K*mu star j estimate}) for all $j \in J$, we get
\begin{equation}\label{K*mu star estimate outside}
 K*\mu^*(\tilde x) \leq K*\mu'(\tilde x).
\end{equation}
If $\tilde x \in B_{\rho}(x_{j_0},\eta_p^*(x_{j_0},r_{j_0}))$ for some $j_0 \in J$, then, by the disjointness of the family $\{ B_{\rho}(x_j,\eta_p^*(x_j,r_j))\}_{j \in J}$, we have 
\[\tilde x \not \in \bigcup_{j\neq j_0} B_{\rho}(x_j,\eta_p^*(x_j,r_j)),\]
so we get
\begin{align*}
K*\mu^*(\tilde x)=&K*\mu_{j_0}^*(\tilde x)+\sum_{j \neq j_0} K*\mu^*_j(\tilde x)\leq\\
&K*\mu^*_{j_0}(\tilde x)+\sum_{j \neq j_0} K * \mu'_j(\tilde x)\leq\\
&K*\mu^*_{j_0}(\tilde x)+\sum_{j} K * \mu'_j(\tilde x)=\\
&K*\mu^*_{j_0}(\tilde x)+K * \mu'(\tilde x).
\end{align*}
So we proved that, if $\tilde x \in  B_{\rho}(x_{j_0},\eta_p^*(x_{j_0},r_{j_0}))$ for some $j_0 \in J$, then
\begin{equation}\label{K*mu star estimate inside}
K*\mu^*(\tilde x) \leq K*\mu^*_{j_0}(\tilde x)+K * \mu'(\tilde x).
\end{equation}
\ \\

Now we are going to prove that
\begin{equation}
\|K*\mu^*\|^{p'}_{L^{p'}(X,m)} \leq C \cdot \sum_j C_{k,p}(E_j),
\end{equation}
where $C=C(X,K,p)$ is a constant depending only on $X$, $K$ and $p$.\\
We use (\ref{K*mu star estimate outside}), (\ref{K*mu star estimate inside}) and the disjointness of the sets $\{B_{\rho}(x_j,\eta_p^*(x_j,r_j))\}_j$ to estimate
\begin{align*}
\|K*\mu^*\|^{p'}_{L^{p'}(X, m)} =& \int_{X \backslash \bigcup_j  B_{\rho}(x_j,\eta_p^*(x_j,r_j))}(K*\mu^*(x))^{p'} d m(x) +\\
&\int_{\bigcup_j  B_{\rho}(x_j,\eta_p^*(x_j,r_j))}(K*\mu^*(x))^{p'} d m(x)\leq\\
& \int_{X \backslash \bigcup_j  B_{\rho}(x_j,\eta_p^*(x_j,r_j))}(K*\mu'(x))^{p'} d m(x) +\\
&\sum_{j} \int_{B_{\rho}(x_j,\eta_p^*(x_j,r_j))}(K*\mu^*(x))^{p'} d m(x)\leq\\
&\|K*\mu'\|^{p'}_{L^{p'}(X, m)} +\sum_{j} \int_{B_{\rho}(x_j,\eta_p^*(x_j,r_j))}(K*\mu^*_j(x)+K*\mu'(x))^{p'} d m(x).
\end{align*}
We are going to apply Jensen's inequality for finite sums to the last inequality. Let $n\geq 2$ be a natural number, let ${p'}>1$, let $a_i\geq 0$ for $1\leq i \leq n$. By Jensen's inequality we have
\begin{equation}\label{Jensen for sums}
\bigg(\sum_{i=1}^n a_i \bigg)^{p'} \leq n^{{p'}-1} \sum_{i=1}^n a_i^{p'}.
\end{equation}
For each fixed $j \in J$ we apply (\ref{Jensen for sums}) to the last estimate for the ${p'}$-norm of $K*\mu^*$, where $n=2$, $a_1=K*\mu^*_j(x)$ and $a_2=K*\mu'(x)$, and we get
\begin{equation}
\|K*\mu^*\|^{p'}_{L^{p'}(X, m)} \leq \|K*\mu'\|^{p'}_{L^{p'}(X, m)} +2^{{p'}-1}\sum_{j} \int_{B_{\rho}(x_j,\eta_p^*(x_j,r_j))}(K*\mu^*_j(x)^{p'}+K*\mu'(x)^{p'}) d m(x).
\end{equation}
So we get the estimate
\begin{equation}\label{last estimate for K*mu*}
\|K*\mu^*\|^{p'}_{L^{p'}(X, m)} \leq \|K*\mu'\|^{p'}_{L^{p'}(X, m)} +2^{{p'}-1}\sum_{j}\|K*\mu^*_j\|^{p'}_{L^{p'}(X,m)}+2^{{p'}-1}\|K*\mu'\|^{p'}_{L^{p'}(X,m)}.
\end{equation}
We observe that, by definition of $\mu^*$, we have
\begin{equation}
\|K*\mu_j^*\|^{p'}_{L^{p'}(X, m)}=\|K*(C_{K,p}(E_j)^{\frac{1}{p'}}\mu_j)\|^{p'}_{L^{p'}(X, m)}=C_{K,p}(E_j)\|K*\mu_j\|^{p'}_{L^{p'}(X, m)},
\end{equation}
but $\|K*\mu_j\|_{L^{p'}(X, m)}=1$ by (\ref{epsilon optimality}), so we get
\begin{equation}\label{norm of K*mu j is capacity}
\|K*\mu_j^*\|^{p'}_{L^{p'}(X, m)}=C_{K,p}(E_j).
\end{equation}
So we use (\ref{K bound}) and (\ref{norm of K*mu j is capacity}) and the estimate (\ref{last estimate for K*mu*}) to get
\begin{equation}
\|K*\mu^*\|^{p'}_{L^{p'}(X, m)} \leq  \|K\|^{p'}_1 \sum_j C_{K,p}(E_j) +2^{{p'}-1}\sum_{j} C_{K,p}(E_j)+ 2^{{p'}-1}\|K\|^{p'}_1 \sum_j C_{K,p}(E_j).
\end{equation}
So we proved that
\begin{equation}\label{last estimate norm K*mu^*}
\|K*\mu^*\|^{p'}_{L^{p'}(X, m)} \leq\bigg[ (2^{{p'}-1}+1) \|K\|^{p'}_1 +2^{{p'}-1}\bigg]\sum_j C_{K,p}(E_j),
\end{equation}
where $C=[ (2^{{p'}-1}+1) \|K\|^{p'}_1 +2^{{p'}-1}]<+\infty$ is a constant depending only on $X$, $K$ and $p$.\\
\ \\

Now we are going to finish the proof by defining a proper normalized measure.\\
We define the measure
\begin{equation}
\mathring{\mu}:= \bigg[C \cdot \sum_j C_{K,p}(E_j)\bigg]^{-\frac{1}{p'}} \mu^*.
\end{equation}
From (\ref{last estimate norm K*mu^*}) we get
\begin{equation}\label{final estimate norm K*mu^*}
\|K*\mu^*\|_{L^{p'}(X, m)} \leq \bigg[C \cdot \sum_j C_{K,p}(E_j)\bigg]^{\frac{1}{p'}}.
\end{equation}
By construction $\mathring{\mu}$ is concentrated on $E$, and, using (\ref{final estimate norm K*mu^*}), we get
\begin{align*}
\|K*\mathring{\mu}\|_{L^{p'}(X, m)}=&\bigg\| K*\bigg[\big[C \cdot \sum_j C_{K,p}(E_j)\big]^{-\frac{1}{p'}} \mu^*\bigg]\bigg\|_{L^{p'}(X, m)}=\\
&\big[C \cdot \sum_j C_{K,p}(E_j)\big]^{-\frac{1}{p'}}\|K*\mu^*\|_{L^{p'}(X, m)}\leq\\
&\big[C \cdot \sum_j C_{K,p}(E_j)\big]^{-\frac{1}{p'}}\big[C \cdot \sum_j C_{K,p}(E_j)\big]^{\frac{1}{p'}}=1.
\end{align*}
So we proved that the measure $\mathring{\mu}$ is a test measure for the dual definition of capacity, so we get
\begin{equation}\label{Capacity of E estimate start}
C_{K,p}(E)\geq \|\mathring{\mu}\|^p.
\end{equation}
By computation, using (\ref{fourth epsilon optimality}), we get
\begin{align*}
\|\mathring{\mu}\|=&\bigg[C \cdot \sum_j C_{K,p}(E_j)\bigg]^{-\frac{1}{p'}}\| \mu^*\|=\\
&\bigg[C \cdot \sum_j C_{K,p}(E_j)\bigg]^{-\frac{1}{p'}}\sum_j \| \mu_j^*\|=\\
&\bigg[C \cdot \sum_j C_{K,p}(E_j)\bigg]^{-\frac{1}{p'}}\sum_j C_{K,p}(E_j).
\end{align*}
So, using the last equation and (\ref{Capacity of E estimate start}), we get
\begin{align*}
C_{K,p}(E)\geq C^{-\frac{p}{p'}} \bigg[ \sum_j C_{K,p}(E_j)^{1-\frac{1}{p'}}\bigg]^p,
\end{align*}
and, since $1-\frac{1}{p'} = \frac 1 p$, we get
\begin{equation}
\sum_j C_{K,p}(E_j)\leq C^{\frac{p}{p'}} \cdot C_{K,p}(E).
\end{equation}
So, since $\frac{p}{p'}=\frac{1}{{p'}-1}$, we proved that
\begin{equation}
\sum_j C_{K,p}(E_j)\leq \bigg[ (2^{{p'}-1}+1) \|K\|^{p'}_1 +2^{{p'}-1}\bigg]^\frac{1}{{p'}-1} C_{K,p}(E),
\end{equation}
where the constant $A=[ (2^{{p'}-1}+1) \|K\|^{p'}_1 +2^{{p'}-1}]^\frac{1}{{p'}-1}<+\infty$ is a constant depending only on $X$, $K$ and $p$, ending the proof.
\endproof
\begin{Rem}
The previous theorem also holds (up to modifying the constant $A$) for a generic non-radial kernel $K=K(x,y)\geq 0$ such that:
\begin{itemize}
\item The function $x\mapsto K(x,y_0)$ is lower semi-continuous for any $y_0 \in X$.
\item The function $y\mapsto K(x_0,y)$ is measurable for any $x_0 \in X$.
\item The kernel is globally integrable, i.e.
\[\|K\|_1:= \max\bigg\{\sup_{x \in X}\int_{X} K(x,y)dm(y),  \sup_{y \in X}\int_{X} K(x,y)dm(x)\bigg \}<+\infty.\]
\item There exists a constant $C=C(X,K,p)$ such that at least one of the following conditions holds:
\begin{enumerate}
\item \[K*\mu^{E_j}(x)\leq C \cdot K*\left(\mu^{E_j}\right)'(x),\]
for all $x \in  B_{\rho}(x_j,\eta_p^*(x_j,r_j))$, for all $j \in J$, where $\mu^{E_j}$ is the equilibrium measure for the set $E_j:=E \cap B_{\rho}(x_j,\eta_p^*(x_j,r_j))$, and $\left(\mu^{E_j}\right)'$ is the measure defined by
\[d\left(\mu^{E_j}\right)'(y):= \frac{\|\mu^{E_j}\|}{m(B_{\rho}(x_j,\eta^*_p(x_j,r_j)))} \chi_{B_{\rho}(x_j,\eta^*_p(x_j,r_j))}(y)d m(y).\]
\item \[\sup_{y \in  B_{\rho}(x_j,\eta_p^*(x_j,r_j))} K(\tilde x,y) \leq C \cdot \inf_{y \in  B_{\rho}(x_j,\eta_p^*(x_j,r_j))} K(\tilde x,y),\]
for all $\tilde x \not \in B_{\rho}(x_j,\eta_p^*(x_j,r_j))$, for any $j \in J$.
\end{enumerate}
We observe that the condition $2$ entails the condition $1$.
\end{itemize}
\end{Rem}

\section{Ahlfors-regular spaces}
In this section we prove a quasi-additivity formula for the Riesz capacity in the setting of the Ahlfors-regular spaces.
\subsection{Setting of the problem}
\begin{Def}
Let $(X,d,m)$ be a compact metric measure space. Let $Q>0$. We say that $X$ is a $Q$-regular Ahlfors space if there exist constants $0<C_1\leq C_2$ such that
\begin{equation}\label{ahlfors regular condition for X}
C_1 r^Q \leq m(B(x,r)) \leq C_2 r^Q
\end{equation}
for all $x \in X$, for all $0<r<\text{diam}(X)$.\\
We will say that $X$ is an Ahlfors-regular space without mentioning the dimension $Q$ when that paremeter is not relevant to the discussion.
\end{Def}
Our work focuses on $Q$-regular Ahlfors spaces such that the measure $m$ is the Hausdorff measure of dimensional parameter $Q$. Hausdorff measures are known to be regular measures, so we will assume regularity of the measure $m$ in the following part of this work.\\
The following theorem and definitions allow us construct a tree structure starting from an Ahlfors-regular space (see \cite{ARSW1} for more details).
\begin{Teo}[Christ decomposition]\label{Christ decomposition}
Let $(X,d,m)$ be a compact $Q$-regular Ahlfors space. There exists a collection of qubes $\{ Q_{\alpha}^k \subseteq X \; | \; \alpha \in I_k, \, k \in \mathbb N\}$, where $I_k$ is a set of indices, and there exist constants $0<\delta<1$, $C_3>0$ and $C_4>0$ such that:
\begin{itemize}
\item[i)] $m\left(X\backslash\bigcup_{\alpha \in I_k} Q_{\alpha}^k\right)=0$ for all $k \in \mathbb N$.
\item[ii)] If $l\geq k$ then $\forall \alpha \in I_k$, $\forall \beta \in I_l$ we have either $Q_{\beta}^l \subseteq Q_{\alpha}^k$ or $Q_{\beta}^l \cap Q_{\alpha}^k = \emptyset$.
\item[iii)] For all $l \in \mathbb N$, for all $\beta \in I_l$ and for all $k<l$ there exists a unique $\alpha \in I_k$ such that $Q_{\beta}^l\subseteq Q_{\alpha}^k$.
\item[iv)] $\text{diam}(Q_{\alpha}^k) \leq C_4 \cdot \delta^k$.
\item[v)] For all $k \in \mathbb N$, for all $\alpha \in I_k$, there exist $z_{\alpha}^k \in Q_{\alpha}^k$ such that $B(z_{\alpha}^k,C_3 \cdot \delta^k) \subseteq Q_{\alpha}^k$.
\end{itemize}
Moreover, we may assume that $I_0=\{1\}$, and $Q_1^0=X$.
\end{Teo}
\begin{Def}
Let $(X,d,m)$ be an Ahlfors $Q$-regular space as above. Consider a Christ decomposition 
\begin{equation}
\{ Q_{\alpha}^k \subseteq X \; | \; \alpha \in I_k, \, k \in \mathbb N\},
\end{equation}
like the one in Theorem \ref{Christ decomposition}, such that $I_0=\{1\}$, and $Q_1^0=X$.\\
We define a tree structure $T \equiv (V(T),E(T))$, where the set of the vertices of $T$ is
\begin{equation}
V(T):=\{ Q_{\alpha}^k \subseteq X \; | \; \alpha \in I_k, \, k \in \mathbb N\},
\end{equation}
and we define the set of edges $E(T)$ in the following way: for all $Q_{\alpha}^k \in V(T)$ and for all $Q_{\beta}^l \in V(T)$ then $(Q_{\alpha}^k, Q_{\beta}^l) \in E(T)$ if and only if $k = l-1$ and $\alpha$ is the unique index in $I_k$ such that $Q_{\beta}^l \subseteq Q_{\alpha}^k$ defined by property iii) in Theorem \ref{Christ decomposition}.\\
The structure $T \equiv(V(T),E(T))$ is a tree structure such that $o:=Q_1^0=X$ is the root of the tree $T$.
\end{Def}
\begin{Def}\label{boundary of the tree definition}
$(X,d,m)$ be a $Q$-regular Ahlfors space, and let $T\equiv (V(T),E(T))$ be the tree previously defined. We define the boundary $\partial T$ of the tree $T$ as the set of half infinite geodesics starting at the origin $o$, i.e.
\begin{equation}
\partial T:=\bigg\{\left(Q_{\alpha_0}^0,Q_{\alpha_1}^1,\dots,Q_{\alpha_{j-1}}^{j-1},Q_{\alpha_j}^j,Q_{\alpha_{j+1}}^{j+1}, \dots \right)\; \bigg| \; Q_{\alpha_0}^0=o, \; (Q_{\alpha_k}^k,Q_{\alpha_{k+1}}^{k+1}) \in E(T) \; \forall k \in N \bigg\}.
\end{equation}
We define the map
\begin{equation}\label{Lambda map definition}
\Lambda: \partial T \longrightarrow X
\end{equation}
\begin{equation*}
\Lambda((o,Q_{\alpha_1}^1,Q_{\alpha_2}^2,\dots))=\bigcap_{n \in N} \left( \overline{Q_{\alpha_n}^n}\right).
\end{equation*}
The map $\Lambda$ identifies $X$ and $\partial T$ (see \cite{ARSW1} for more details).\\
Let $x=(o,Q_{\alpha_1}^1,Q_{\alpha_2}^2,\dots) \in \partial T$ and $y=(o,Q_{\beta_1}^1,Q_{\beta_2}^2,\dots) \in \partial T$. We define $x \wedge y \in V(T) \cup \partial T$ in the following way: if $x=y$ then $x \wedge y := x=y$, if $x \neq y$ then $x\wedge y := Q_{\gamma}^j$, where $j= \max\{k \in \mathbb N \; | \; \alpha_l=\beta_l \; \;\text{for all } l\leq k\}$, and $\gamma = \alpha_k=\beta_k$.\\
 Let us define the distance 
\begin{equation}\label{rho distance def}
\rho: \partial T \times \partial T \longrightarrow \mathbb R
\end{equation}
\begin{equation*}
\rho(x,y)=\delta^{\text{count}(x \wedge y,o)},
\end{equation*}
where $0<\delta<1$ is the constant defined in Theorem \ref{Christ decomposition}, and $\text{count}(x \wedge y,o)$ is the distance that counts how many edges of the geodesic that connects $x \wedge y$ and $o$ are in between $o$ and $x \wedge y$. In particular, if $x \wedge y=Q_{\alpha}^k \in V(T)$ then $\text{count}(x \wedge y,o)=k$, otherwise $x \wedge y=x=y \in \partial T$ and $\text{count}(x \wedge y,o)=+\infty$ and $\rho(x,y)=0$. We have $\text{diam}(\partial T)=1<+\infty$.
\end{Def}
Let ${\mathcal H}_{\rho}^Q$ denote the $Q$-dimensional Hausdorff measure on $\partial T$ with respect to the distance $\rho$. The space $\partial T$ endowed with the distance $\rho$ and the measure ${\mathcal H}_{\rho}^Q$ is a compact $Q$-regular Ahlfors space, so there exist constants $0<K_1<K_2$ such that
\begin{equation}\label{Ahlfors regular condition on partial T}
K_1  \cdot r^Q \leq \mathcal H_{\rho}^Q\left(B_{\rho}(x,r)\right) \leq K_2 \cdot r^Q
\end{equation}
for all $x \in \partial T$, for all $0<r<\text{diam}_{\rho}(\partial T)$. Here $\left(B_{\rho}(x,r)\right)$ denotes the metric ball of center $x$ and radius $r$ in $\partial T$ with respect to the metric $\rho$, and $\text{diam}_{\rho}$ denotes the diameter with respect to the metric $\rho$.
\subsection{Riesz potential}
In this subsection we define the Riesz capacity in our setting and we enunciate some properties we will use later in this work.
\begin{Def}[Riesz potential on $X$]
Let $(X,d,m)$ be a compact $Q$-regular Ahlfors space. Let $1<p<+\infty$ and $\frac 1 p+\frac{1}{p'}=1$. Let $\frac{1}{p'}<s<1$. We define the Riesz kernel
\begin{align}
K_{X,s}:X\times X &\longrightarrow \mathbb R\\
(x,y) &\longmapsto \frac{1}{d(x,y)^{Q \cdot s}}. \nonumber
\end{align}
Let $E \subseteq X$ be a compact set. We define the $L^p$ capacity of $E$ associated to the kernel $K_{X,s}$:
\begin{equation}
C_{K_{X,s},p}(E):= \inf \left\{\| f\|^p_{L^p(X,m)} \; \big| \; f \in L^p(X,m), \; Gf(x)\geq 1 \; \forall x \in E\right \},
\end{equation}
where $Gf$ denotes the potential of $f$, and it is defined by
\begin{equation}
Gf(x):= \int_X K_{X,s}(x,y)f(y)dm(y).
\end{equation}
\end{Def}
\begin{Def}[Riesz potential on $\partial T$]
Let $(X,d,m)$ be a compact $Q$-regular Ahlfors space. Let $\partial T$ be the boundary of the tree $T$ associated to $X$. Let $\partial T$ be endowed with the distance $\rho$ defined in \ref{rho distance def} and with the Hausdorff measure ${\mathcal H}_{\rho}^Q$. Let $1<p<+\infty$ and $\frac 1 p+\frac{1}{p'}=1$. Let $\frac{1}{p'}<s<1$. We define the Riesz kernel
\begin{align}
K_{\partial T,s}:\partial T\times \partial T &\longrightarrow \mathbb R\\
(x,y) &\longmapsto \frac{1}{\rho(x,y)^{Q \cdot s}}. \nonumber
\end{align}
Let $E \subseteq \partial T$ be a compact set. We define the $L^p$ capacity of $E$ associated to the kernel $K_{\partial T,s}$:
\begin{equation}
C_{K_{\partial T,s},p}(E):= \inf \left\{\| f\|^p_{L^p(\partial T, \rho)} \; \big| \; f \in L^p(\partial T, \rho), \; If(x)\geq 1 \; \forall x \in E\right \},
\end{equation}
where $If$ denotes the potential of $f$, and it is defined by
\begin{equation}
If(x):= \int_{\partial T} K_{\partial T,s}(x,y)f(y)d{\mathcal H
}_{\rho}^Q(y).
\end{equation}
\end{Def}
The following estimate for the capacity of a ball in an Ahlfors-regular space is a known result from the general theory of potential on Ahlfors-regular spaces, and it can be found in \cite{ARSW1}, Corollary 28. Unfortunately, the statement provided in \cite{ARSW1} is filled with typos, and it should be reformulated in the following way.
\begin{Prop}[Estimate for the capacity of a ball in an Ahlfors-regular space]\label{capacity of ball general formula proposition}
Let $(X,d,m)$ be a $Q$-regular Ahlfors space, let $p>1$, $\frac{1}{p}+\frac{1}{p'}=1$, let $\frac{1}{p'}\leq s <1$. Then there exist constants $0<\tilde C_1 <\tilde C_2$, $0<\tilde K_1 <\tilde K_2$, $r_0>0$ which depend only on $X$, $Q$, $p$ and $s$ such that for all $x \in X$, for all $r<r_0$ the following formulas hold:
\begin{itemize}
\item Case $\frac{1}{p'}<s<1$: 
\begin{equation}
\tilde C_1 \cdot  r^{Qp\left(s-\frac{1}{p'}\right)}\leq C_{K_{X,s},p}(B_d(x,r)) \leq \tilde C_2 \cdot  r^{Qp\left(s-\frac{1}{p'}\right)}.
\end{equation}
\item Case $s=\frac{1}{p'}$: 
\begin{equation}
\tilde K_1 \cdot \frac{1}{\log\left(\frac{1}{r}\right)}\leq C_{K_{X,s},p}(B_d(x,r)) \leq \tilde K_2 \cdot \frac{1}{\log\left(\frac{1}{r}\right)}.
\end{equation}
\end{itemize}
\end{Prop}
The following theorem is a known result from the theory of capacity on trees and Ahlfors-regular spaces  (see \cite{ARSW1}), and it will be used in the proof of the second main result of this work.
\begin{Teo}[Comparing the capacities on $X$ and $\partial T$]\label{comparing capacities theorem}
Let $(X,d,m)$ be a $Q$-regular Ahlfors space, let $(\partial T, \rho,{\mathcal H}_{\rho}^Q)$ be the $Q$-regular Ahlfors boundary of the associated tree $T$. Then there exist constants $A_1,A_2$ such that $0<A_1\leq A_2$, and such that for every closed set $F \subseteq \partial T$ and for every closed set $G \subseteq X$ we have
\begin{enumerate}
\item \begin{equation*}
A_1 \cdot C_{K_{\partial T,s},p}(F)\leq C_{K_{X,s},p}(\Lambda(F)) \leq A_2 \cdot C_{K_{\partial T,s},p}(F).
\end{equation*}
\item \begin{equation*}
A_1 \cdot C_{K_{X,s},p}(G)\leq C_{K_{\partial T,s},p}(\Lambda^{-1}(G)) \leq A_2 \cdot C_{K_{X,s},p}(G).
\end{equation*}
\end{enumerate}
\end{Teo}
\subsection{Quasi-additivity on compact Ahlfors-regular spaces}
The following theorem is the second result in this work, and it will be used later in the proof of the main results in this work.
\begin{Teo}[Quasi-additivity for Riesz capacity on compact Ahlfors-regular spaces] \label{Quasi additivity Ahlfors space theorem}
Let $(X,d,m)$ be a compact $Q$-regular Ahlfors space. Let $1<p<+\infty$ and $\frac 1 p + \frac{1}{p'}=1$. Let $\frac{1}{p'}\leq s<1$. For every $x \in X$ and $r>0$ let us define (when it exists) the radius
\begin{equation}
\eta_{X,p}(x,r):= \inf \left\{ R>0 \; | \; m(B_d(x,R))\geq C_{K_{X,s},p}(B_d(x,r))\right\},
\end{equation}
and let us define
\begin{equation}\label{eta star in X definition}
\eta^*_{X,p}(x,r):=\max \{ r,\eta_{X,p}(x,r)\}.
\end{equation}
Then there exists a constant $\Psi=\Psi(X,p,s)\geq 1$ such that for all $M\geq 1$ there exists a constant $1<\tilde A<+\infty$ such that, for any countable family $\{B_d(x_k,r_k)\}_{k\in \mathcal F}$ of balls in $X$ such that the family $\{B_d(x_k,\Psi \cdot \eta^*_{X,p}(x_k,M \cdot r_k))\}_{k\in \mathcal F}$ is disjoint, for any compact set $E \subseteq X$ such that $E=\bigcup_k E_k$ and $E_k \subseteq B_d(x_k,r_k) \; \forall k$, we have
\begin{equation}\label{thesis quasi additivity ahlfors}
\sum_{k \in \mathcal F} C_{K_{X,s},p}(E_k) \leq \tilde A \cdot C_{K_{X,s},p}(E).
\end{equation}
The constant $\tilde A$ depend only on the choice of the space $X$ and the of the parameters $p$, $s$ and $M$.
\end{Teo}
\proof
We begin the proof with the case $M>1$. Let $M>1$ be arbitrary. We can write $M$ in the following way
\begin{equation}\label{epsilon definition}
M=1+\epsilon>1, \quad \text{where } \epsilon >0.
\end{equation}
Let $\{B_d(x_k,r_k)\}_{k \in \mathcal F}$ be a family of balls in $X$. Let $E= \bigcup_{k \in \mathcal F}E_k$ be a compact set such that $E_k$ is a compact subset of $B_d(x_k,r_k)$ for all $k \in \mathcal F$. Consider a Christ decomposition
\begin{equation}
\{ Q_{\alpha}^k \subseteq X \; | \; \alpha \in I_k, \, k \in \mathbb N\} 
\end{equation}
like the one defined in Theorem \ref{Christ decomposition}. Consider le associated tree structure $T$ and the boundary of the tree $(\partial T, \rho,{\mathcal H}_{\rho}^Q)$ defined in definition \ref{boundary of the tree definition}.\\
Without loss of generality we may assume $r_k\leq C_3$ for all $k \in \mathcal F$. Consider a fixed $k \in \mathcal F$. We define
\begin{equation}\label{j(k) level definition}
j(k):=\max\{j \in \mathbb N \; | \; C_3 \cdot \delta ^j \geq r_k \},
\end{equation}
where $0<\delta<1$ and $C_3>0$ are the constants defined in Theorem \ref{Christ decomposition}.\\
By the statement v) in Theorem \ref{Christ decomposition} we have that for every $k \in \mathcal F$ for every $\alpha \in I_{j(k)}$ there exists $z_{\alpha}^{j(k)} \in Q_{\alpha}^{j(k)}$ such that $B(z_{\alpha}^{j(k)},C_3 \cdot \delta^{j(k)}) \subseteq Q_{\alpha}^{j(k)}$.\\
Moreover, by the statement ii) in Theorem \ref{Christ decomposition} we have 
\begin{equation}Q_{\alpha_1}^{j(k)} \cap  Q_{\alpha_2}^{j(k)}=\emptyset \quad \text{for all } \alpha_1,\alpha_2 \in I_{j(k)}\; \text{ such that }\alpha_1\neq \alpha_2.
\end{equation}
Let us consider the family $\mathcal G(k) \subseteq I_{j(k)}$ of indices such that the family of qubes
\begin{equation}
\left\{ Q_{\alpha}^{j(k)} \; | \; \alpha \in \mathcal G(k)\right\}
\end{equation}
is the minimal covering of qubes at the level $j(k)$ of the set $E_k$, i.e.
\begin{equation}\label{qubes covering definition}
E_k \subseteq \bigcup_{\alpha \in \mathcal G(k)} \overline{Q_{\alpha}^{j(k)}} \quad \text{and }\quad E_k \cap \overline{Q_{\alpha}^{j(k)}}\neq \emptyset \; \forall \alpha \in \mathcal G(k).
\end{equation}
Consider an arbitrary $\alpha \in \mathcal G(k)$. Let $w(k)\in \overline{Q_{\alpha}^{j(k)}}$ be arbitrary. By (\ref{qubes covering definition}) there exists $y(k) \in E_k \cap \overline{Q_{\alpha}^{j(k)}}$. By construction $E_k \subseteq B_d(x_k,r_k)$, so $d(x_k,y(k))\leq r_k$. Moreover, by statement iv) in Theorem \ref{Christ decomposition} we have $\text{diam}(\overline{Q_{\alpha}^{j(k)}})\leq C_4 \cdot \delta^{j(k)}$. By definition \ref{j(k) level definition} we have $\delta^{j(k)}\leq \delta^{-1} \cdot r_k$, so we get
\begin{equation}
\text{diam}(\overline{Q_{\alpha}^{j(k)}}) \leq C_4 \cdot \delta^{-1} \cdot r_k.
\end{equation}
So, by triangle inequality, we get
\begin{equation}
\overline{Q_{\alpha}^{j(k)}} \subseteq B_d(x_k,r_k(1+ C_4 \cdot \delta^{-1})).
\end{equation}
We are now going to repeat the previous steps, but instead of qubes at the level $j(k)$ we will use qubes at the level $j(k)+n$, where $n \in \mathbb N$ will be fixed later.\\
We define the family $\mathcal H(k,n) \subseteq I_{j(k)+n}$ such that
\begin{equation}\label{thinner cubes covering definition}
E_k \subseteq \bigcup_{\alpha \in \mathcal H(k,n)} \overline{Q_{\alpha}^{j(k)+n}} \quad \text{and }\quad E_k \cap \overline{Q_{\alpha}^{j(k)+n}}\neq \emptyset \; \forall \alpha \in \mathcal H(k,n).
\end{equation}
For all $\alpha \in \mathcal H(k,n)$ we have that there exists $z_{\alpha}^{j(k)+n} \in Q_{\alpha}^{j(k)+n}$ such that \begin{equation}\label{j(k)+n ball inside the qube exists}
B(z_{\alpha}^{j(k)+n},C_3 \cdot \delta^{j(k)+n}) \subseteq Q_{\alpha}^{j(k)+n},
\end{equation}
and such that
\begin{equation}\label{j(k)+n qubes are not intersected}
Q_{\alpha_1}^{j(k)+n} \cap  Q_{\alpha_2}^{j(k)+n}=\emptyset \quad \text{for all } \alpha_1,\alpha_2 \in I_{j(k)+n}\; \text{ such that }\alpha_1\neq \alpha_2.
\end{equation}
We still have $\delta^{j(k)}\leq \delta^{-1} \cdot r_k$, and we also have $\text{diam}(\overline{Q_{\alpha}^{j(k)+n}})\leq C_4 \cdot \delta^{j(k)+n}$, so we get
\begin{equation}
\text{diam}(\overline{Q_{\alpha}^{j(k)+n}}) \leq C_4 \cdot \delta^{n-1} \cdot r_k.
\end{equation}
By triangle inequality we get
\begin{equation}
\overline{Q_{\alpha}^{j(k)+n}} \subseteq B_d(x_k,r_k(1+ C_4 \cdot \delta^{n-1})).
\end{equation}
So, by taking
\begin{equation}\label{definition of n extra levels}
n:= \inf\{\tilde n \in \mathbb N \; | \; C_4 \cdot \delta^{\tilde n-1}<\epsilon\}
\end{equation}
and using (\ref{epsilon definition}), we have
\begin{equation}\label{j(k)+n covering isn't big}
\overline{Q_{\alpha}^{j(k)+n}} \subseteq B_d(x_k,M \cdot r_k).
\end{equation}
We observe that the definition (\ref{definition of n extra levels}) does not depend on the choice of $k \in \mathcal F$.\\
Equation (\ref{j(k)+n covering isn't big}) holds for all $k \in \mathcal F$ and for all $\alpha \in \mathcal H(k,n)$, and we get
\begin{equation}
 B_d(x_k,r_k) \subseteq \bigcup_{\alpha \in \mathcal H(k,n)} \overline{Q_{\alpha}^{j(k)+n}}  \subseteq  B_d(x_k,M \cdot r_k) \quad \forall k \in \mathcal F.
\end{equation}
We observe that there exists $\tilde N \in \mathbb N$ such that
\begin{equation}\label{claim finite covering}
|\mathcal H(k,n)| \leq \tilde N \quad \forall k \in \mathcal F.
\end{equation}
Indeed, for any $k \in \mathcal F$, using (\ref{j(k)+n ball inside the qube exists}) and (\ref{j(k)+n qubes are not intersected}), we get that the set
\begin{equation}
\mathcal S(k) := \bigcup_{\alpha \in \mathcal H(k,n)}B_d(z_{\alpha}^{j(k)+n},C_3 \cdot \delta^{j(k)+n})
\end{equation}
is a disjoint intersection of metric balls, and 
\begin{equation}
\mathcal S(k)\subseteq  \bigcup_{\alpha \in \mathcal H(k,n)} \overline{Q_{\alpha}^{j(k)+n}}\subseteq B_d(x_k,M \cdot r_k).
\end{equation}
However, the space $X$ is a $Q$-regular Ahlfors space, so, using (\ref{ahlfors regular condition for X}) and (\ref{j(k) level definition}), we get
\begin{align}
m(\mathcal S(k))=&\sum_{\alpha \in \mathcal H(k,n)}m(B_d(z_{\alpha}^{j(k)+n},C_3 \cdot \delta^{j(k)}\cdot \delta^n)) \geq \\
&\sum_{\alpha \in \mathcal H(k,n)}m(B_d(z_{\alpha}^{j(k)+n}, r^k \cdot \delta^n)) \geq \nonumber\\
&\sum_{\alpha \in \mathcal H(k,n)}C_1 \cdot r_k^Q\cdot \delta^{nQ}=\nonumber\\
&|\mathcal H(k,n)| \cdot C_1 \cdot r_k^Q\cdot \delta^{nQ}.\nonumber
\end{align}
On the other hand, applying the estimate (\ref{ahlfors regular condition for X}) to $B_d(x_k,M \cdot r_k)$ gives us
\begin{equation}
m(B_d(x_k,M \cdot r_k))\leq C_2 \cdot M^Q \cdot r_k^Q.
\end{equation}
However, $\mathcal S(k) \subseteq B_d(x_k,M \cdot r_k)$, so $m(\mathcal S(k))\leq m(B_d(x_k,M \cdot r_k))$ and we get
\begin{equation}
|\mathcal H(k,n)| \leq \frac{C_2}{ C_1 } \cdot M^Q \cdot\delta^{-nQ}\quad \forall k \in \mathcal F,
\end{equation}
proving the claim (\ref{claim finite covering}).\\

Consider $k \in \mathcal F$ and $\alpha \in I_k$. By definition \ref{Lambda map definition} we have
\begin{equation}
X \supseteq Q_{\alpha}^k \longmapsto \Lambda^{-1}(Q_{\alpha}^k)\subseteq \partial T.
\end{equation}
Let $S\subseteq X$ be a closed set. Let us define (with a slight abuse of notation)
\begin{equation}
\tilde \Lambda^{-1}(Q_{\alpha}^k \cap S):=\bigg\{x=\left(Q_{\beta_0}^0,Q_{\beta_1}^1,Q_{\beta_2}^2,\dots\right) \in \partial T \; \bigg| \; \Lambda(x)\in \overline{Q_{\alpha}^k} \cap S, \; \beta_k = \alpha\bigg\},
\end{equation}
and we also define $\tilde \Lambda^{-1}(Q_{\alpha}^k):=\tilde \Lambda^{-1}(Q_{\alpha}^k \cap X)$.\\
We have
\begin{equation}\label{Lamda tilde properties}
\tilde \Lambda^{-1}(Q_{\alpha}^k\cap S) \subseteq \Lambda^{-1}(Q_{\alpha}^k\cap S), \quad \text{and } \quad \Lambda(\tilde \Lambda^{-1}(Q_{\alpha}^k\cap S))=\overline{Q_{\alpha}^k}\cap S.
\end{equation}
Let $k \in \mathcal F$. Let us denote $N(k):=| H(k,n)|<\tilde N$. To simplify the notation let us denote the family of indices $H(k,n)$ by
\begin{equation}
H(k,n)=\{\alpha_1(k), \alpha_2(k), \dots , \alpha_{N(k)}(k)\}.
\end{equation}
For every $k \in \mathcal F$, for every $i=1,2,\dots, N(k)$ we define, when they exist,
\begin{equation}\label{definition of bigger qube capacity wise}
\overline{\eta^*}(k,i) \in \mathbb N \quad \text{and } \quad  \beta(k,i) \in I_{\overline{\eta^*}(k,i)},
\end{equation}
so that the qube $Q_{\beta(k,i)}^{\overline{\eta^*}(k,i)}$ is the smallest qube such that $Q_{\beta(k,i)}^{\overline{\eta^*}(k,i)}\supseteq Q_{\alpha_i(k)}^{j(k)+n}$ and such that
\begin{equation}\label{measure bigger than capacity for bigger qube}
\mathcal H_{\rho}^Q\left ( \tilde \Lambda^{-1}\left(Q_{\beta(k,i)}^{\overline{\eta^*}(k,i)}\right)\right) \geq C_{K_{\partial T,s},p}\left(\tilde \Lambda^{-1}\left(Q_{\alpha_i(k)}^{j(k)+n}\right)\right).
\end{equation}
Statements ii) and iii) in Theorem \ref{Christ decomposition} prove that $Q_{\beta(k,i)}^{\overline{\eta^*}(k,i)}$ is uniquely determined.\\
We claim that there exists a universal constant $\Psi_1>0$ such that for all $k \in \mathcal F$, for all $i=1,2,\dots,N(k)$ we have
\begin{equation}\label{claim about diameter estimate for bigger qubes}
\text{diam}_d\left(Q_{\beta(k,i)}^{\overline{\eta^*}(k,i)}\right) \leq \Psi_1 \cdot C_{K_{\partial T,s},p}\left(\tilde \Lambda^{-1}\left(Q_{\alpha_i(k)}^{j(k)+n}\right)\right)^{\frac 1 Q}.
\end{equation}
Indeed, by the Ahlfors $Q$-regularity of $X$ and $\partial T$ we can prove that there exist universal constants $0<\Psi_2<\Psi_3$ such that
\begin{equation}\label{the two measures are equal up to a product with a constant}
\Psi_2 \cdot  m\left(Q_{\beta(k,i)}^{\overline{\eta^*}(k,i)}\right)\leq\mathcal H_{\rho}^Q\left(\tilde \Lambda^{-1}\left(Q_{\beta(k,i)}^{\overline{\eta^*}(k,i)}\right)\right)\leq\Psi_3 \cdot m\left(Q_{\beta(k,i)}^{\overline{\eta^*}(k,i)}\right).
\end{equation}
Moreover, using the definition (\ref{definition of bigger qube capacity wise}) and using the Ahlfors $Q$-regular property (\ref{Ahlfors regular condition on partial T}) we can also prove that
\begin{equation}\label{capacity and measure of the bigger ball on boundary of T are equal up to a product with a constant}
C_{K_{\partial T,s},p}\left(\tilde \Lambda^{-1}\left(Q_{\alpha_i(k)}^{j(k)+n}\right)\right)\leq \mathcal H_{\rho}^Q\left(\tilde \Lambda^{-1}\left(Q_{\beta(k,i)}^{\overline{\eta^*}(k,i)}\right)\right) \leq \frac {K_2}{K_1} \delta^{-Q} \cdot C_{K_{\partial T,s},p}\left(\tilde \Lambda^{-1}\left(Q_{\alpha_i(k)}^{j(k)+n}\right)\right).
\end{equation}
By combining (\ref{the two measures are equal up to a product with a constant}) and (\ref{capacity and measure of the bigger ball on boundary of T are equal up to a product with a constant}) we get
\begin{equation}\label{first of 3 equation to combine in ahlfors proof}
 m\left(Q_{\beta(k,i)}^{\overline{\eta^*}(k,i)}\right) \leq  \frac{1}{\Psi_2}\frac {K_2}{K_1} \delta^{-Q} \cdot C_{K_{\partial T,s},p}\left(\tilde \Lambda^{-1}\left(Q_{\alpha_i(k)}^{j(k)+n}\right)\right).
\end{equation}
However, by statements iv) and v) in Theorem \ref{Christ decomposition} we have that there exists a ball \\$B_d\left(z_{\beta(k,i)}^{\overline{\eta^*}(k,i)},\tilde r_{k,i}\right)$ such that
\begin{equation}\label{second of 3 equation to combine in ahlfors proof}
B_d\left(z_{\beta(k,i)}^{\overline{\eta^*}(k,i)},\tilde r_{k,i}\right) \subseteq Q_{\beta(k,i)}^{\overline{\eta^*}(k,i)} \quad \text{and }\quad \text{diam}_d\left( Q_{\beta(k,i)}^{\overline{\eta^*}(k,i)}\right) \leq \frac{C_4}{C_3}\tilde r_{k,i}.
\end{equation}
By the Ahlfors $Q$-regularity condition (\ref{ahlfors regular condition for X}) we get
\begin{equation}\label{third of 3 equation to combine in ahlfors proof}
\tilde r_{k,i} \leq \frac{1}{C_1^{\frac{1}{Q}}} \cdot m\left(B_d\left(z_{\beta(k,i)}^{\overline{\eta^*}(k,i)},\tilde r_{k,i}\right) \right)^{\frac 1 Q} \leq \frac{1}{C_1^{\frac{1}{Q}}} \cdot m\left(Q_{\beta(k,i)}^{\overline{\eta^*}(k,i)}\right)^{\frac 1 Q}.
\end{equation}
Finally, by combining (\ref{first of 3 equation to combine in ahlfors proof}), (\ref{second of 3 equation to combine in ahlfors proof}) and (\ref{third of 3 equation to combine in ahlfors proof}), we get
\begin{equation}
\text{diam}_d\left( Q_{\beta(k,i)}^{\overline{\eta^*}(k,i)}\right)  \leq \frac{C_4}{C_3} \left(\frac{1}{C_1 \Psi_2}\frac {K_2}{K_1} \right)^{\frac 1 Q}\delta^{-1} \cdot C_{K_{\partial T,s},p}\left(\tilde \Lambda^{-1}\left(Q_{\alpha_i(k)}^{j(k)+n}\right)\right)^{\frac 1 Q}.
\end{equation}
By choosing 
\begin{equation}
\Psi_1:=\frac{C_4}{C_3} \left(\frac{1}{C_1 \Psi_2}\frac {K_2}{K_1} \right)^{\frac 1 Q}\delta^{-1}
\end{equation}
we prove the claim (\ref{claim about diameter estimate for bigger qubes}). The constant $\Psi_1$ is universal and does not depend on the choice of $k \in \mathcal F$ and $i=1,2,\dots,N(k)$.\\

Now we are going to define the value of the constant $\Psi$ defined in the hypotheses of this theorem.\\
Let $k \in \mathcal F$. We define
\begin{equation}
B^*_k:= B_d (x_k, \eta^*_{X,p}(x_k,M \cdot r_k)).
\end{equation}
By (\ref{j(k)+n covering isn't big}) we have $B_d (x_k, \eta^*_{X,p}(x_k,M \cdot r_k))\supseteq \overline{Q_{\alpha_i}^{j(k)+n}}$ for all $i=1,2,\dots,N(k)$. So, by the definition (\ref{eta star in X definition}) we get
\begin{align}
m(B^*_k)\geq& C_{K_{X,s},p}\left(B_d (x_k, \eta^*_{X,p}(x_k,M \cdot r_k))\right)\geq\\
&C_{K_{X,s},p}\left(\overline{Q_{\alpha_i(k)}^{j(k)+n}}\right) \quad \quad \quad \forall i=1,2,\dots,N(k).\nonumber
\end{align}
However, by construction, we have
\begin{equation}
\Lambda\left(\tilde \Lambda^{-1}\left(Q_{\alpha_i(k)}^{j(k)+n}\right)\right)=\overline{Q_{\alpha_i}^{j(k)+n}},
\end{equation}
and the set $\tilde \Lambda^{-1}\left(Q_{\alpha_i(k)}^{j(k)+n}\right) \subseteq \partial T$ is a closed set, so we can apply Theorem \ref{comparing capacities theorem} to get
\begin{equation}
m(B^*_k)\geq A_1 \cdot C_{K_{\partial T,s},p}\left(\tilde \Lambda^{-1}\left(Q_{\alpha_i(k)}^{j(k)+n}\right)\right)  \quad \forall i=1,2,\dots,N(k).
\end{equation}
Since $B_k^*$ is a metric ball in $X$ we can use the Ahlfors $Q$-regularity condition (\ref{ahlfors regular condition for X}) to get
\begin{equation}\label{bigger radius estimate for end of construction}
\eta^*_{X,p}(x_k,M \cdot r_k)\geq \frac{1}{C_2^{\frac 1 Q}}\cdot  m(B^*_k)^{\frac 1 Q}\geq \left(\frac{A_1}{C_2}\right)^{\frac 1 Q} \cdot  C_{K_{\partial T,s},p}\left(\tilde \Lambda^{-1}\left(Q_{\alpha_i(k)}^{j(k)+n}\right)\right) ^{\frac 1 Q}
\end{equation}
for $i=1,2,\dots,N(k).$\\
Now we observe that $\overline{Q_{\beta(k,i)}^{\overline{\eta^*}(k,i)}} \cap E_k \neq \emptyset$ because $Q_{\beta(k,i)}^{\overline{\eta^*}(k,i)}\supseteq Q_{\alpha_i(k)}^{j(k)+n}$ and because of (\ref{thinner cubes covering definition}), and we also recall that $E_k \subseteq B_d(x_k,r_k)$ by construction, so, by triangle inequality and by using the estimate (\ref{claim about diameter estimate for bigger qubes}), we get that for all $y \in \overline{Q_{\beta(k,i)}^{\overline{\eta^*}(k,i)}} $ we have
\begin{equation}\label{distance of bigger qubes from center of ball}
d(y,x_k)\leq r_k+\text{diam}_d\left(Q_{\beta(k,i)}^{\overline{\eta^*}(k,i)}\right)\leq r_k +  \Psi_1 \cdot C_{K_{\partial T,s},p}\left(\tilde \Lambda^{-1}\left(Q_{\alpha_i(k)}^{j(k)+n}\right)\right)^{\frac 1 Q}.
\end{equation}
Finally, we are going to define the constant $\Psi\geq 1$ mentioned in the thesis of this theorem by
\begin{equation}
\Psi:=1+ \Psi_1 \cdot \left(\frac{C_2}{A_1}\right)^{\frac 1 Q}.
\end{equation}
For all $k \in \mathcal F$ consider the ball $B_d(x_k,\Psi \cdot \eta^*_{X,p}(x_k,M \cdot r_k))$. \\
By construction $\eta^*_{X,p}(x_k,M \cdot r_k))\geq M \cdot r_k\geq r_k$, and we proved the estimate (\ref{bigger radius estimate for end of construction}), so we get
\begin{align}
\Psi \cdot \eta^*_{X,p}(x_k,M \cdot r_k)) =& \left(1+ \Psi_1\left(\frac{C_2}{A_1}\right)^{\frac 1 Q}\right) \cdot  \eta^*_{X,p}(x,M \cdot r_k)) \geq \\
&  r_k + \Psi_1\cdot \left(\frac{C_2}{A_1}\right)^{\frac 1 Q} \cdot \left(\frac{A_1}{C_2}\right)^{\frac 1 Q} \cdot  C_{K_{\partial T,s},p}\left(\tilde \Lambda^{-1}\left(Q_{\alpha_i(k)}^{j(k)+n}\right)\right) ^{\frac 1 Q}\geq \nonumber \\
&r_k +  \Psi_1 \cdot C_{K_{\partial T,s},p}\left(\tilde \Lambda^{-1}\left(Q_{\alpha_i(k)}^{j(k)+n}\right)\right)^{\frac 1 Q}.\nonumber
\end{align}
However, we proved (\ref{distance of bigger qubes from center of ball}), so we get
\begin{equation}\label{big qubes subsets of big ball}
Q_{\beta(k,i)}^{\overline{\eta^*}(k,i)} \subseteq B_d(x_k,\Psi \cdot \eta^*_{X,p}(x_k,M \cdot r_k))
\end{equation}
for every $k \in \mathcal F$, for every $i=1,2,\dots,N(k)$.\\

Now we finish the proof by proving (\ref{thesis quasi additivity ahlfors}).\\
Let $\{B_d(x_k,r_k)\}_{k\in \mathcal F}$ be a family of balls in $X$ such that $\{B_d(x_k,\Psi \cdot \eta^*_{X,p}(x_k,M \cdot r_k))\}_{k\in \mathcal F}$ is disjoint, Let $E \subseteq X$ be a compact set such that $E=\bigcup_k E_k$ and $E_k \subseteq B_d(x_k,r_k) \; \forall k$.\\
By construction 
\begin{equation}
E_k \subseteq \bigcup_{i=1}^{N(k)} \overline{Q_{\alpha_i(k)}^{j(k)+n}}
\end{equation}
and there exists a universal $\tilde N \in \mathbb N$ such that $N(k)\leq \tilde N$ for all $k \in \mathcal F$.\\
We are going to define for all $k \in \mathcal F$ and for all $i=1,2,\dots,\tilde N$
\begin{equation}
E_{k,i}:=\begin{cases}
E_k \cap \overline{Q_{\alpha_i(k)}^{j(k)+n}} \quad \text{if } i\leq N(k),\\
\emptyset \quad\quad\quad\quad\quad\quad\, \text{otherwise.}
\end{cases}
\end{equation}
It follows that 
\begin{equation}\label{double union definition for E}
E_k=\bigcup_{i=1}^{\tilde N} E_{k,i} \quad \forall k \in \mathcal F.
\end{equation}
Using (\ref{Lamda tilde properties}) we get
\begin{equation}
E_{k,i}=\begin{cases}
\Lambda\left(\tilde \Lambda^{-1}\left(E_k \cap Q_{\alpha_i(k)}^{j(k)+n}\right)\right) \quad \text{if }i\leq N(k),\\
\emptyset \quad\quad\quad\quad\quad\quad\quad\quad\quad\quad\quad\;\;\text{otherwise.}
\end{cases}
\end{equation}
To simplify the notation let us denote
\begin{equation}
E_k \cap Q_{\alpha_i(k)}^{j(k)+n}:= \emptyset, \quad \Lambda^{-1}\left(E_k \cap Q_{\alpha_i(k)}^{j(k)+n}\right):= \emptyset
\end{equation}
whenever $i>N(k)$.\\
So we apply the subadditivity of the capacity to (\ref{double union definition for E}) to get
\begin{equation}
\sum_{k \in \mathcal F} C_{K_{X,s},p}\left( E_{k}\right)= \sum_{k \in \mathcal F}C_{K_{X,s},p}\left(\bigcup_{i=1}^{\tilde N}E_{k,i}\right)\leq \sum_{i=1}^{\tilde N}\sum_{k \in \mathcal F}C_{K_{X,s},p}\left(E_{k,i}\right).
\end{equation}
Let $i<\tilde N$ be a fixed index. Let $k \in \mathcal F$.\\
If $i>N(k)$ then, by definition 
\begin{equation}
C_{K_{\partial T,s},p}\left(\Lambda^{-1}\left(E_k \cap Q_{\alpha_i(k)}^{j(k)+n}\right)\right)=C_{K_{\partial T,s},p}(\emptyset)=0.
\end{equation}
Otherwise, we have $E_{k,i}=\Lambda\left(\tilde \Lambda^{-1}\left(E_k \cap Q_{\alpha_i(k)}^{j(k)+n}\right)\right)$, so we can apply the estimate in Theorem \ref{comparing capacities theorem} to get
\begin{equation}\label{sum to estimate with theorem 5}
\sum_{k \in \mathcal F} C_{K_{X,s},p}\left( E_{k}\right)\leq  A_2 \cdot  \sum_{i=1}^{\tilde N}\sum_{k \in \mathcal F}C_{K_{\partial T,s},p}\left(\tilde \Lambda^{-1}\left( E_k \cap Q_{\alpha_i(k)}^{j(k)+n}\right)\right).
\end{equation}
We observe that, for any fixed $i<\tilde N$, the sum
\begin{equation}
\sum_{k \in \mathcal F}C_{K_{\partial T,s},p}\left(\tilde \Lambda^{-1}\left( E_k \cap Q_{\alpha_i(k)}^{j(k)+n}\right)\right)
\end{equation}
can be estimated using Theorem \ref{Quasi additivity for boundary of tree theorem}.\\
Indeed, we have
\begin{equation}
\Lambda^{-1}\left( E_k \cap Q_{\alpha_i(k)}^{j(k)+n}\right)=\bigg \{x=(Q_{\beta_0}^0,Q_{\beta_1}^1,\dots) \in \partial T \; \bigg | \;\Lambda(x)\in  E_k \cap \overline{Q_{\alpha_i(k)}^{j(k)+n}}, \; \beta_{j(k)+n}=\alpha_i(k)\bigg\},
\end{equation}
so it follows that, for any choice of $w_{k,i}\in\Lambda^{-1}\left( E_k \cap Q_{\alpha_i(k)}^{j(k)+n}\right)$, we have
\begin{equation}
\Lambda^{-1}\left( E_k \cap Q_{\alpha_i(k)}^{j(k)+n}\right) \subseteq B_{\rho}(w_{k,i},\delta^{j(k)+n}) \subseteq \partial T.
\end{equation}
Moreover, we have
\begin{equation}
\Lambda^{-1}\left( Q_{\alpha_i(k)}^{j(k)+n}\right) =B_{\rho}(w_{k,i},\delta^{j(k)+n}) \subseteq \partial T \quad \text{and } \quad \Lambda^{-1}\left(Q_{\beta(k,i)}^{\overline{\eta^*}(k,i)} \right) = B_{\rho}(w_{k,i},\delta^{\overline{\eta^*}(k,i)}) \subseteq \partial T.
\end{equation}
Using (\ref{measure bigger than capacity for bigger qube}) we get that $B_{\rho}(w_{k,i},\delta^{\overline{\eta^*}(k,i)})$ is the smallest possible ball such that 
\begin{equation}
B_{\rho}(w_{k,i},\delta^{\overline{\eta^*}(k,i)})\supseteq B_{\rho}(w_{k,i},\delta^{j(k)+n}),
\end{equation}
and
\begin{equation}
\mathcal H_{\rho}^Q\left ( B_{\rho}(w_{k,i},\delta^{\overline{\eta^*}(k,i)}) \right) \geq C_{K_{\partial T,s},p}\left(B_{\rho}(w_{k,i},\delta^{j(k)+n}) \right).
\end{equation}
So, using definition (\ref{bigger radius in tree boundary definition}),  we have
\begin{equation}
\eta^*_p(w_{k,i},\delta^{j(k)+n})=\delta^{\overline{\eta^*}(k,i)}.
\end{equation}
In (\ref{big qubes subsets of big ball}) we proved
\begin{equation}
Q_{\beta(k,i)}^{\overline{\eta^*}(k,i)} \subseteq B_d(x_k,\Psi \cdot \eta^*_{X,p}(x,M \cdot r_k)) \quad \forall k \in \mathcal F, \; \forall i\leq \tilde N,
\end{equation}
but the family $\{B_d(x_k,\Psi \cdot \eta^*_{X,p}(x,M \cdot r_k))\}_{k \in \mathcal F}$ is disjoint by hypothesis, so it follows that the family
\begin{equation}
\left\{ \tilde \Lambda^{-1}\left(Q_{\beta(k,i)}^{\overline{\eta^*}(k,i)} \right)\right\}_{k \in \mathcal F}=\left\{ B_{\rho}(w_{k,i},\delta^{\overline{\eta^*}(k,i)})\right\}_{k \in \mathcal F}
\end{equation}
is disjoint for every fixed $i\leq \tilde N$.\\
So the hypotheses of Theorem \ref{Quasi additivity for boundary of tree theorem} are satisfied, and we may apply (\ref{Thesis}) to (\ref{sum to estimate with theorem 5}) to get
\begin{equation}\label{sum to fix with subadditivity}
\sum_{k \in \mathcal F} C_{K_{X,s},p}\left( E_{k}\right)\leq A \cdot A_2 \cdot  \sum_{i=1}^{\tilde N} C_{K_{\partial T,s},p}\left(\bigcup_{k \in \mathcal F}\tilde \Lambda^{-1}\left( E_k \cap Q_{\alpha_i(k)}^{j(k)+n}\right)\right).
\end{equation}
We recall the following quasi-additivity formula for the capacity of a finite union of sets:
\begin{equation}
\sum_{i=1}^{\tilde N} C_{K_{\partial T,s},p}\left(S_i\right)\leq \tilde N \cdot C_{K_{\partial T,s},p}\left(\bigcup_{i=1}^{\tilde N} S_i\right).
\end{equation}
Applying this formula to (\ref{sum to fix with subadditivity}) we get
\begin{equation}\label{second to last equation in ahlfors proof}
\sum_{k \in \mathcal F} C_{K_{X,s},p}\left( E_{k}\right)\leq A \cdot A_2 \cdot \tilde N \cdot  C_{K_{\partial T,s},p}\left(\bigcup_{i=1}^{\tilde N}\bigcup_{k \in \mathcal F}\tilde \Lambda^{-1}\left( E_k \cap Q_{\alpha_i(k)}^{j(k)+n}\right)\right).
\end{equation}
Now, to finish the proof, we observe that
\begin{equation}
\Lambda\left(\bigcup_{i=1}^{\tilde N}\bigcup_{k \in \mathcal F}\tilde \Lambda^{-1}\left( E_k \cap Q_{\alpha_i(k)}^{j(k)+n}\right)\right)=\bigcup_{i=1}^{\tilde N}\bigcup_{k \in \mathcal F}  E_k \cap \overline{Q_{\alpha_i(k)}^{j(k)+n}}=\bigcup_{k \in \mathcal F} E_k=E.
\end{equation}
So we may apply Theorem \ref{comparing capacities theorem} to (\ref{second to last equation in ahlfors proof}) to get
\begin{equation}\label{thesis proved in quasi additivity ahlfors regular}
\sum_{k \in \mathcal F} C_{K_{X,s},p}\left( E_{k}\right)\leq A \cdot A_2^2 \cdot \tilde N \cdot  C_{K_{X,s},p}\left(E\right).
\end{equation}
The constant $\tilde A :=A \cdot A_2^2 \cdot \tilde N$ only depends on $X$, $s$, $p$, $M$ and on the chosen Christ decomposition, finishing the proof for the case $M>1$.\\

The case $M=1$ follows as a corollary from the case $M>1$.\\
Indeed, choose an arbitrary $\tilde M>1$. Let $\{B_d(x_k,r_k)\}_{k\in \mathcal F}$ be a countable family of balls in $X$ such that the family $\{B_d(x_k,\tilde \Psi \cdot \eta^*_{X,p}(x_k, r_k))\}_{k\in \mathcal F}$ is disjoint. Here the constant $\tilde \Psi=\tilde \Psi(X,p,s,\tilde M)>\Psi$ will be fixed later. Then there exists $\hat N = \hat N(X,s,p,\tilde M)\in \mathbb N$ such that, up to a proper choice of the indices of the family $\{B_d(x_k,r_k)\}_{k\in \mathcal F}$, we have
\begin{equation}\label{property to generalize quasi additivity sum of infinite things}
\bigcap_{k=\hat N}^{+\infty} B_d(x_k,\tilde \Psi \cdot \eta^*_{X,p}(x_k,\tilde M \cdot r_k))= \emptyset.
\end{equation}
This property holds if we choose $\tilde \Psi$ such that 
\begin{equation}
\tilde \Psi \cdot \eta^*_{X,p}(x_k,  r_k)\geq \Psi \cdot \eta^*_{X,p}(x_k,\tilde M \cdot r_k)
\end{equation}
for all $k\geq \hat N$.\\
We claim that such $\hat N$ exists because of the properties of the Riesz capacity (see Proposition \ref{capacity of ball general formula proposition}). Indeed, there exists $\hat C>1$ and $\hat r>0$ universal constants such that, for all $x \in X$, for all $0<r<\hat r$ we have
\begin{equation}
\tilde C \cdot C_{K_{X,s},p}\left( B_d(x,r)\right)\geq  C_{K_{X,s},p}\left( B_d(x,\tilde M \cdot r)\right),
\end{equation}
which entails
\begin{equation}
\eta^*_{X,p}(x_k, r_k)\gtrsim_{X,s,p,\tilde M} \tilde C^{Q} \cdot \eta^*_{X,p}(x_k,\tilde M \cdot r_k),
\end{equation}
for all $ r_k<\hat r$. The claim follows from the compactness of $X$, because there exists $\hat N \in \mathbb N$ such that $r_k<\hat r$ for all $k\geq \hat N$. If such $\hat N$ did not exist then we would have
\begin{equation}
m(X)\geq m\left(\bigcup_{k\in \mathbb N} B_d(x_k,r_k)\right)=\sum _{k\in \mathbb N}m( B_d(x_k,r_k))=+\infty,
\end{equation}
which contradicts the compactness of $X$.\\
So we proved that (\ref{property to generalize quasi additivity sum of infinite things}) holds as long as we choose 
\begin{equation}
\tilde \Psi := \tilde C \cdot \Psi.
\end{equation}
Now we consider an arbitrary compact set $E \subseteq X$ such that $E=\bigcup_k E_k$ and $E_k \subseteq B_d(x_k,r_k) \; \forall k$. We apply the quasi-additivity formula (\ref{thesis proved in quasi additivity ahlfors regular}) to the family $\{E_k\}_{k \in \mathbb N}$ and we get
\begin{equation}
\sum_{k =\hat N}^{+\infty} C_{K_{X,s},p}\left( E_{k}\right)\leq \tilde A \cdot  C_{K_{X,s},p}\left(\bigcup_{k =\hat N}^{+\infty}E_k\right).
\end{equation}
Finally, we apply the finite quasi-additivity formula and we get
\begin{eqnarray}
\sum_{k \in \mathcal F} C_{K_{X,s},p}\left( E_{k}\right)&=&\sum_{k =1}^{\hat N-1} C_{K_{X,s},p}\left( E_{k}\right)+\sum_{k =\hat N}^{+\infty} C_{K_{X,s},p}\left( E_{k}\right)\leq\\
& &\sum_{k =1}^{\hat N-1} C_{K_{X,s},p}\left( E_{k}\right)+\tilde A \cdot  C_{K_{X,s},p}\left(\bigcup_{k =\hat N}^{+\infty}E_k\right)\leq\nonumber\\
& &\tilde A\left[\sum_{k =1}^{\hat N-1} C_{K_{X,s},p}\left( E_{k}\right)+  C_{K_{X,s},p}\left(\bigcup_{k =\hat N}^{+\infty}E_k\right)\right]\leq\nonumber\\
 & & \tilde A \cdot \hat N \cdot  C_{K_{X,s},p}\left(\bigcup_{k =1}^{\hat N-1}E_k \cup \bigcup_{k =\hat N}^{+\infty}E_k\right),\nonumber
\end{eqnarray}
so we get 
\begin{equation}
\sum_{k \in \mathcal F} C_{K_{X,s},p}\left( E_{k}\right)\lesssim_{(X,s,p)}  C_{K_{X,s},p}\left(E\right),
\end{equation}
so, up to choosing the new values of the constants $\Psi$ and $\tilde A$ for the case $M=1$, the theorem is proved.
\endproof
\section{Harmonic extension}
In this section we define the harmonic extension of a function defined over an Ahlfors-regular space and we enunciate and prove several properties of the Harmonic extension.\\
\subsection{Dyadic Poisson Integral and Riesz kernel}
Let $(X,d,m)$ be an Ahlfors $Q$-regular space. In the following part of this work we will be considering the space $X \times (0,+\infty)$ with the metric
\begin{equation}
\rho\left((x_1,y_1),(x_2,y_2)\right):= \max\{d(x_1,x_2),|y_1-y_2|\}.
\end{equation}

We are now going to define the harmonic extension on $X \times (0,+\infty)$.
\begin{Def}[Poisson Integral in $X \times (0,+\infty)$]\label{PI definition}
Let $f \in L^p(X)$. We define the Poisson Integral
\begin{equation}
PI(f)(x,y):= \int_{X} C(x,y) \cdot \frac{1}{y^Q} \sum_{k=0}^{+\infty} \frac{\chi_{B_d(x,2^k y)}(z)}{2^{(Q+1)k}}f(z)dm(z).
\end{equation}
Here $C(x,y)$ is the constant that normalizes the Poisson kernel, i.e.
\begin{equation}
C(x,y):=\left[ \int_{X}  \frac{1}{y^Q} \sum_{k=0}^{+\infty} \frac{\chi_{B_d(x,2^k y)}(z)}{2^{(Q+1)k}}dm(z)\right]^{-1}.
\end{equation}
\end{Def}
\begin{Rem}\label{constant in Poisson Kernel is bounded}
There exist constants $0<C_1\leq C_2<+\infty$ depending only on the choice of $X$ and $Q$ such that
\begin{equation}\label{boundedness of constant in Poisson Integral}
C_1\leq C(x,y)\leq C_2 \quad \forall x \in X, \; \forall y>0,
\end{equation}
i.e.
\begin{equation}
C(x_1,y_1)\approx_{(X,Q)} C(x_2,y_2) \quad \text{for all }(x_1,y_1),(x_2,y_2) \in X \times(0,+\infty).
\end{equation}
\end{Rem}

We will use a dyadicization of the Riesz kernel to prove a property of the Poisson Integral.
\begin{Lem}[Discretization of the Riesz kernel]\label{Discretization of the Riesz kernel lemma}
Let $g \in L^p(X)$. Then
\begin{equation}
K_{X,s} * g(x) \approx_{(Q,s)} \int_X \sum_{j=-\infty}^{+\infty} \frac{\chi_{B_d(x,2^j)}(z)}{2^{(Qs)j}} g(z)dm(z)
\end{equation}
for every $x \in X$.
\end{Lem}
\proof
Let $x \in X$. For every $z\neq x$ we have
\begin{equation}
\sum_{j=-\infty}^{+\infty} \frac{\chi_{B_d(x,2^j)}(z)}{2^{(Qs)j}}=\sum_{\underset{j>\log_2(d(x,z))}{j \text{ such that}}}\frac{1}{2^{(Qs)j}},
\end{equation}
and we have
\begin{equation}
\sum_{\underset{j>\log_2(d(x,z))}{j \text{ such that}}}\frac{1}{2^{(Qs)j}}= \sum_{\underset{j\geq\lceil\log_2(d(x,z))\rceil}{j \text{ such that}}}\frac{1}{2^{(Qs)j}}\leq \frac{1}{1-2^{Qs}} \frac{1}{2^{(Qs)(\log_2(d(x,z)))}}=\frac{1}{1-2^{Qs}}\frac{1}{d(x,z)^{Qs}},
\end{equation}
\begin{equation}
\begin{aligned}
\sum_{\underset{j>\log_2(d(x,z))}{j \text{ such that}}}\frac{1}{2^{(Qs)j}}=& \sum_{\underset{j\geq\lceil\log_2(d(x,z))\rceil}{j \text{ such that}}}\frac{1}{2^{(Qs)j}}\geq \frac{1}{1-2^{Qs}} \frac{1}{2^{(Qs)(\log_2(d(x,z))+1)}}=\\
&\frac{1}{2^{Qs}-2^{2Qs}}\frac{1}{d(x,z)^{Qs}}.
\end{aligned}
\end{equation}
But $K_{X,s}(x,z)=1/d(x,z)^{Qs}$, so we get
\begin{equation}
K_{X,s}(x,z)\approx_{(Q,s)} \sum_{j=-\infty}^{+\infty} \frac{\chi_{B_d(x,2^j)}(z)}{2^{(Qs)j}}
\end{equation}
for all $z\neq x$. Since $m(\{x\})=0$ we get
\begin{equation}
K_{X,s} * g(x) \approx_{(Q,s)} \int_X \sum_{j=-\infty}^{+\infty} \frac{\chi_{B_d(x,2^j)}(z)}{2^{(Qs)j}} g(z)dm(z)
\end{equation}
for all $g \in L^p(X)$, for all $x \in X$, and the lemma is proved.
\endproof
\subsection{Commutative convolution-like property}
Now we will prove that the order of the Poisson Integral and the Riesz potential can be exchanged up to a universal multiplicative constant. The proof of this property in $\mathbb R^{n+1}_+$ trivially follows from the commutative property of the convolution. In the setting of compact Ahlfors-regular spaces this proof is based on a geometrical property of the dyadicization of the Riesz kernel, and on the dyadic nature of the Poisson Integral we defined.
\begin{Lem}[Exchanging the order of the Poisson Integral and of the Riesz potential]\label{Lemma exchange PI and Riesz potential}
Let $f \in L^p(X)$, let $y>0$. Then
\begin{equation}\label{thesis associativity convolution lemma}
K_{X,s}*(PI(f)(\cdot,y))(x)\approx_{(X,Q,s)} PI(K_{X,s} * f)(x,y)
\end{equation}
for all $x \in X$.
\end{Lem}
\proof
By Lemma \ref{Discretization of the Riesz kernel lemma} we know that
\begin{equation}
K_{X,s}*(PI(f)(\cdot,y))(x)\approx_{(Q,s)} \int_X\sum_{j=-\infty}^{+\infty} \frac{\chi_{B_d(x,2^j)}(z)}{2^{(Qs)j}} PI(f)(z,y)dm(z),
\end{equation}
and
\begin{equation}
PI(K_{X,s} * f)(x,y)\approx_{(Q,s)}PI\left(\int_X\sum_{j=-\infty}^{+\infty} \frac{\chi_{B_d(\cdot,2^j)}(z)}{2^{(Qs)j}}f(z)dm(z)\right)(x,y),
\end{equation}
so we will prove the statement by proving 
\begin{equation}\label{thing to prove for associativity convolution}
\underset{X}{\int}\sum_{j=-\infty}^{+\infty} \frac{\chi_{B_d(x,2^j)}(z)}{2^{(Qs)j}} PI(f)(z,y)dm(z)\approx_{(X,Q,s)} PI\left(\underset{X}{\int}\sum_{j=-\infty}^{+\infty} \frac{\chi_{B_d(\cdot,2^j)}(z)}{2^{(Qs)j}}f(z)dm(z)\right)(x,y)
\end{equation}
for all $f \in L^p(X)$, for all $x,z \in X$, for all $y>0$.\\
Let $L.H.S.$ and $R.H.S.$ denote the left hand side and the right hand side of (\ref{thing to prove for associativity convolution}) respectively. We compute
\begin{align}
L.H.S.=&\int_X\sum_{j=-\infty}^{+\infty} \frac{\chi_{B_d(x,2^j)}(z)}{2^{(Qs)j}}\int_X\frac{ C(z,y)}{y^Q} \sum_{k=0}^{+\infty} \frac{\chi_{B_d(z,2^k y)}(w)}{2^{(Q+1)k}}f(w)dm(w)dm(z)=\\
&\int_X\int_X\frac{ C(z,y)}{y^Q}\sum_{j=-\infty}^{+\infty} \sum_{k=0}^{+\infty}\frac{\chi_{B_d(x,2^j)}(z)\chi_{B_d(z,2^k y)}(w)}{2^{[(Q+1)k+(Qs)j]}}f(w)dm(z)dm(w).\nonumber
\end{align}
Remark \ref{constant in Poisson Kernel is bounded} entails $C(z,y)\approx_{(X)} C(x,y)$, and we have 
\begin{equation*}
\chi_{B_d(z,2^k y)}(w)=\chi_{B_d(w,2^k y)}(z),
\end{equation*}
so we get
\begin{align}
L.H.S.\approx_{(X)}\frac{ C(x,y)}{y^Q}\int_X\left[\int_X\sum_{j=-\infty}^{+\infty} \sum_{k=0}^{+\infty}\frac{\chi_{B_d(x,2^j)}(z)\chi_{B_d(w,2^k y)}(z)}{2^{[(Q+1)k+(Qs)j]}}dm(z)\right]f(w)dm(w).
\end{align}
Now we compute
\begin{align}
R.H.S.=&PI\left(\int_X\sum_{j=-\infty}^{+\infty} \frac{\chi_{B_d(\cdot,2^j)}(z)}{2^{(Qs)j}}f(z)dm(z)\right)(x,y)=\\
&\int_X\frac{ C(x,y)}{y^Q}\sum_{k=0}^{+\infty} \frac{\chi_{B_d(x,2^k y)}(w)}{2^{(Q+1)k}}\int_X  \sum_{j=-\infty}^{+\infty} \frac{\chi_{B_d(w,2^j)}(z)}{2^{(Qs)j}} f(z)dm(z)dm(w)=\nonumber\\
&\int_X\int_X\frac{ C(x,y)}{y^Q}\sum_{j=-\infty}^{+\infty} \sum_{k=0}^{+\infty}\frac{\chi_{B_d(x,2^k y)}(w)\chi_{B_d(w,2^j)}(z)}{2^{[(Q+1)k+(Qs)j]}}f(z)dm(z)dm(w)=\nonumber \\
&\frac{ C(x,y)}{y^Q}\int_X\left[\int_X\sum_{j=-\infty}^{+\infty} \sum_{k=0}^{+\infty}\frac{\chi_{B_d(x,2^k y)}(w)\chi_{B_d(z,2^j)}(w)}{2^{[(Q+1)k+(Qs)j]}}dm(w)\right] f(z)dm(z).\nonumber
\end{align}
We rename the bound variables in the last equation to get
\begin{align}
R.H.S.=\frac{ C(x,y)}{y^Q}\int_X\left[\int_X\sum_{j=-\infty}^{+\infty} \sum_{k=0}^{+\infty}\frac{\chi_{B_d(x,2^k y)}(z)\chi_{B_d(w,2^j)}(z)}{2^{[(Q+1)k+(Qs)j]}}dm(z)\right] f(w)dm(w).
\end{align}
So, to prove that $L.H.S. \approx_{(X,Q,s)} R.H.S.$ for all $f \in L^p(X)$, for all $x \in X$, for all $y>0$, it is sufficient to prove that
\begin{equation}
\begin{aligned}
&\int_X\sum_{j=-\infty}^{+\infty} \sum_{k=0}^{+\infty}\frac{\chi_{B_d(x,2^j)}(z)\chi_{B_d(w,2^k y)}(z)}{2^{[(Q+1)k+(Qs)j]}}dm(z)\approx _{_{(X,Q,s)}} \\
&\int_X\sum_{j=-\infty}^{+\infty} \sum_{k=0}^{+\infty}\frac{\chi_{B_d(x,2^k y)}(z)\chi_{B_d(w,2^j)}(z)}{2^{[(Q+1)k+(Qs)j]}}dm(z)
\end{aligned}
\end{equation}
for all $x,w \in X$, for all $y>0$.\\
We will prove the stronger statement
\begin{equation}\label{literally final thing to prove}
\sum_{j=-\infty}^{+\infty} \sum_{k=0}^{+\infty}\frac{m(B_d(x,2^j)\cap B_d(w,2^k y))}{2^{[(Q+1)k+(Qs)j]}}\approx _{(X,Q,s)}  \sum_{j=-\infty}^{+\infty} \sum_{k=0}^{+\infty}\frac{m( B_d(w,2^j)\cap B_d(x,2^k y))}{2^{[(Q+1)k+(Qs)j]}}
\end{equation}
for all $x,w \in X$, for all $y>0$.\\
The claim is trivial when $x=w$. Let us consider $x\neq w$. So we have $d(x,w)>0$, and we may define
\begin{equation}\label{definition of tilde j}
\tilde{j}(x,w)=\left\lceil \log_2(d(x,w)) \right\rceil,
\end{equation}
which is the first index $j$ such that $w \in \overline{B_d(x,2^j)}$.\\
Let us consider the following equation
\begin{equation}
\sum_{j=-\infty}^{+\infty} \sum_{k=0}^{+\infty}\frac{m(B_d(x,2^j)\cap B_d(w,2^k y))}{2^{[(Q+1)k+(Qs)j]}}=(I)+(II)+(III),
\end{equation}
where
\begin{align}
(I)=&\sum_{j=-\infty}^{\tilde{j}(x,w)-2} \sum_{k=0}^{+\infty}\frac{m(B_d(x,2^j)\cap B_d(w,2^k y))}{2^{[(Q+1)k+(Qs)j]}},\\
(II)=&\sum_{j=\tilde{j}(x,w)-1}^{\tilde{j}(x,w)} \sum_{k=0}^{+\infty}\frac{m(B_d(x,2^j)\cap B_d(w,2^k y))}{2^{[(Q+1)k+(Qs)j]}},\\
(III)=&\sum_{j=\tilde{j}(x,w)+1}^{+\infty} \sum_{k=0}^{+\infty}\frac{m(B_d(x,2^j)\cap B_d(w,2^k y))}{2^{[(Q+1)k+(Qs)j]}}.
\end{align}
Now we are going to give estimates of the values of $(I)$, $(II)$ and $(III)$.\\

We start from estimating $(I)$. Consider $j\leq\tilde{j}(x,w)-2$. We define
\begin{equation}
\overline{k}(j)=\max\left\{\left\lfloor \log_2\left(\frac{d(x,w)-2^j}{y}\right)\right\rfloor+1,0\right\},
\end{equation}
which is the smallest index $k\geq0$ such that $B_d(x,2^j) \cap B_d(w,2^k y) \neq \emptyset$.\\
We write
\begin{align}
(I)=\sum_{j=-\infty}^{\tilde{j}(x,w)-2}\bigg[& \sum_{k=0}^{\overline{k}(j)-1}\frac{m(B_d(x,2^j)\cap B_d(w,2^k y))}{2^{[(Q+1)k+(Qs)j]}}+\\
&\sum_{k=\overline{k}(j)}^{\overline{k}(j)+1}\frac{m(B_d(x,2^j)\cap B_d(w,2^k y))}{2^{[(Q+1)k+(Qs)j]}}+\nonumber\\
&\sum_{k=\overline{k}(j)+2}^{+\infty}\frac{m(B_d(x,2^j)\cap B_d(w,2^k y))}{2^{[(Q+1)k+(Qs)j]}}\bigg],\nonumber
\end{align}
with the convention $\sum_{k=a}^{b}\varphi(k):=0$ if $b<a$.\\
By definition of $\overline{k}(j)$ we have $B_d(x,2^j)\cap B_d(w,2^k y)= \emptyset$ if $k<\overline{k}(j)$.\\
Now we consider $k\geq \overline{k}(j)+2$. We have
\begin{align}
k\geq\overline{k}(j)+2=\max\left\{\left\lfloor \log_2\left(\frac{d(x,w)-2^j}{y}\right)\right\rfloor+1,0\right\}+2 \geq \log_2\left(\frac{d(x,w)-2^j}{y}\right)+2,
\end{align}
hence we get
\begin{equation}\label{first step triangle inequality for (I)}
2^k y\geq  2^{\log_2\left(\frac{d(x,w)-2^j}{y}\right)+2}y=4 \left(d(x,w)-2^j\right).
\end{equation}
Now we consider an arbitrary point $v \in B_d(x,2^j)$. By triangle inequality we have
\begin{align}\label{third step triangle inequality for (I)}
d(v,w)\leq d(v,x)+d(x,w)\leq 2^j+d(x,w).
\end{align}
However, we set $j\leq\tilde{j}(x,w)-2$, so we get
\begin{equation}
j\leq\tilde{j}(x,w)-2=\left\lceil \log_2(d(x,w)) \right\rceil-2\leq \log_2(d(x,w))-1,
\end{equation}
hence
\begin{equation}\label{second step triangle inequality for (I)}
2^j\leq 2^{\log_2(d(x,w))-1}\leq \frac 1 2  d(x,w).
\end{equation}
So from (\ref{first step triangle inequality for (I)}) and (\ref{second step triangle inequality for (I)}) we get
\begin{equation}\label{fourth step triangle inequality for (I)}
2^k y\geq 2 d(x,w) \geq 2^j+\frac{3}{2}d(x,w).
\end{equation}
From (\ref{third step triangle inequality for (I)}) and (\ref{fourth step triangle inequality for (I)}) we get $d(v,w)\leq2^k y$, so we proved that $B_d(x,2^j)\subseteq B_d(w,2^k y)$ for all $x,w \in X$, for all $j\leq\tilde{j}(x,w)-2$, for all $k\geq \overline{k}(j)+2$.\\
So we proved that
\begin{align}
(I)=\sum_{j=-\infty}^{\tilde{j}(x,w)-2}\bigg[& \sum_{k=0}^{\overline{k}(j)-1}\frac{m(\emptyset)}{2^{[(Q+1)k+(Qs)j]}}+\\
&\sum_{k=\overline{k}(j)}^{\overline{k}(j)+1}\frac{m(B_d(x,2^j)\cap B_d(w,2^k y))}{2^{[(Q+1)k+(Qs)j]}}+\nonumber\\
&\sum_{k=\overline{k}(j)+2}^{+\infty}\frac{m(B_d(x,2^j))}{2^{[(Q+1)k+(Qs)j]}}\bigg].\nonumber
\end{align}
Now we use the equations
\begin{equation}\label{first of 3 equations to estimate measure of balls}
m(\emptyset)=0,
\end{equation}
\begin{equation}\label{second of 3 equations to estimate measure of balls}
0\leq m(B_d(x,2^j))\leq m(B_d(x,2^j)\cap B_d(w,2^k y)),
\end{equation}
\begin{equation}\label{third of 3 equations to estimate measure of balls}
m(B_d(x,r)) \approx_{(X)}  r^Q \quad \forall x \in X, \text{ for } 0<r<\text{diam}(X),
\end{equation}
to get the following equation:
\begin{align}\label{lower estimate for (I): start}
(I)\geq& \sum_{j=-\infty}^{\tilde{j}(x,w)-2}\sum_{k=\overline{k}(j)+2}^{+\infty}\frac{m(B_d(x,2^j))}{2^{[(Q+1)k+(Qs)j]}}\gtrsim_{(X)}\\
&\sum_{j=-\infty}^{\tilde{j}(x,w)-2}\sum_{k=\overline{k}(j)+2}^{+\infty}\frac{2^{jQ}}{2^{[(Q+1)k+(Qs)j]}}=\nonumber\\
&\sum_{j=-\infty}^{\tilde{j}(x,w)-2}\frac{1}{1-\frac{1}{2^{Q+1}}}\frac{1}{2^{(Q+1)(\overline{k}(j)+2)}}\frac{2^{jQ}}{2^{(Qs)j}}.\nonumber
\end{align}
We remark that we may use the estimate (\ref{third of 3 equations to estimate measure of balls}) because for all $j\leq \tilde{j}(x,w)-1$ we have $w \not \in B_d(x,2^j)$, hence $2^j < \text{diam}(X)$ for all $j\leq \tilde{j}(x,w)-2$.\\
Now we estimate the value of $\overline{k}(j)$. Consider $j=\tilde{j}(x,w)-2$. We have
\begin{equation}
\overline{k}(\tilde{j}(x,w)-2)=\max\left\{\left\lfloor \log_2\left(\frac{d(x,w)-2^{\tilde{j}(x,w)-2}}{y}\right)\right\rfloor+1,0\right\}>  \log_2\left(\frac{d(x,w)-2^{\tilde{j}(x,w)-2}}{y}\right).
\end{equation}
So we get
\begin{equation}
2^{\overline{k}(\tilde{j}(x,w)-2)} y> d(x,w)-2^{\tilde{j}(x,w)-2}.
\end{equation}
However we use equation (\ref{definition of tilde j}) to get
\begin{equation}\label{lower estimate (I): how big is the first overline k}
2^{\overline{k}(\tilde{j}(x,w)-2)} y> d(x,w)-2^{\log_2(d(x,w))+1-2}=\frac 1 2 d(x,w).
\end{equation}
Now let $j<\tilde{j}(x,w)-2)$. From (\ref{lower estimate (I): how big is the first overline k}) we get
\begin{equation}\label{lower estimate for (I): how big is the radius}
2^{\overline{k}(\tilde{j}(x,w)-2)+1} y>d(x,w),
\end{equation}
hence
\begin{equation}
B_d(w,2^{\overline{k}(\tilde{j}(x,w)-2)+1} y)\cap B_d(x,2^j)\neq 0 \quad \text{for all }j<\tilde{j}(x,w)-2.
\end{equation}
However, from the definition of $\overline{k}(j)$, it follows that $\overline{k}(j)$ is the smallest index $k$ such that $B_d(w,2^{k} y)\cap B_d(x,2^j)\neq 0$, so we proved 
\begin{equation}\label{lower estimate for (I): how big is overline k for all j}
\overline{k}(j)\leq \overline{k}(\tilde{j}(x,w)-2)+1 \quad \text{for all } j\leq\tilde{j}(x,w)-2.
\end{equation}
So we use equations (\ref{lower estimate for (I): start}) and (\ref{lower estimate for (I): how big is overline k for all j}) to get
\begin{align}
(I)\gtrsim_{(X)}&\frac{1}{1-\frac{1}{2^{Q+1}}}\frac{1}{2^{(Q+1)(\overline{k}(\tilde{j}(x,w)-2)+3)}}\sum_{j=-\infty}^{\tilde{j}(x,w)-2}\frac{2^{jQ}}{2^{(Qs)j}}\gtrsim_{(X,Q)}\\
&\frac{1}{2^{(Q+1)(\overline{k}(\tilde{j}(x,w)-2))}}\sum_{i=2-\tilde{j}(x,w)}^{+\infty}\frac{1}{2^{Q(1-s)i}}=\\
&\frac{1}{2^{(Q+1)(\overline{k}(\tilde{j}(x,w)-2))}}\frac{1}{1-\frac{1}{2^{Q(1-s)}}}\frac{1}{2^{Q(1-s)(2-\tilde{j}(x,w))}}\approx_{(Q,s)}\\
&\frac{1}{2^{(Q+1)(\overline{k}(\tilde{j}(x,w)-2))}}2^{Q(1-s)\tilde{j}(x,w)}.
\end{align}
So we proved the lower estimate
\begin{equation}\label{final lower estimate for (I)}
(I)\gtrsim_{(X,Q,p)}\frac{1}{2^{(Q+1)(\overline{k}(\tilde{j}(x,w)-2))}}2^{Q(1-s)\tilde{j}(x,w)}.
\end{equation}

Now we prove an upper estimate for $(I)$. We compute
\begin{align}
(I)\leq& \sum_{j=-\infty}^{\tilde{j}(x,w)}\sum_{k=\overline{k}(j)}^{+\infty}\frac{m(B_d(x,2^j))}{2^{[(Q+1)k+(Qs)j]}}\lesssim_{(X)}\\
&\sum_{j=-\infty}^{\tilde{j}(x,w)-2}\sum_{k=\overline{k}(j)}^{+\infty}\frac{2^{jQ}}{2^{[(Q+1)k+(Qs)j]}}=\nonumber\\
&\sum_{j=-\infty}^{\tilde{j}(x,w)-2}\frac{1}{1-\frac{1}{2^{Q+1}}}\frac{1}{2^{(Q+1)\overline{k}(j)}}\frac{2^{jQ}}{2^{(Qs)j}}.\nonumber
\end{align}
By monotonicity in the definition of $\overline{k}(j)$ we deduce
\begin{equation}
\overline{k}(j)\geq \overline{k}(\tilde{j}(x,w)-2) \quad \text{for all } j\leq\tilde{j}(x,w)-2,
\end{equation}
so we get
\begin{align}
(I)\lesssim_{(X)}&\frac{1}{1-\frac{1}{2^{Q+1}}}\frac{1}{2^{(Q+1)\overline{k}(\tilde{j}(x,w)-2)}}\sum_{j=-\infty}^{\tilde{j}(x,w)-2}\frac{2^{jQ}}{2^{(Qs)j}}\approx_{(X,Q)}\\
&\frac{1}{2^{(Q+1)\overline{k}(\tilde{j}(x,w)-2)}}\sum_{j=-\infty}^{\tilde{j}(x,w)-2}\frac{2^{jQ}}{2^{(Qs)j}}\approx_{(Q,s)}\nonumber\\
&\frac{1}{2^{(Q+1)\overline{k}(\tilde{j}(x,w)-2)}}2^{Q(1-s)\tilde{j}(x,w)}.\nonumber
\end{align}
So we proved
\begin{equation}
(I)\lesssim_{(X,Q,s)}\frac{1}{2^{(Q+1)\overline{k}(\tilde{j}(x,w)-2)}}2^{Q(1-s)\tilde{j}(x,w)},
\end{equation}
which, alongside (\ref{final lower estimate for (I)}), proves
\begin{equation}\label{final estimate for (I)}
(I)\approx_{(X,Q,s)}\frac{1}{2^{(Q+1)\overline{k}(\tilde{j}(x,w)-2)}}2^{Q(1-s)\tilde{j}(x,w)}.
\end{equation}
\\

Now we are going to estimate the value of $(III)$ and we are going to compare $(II)$ and $(III)$. \\
Let us define
\begin{equation}
j_{\text{diam}}:=\left\lfloor \log_2(\text{diam}(X)) \right \rfloor+1.
\end{equation}
By computation we get
\begin{equation}
2^{j_{\text{diam}}}= 2^{\left\lfloor \log_2(\text{diam}(X)) \right \rfloor+1}> 2^{\log_2(\text{diam}(X))}=\text{diam}(X),
\end{equation}
so we proved
\begin{equation}\label{estimate (III): big j gives whole space}
X=B_d(x,2^j) \quad \text{for all }j \geq j_{\text{diam}}.
\end{equation}
Moreover we get
\begin{equation}
2^{j_{\text{diam}}-1}= 2^{\left\lfloor \log_2(\text{diam}(X)) \right \rfloor}\leq 2^{\log_2(\text{diam}(X))}=\text{diam}(X).
\end{equation}
We may write
\begin{align}
(III)=&\sum_{j=\tilde{j}(x,w)+1}^{j_{\text{diam}}-1} \sum_{k=0}^{+\infty}\frac{m(B_d(x,2^j)\cap B_d(w,2^k y))}{2^{[(Q+1)k+(Qs)j]}}+\nonumber \\
&\sum_{j=j_{\text{diam}}}^{+\infty} \sum_{k=0}^{+\infty}\frac{B_d(w,2^k y)}{2^{[(Q+1)k+(Qs)j]}}
\end{align}

Consider $\tilde{j}(x,w)+1\leq j \leq j_{\text{diam}}-1$. We define
\begin{equation}\label{definition of tilde k}
\tilde{k}(j)=\max\left\{\left\lfloor \log_2\left(\frac{2^{j}-d(x,w)}{y}\right)\right\rfloor+1,0\right\},
\end{equation}
which is the smallest index $k\geq0$ such that $B_d(w,2^k y) \not \subseteq B_d(x,2^j)$.\\
Suppose $k\geq \tilde{k}(j)+2$, suppose $v \in B_d(x,2^j)$. By triangle inequality we get
\begin{equation}
d(v,w)\leq d(v,x)+d(x.w)\leq d(x,w)+2^j.
\end{equation}
However, from (\ref{definition of tilde k}) we get
\begin{equation}\label{estimate (III): triangle inequality for ball subset}
k \geq \tilde{k}(j)+2> \log_2\left(\frac{2^{j}-d(x,w)}{y}\right)+2,
\end{equation}
hence
\begin{equation}\label{estimate of (III): how big is radius of w ball}
2^k y > 4(2^j -d(x,w)).
\end{equation}
But from (\ref{definition of tilde j}) we have
\begin{equation}
j\geq \tilde{j}(x,w)+1\geq \log_2(d(x,w)) +1,
\end{equation}
hence
\begin{equation}\label{estimate of (III): how big is radius of x ball}
2^j \geq 2 d(x,w).
\end{equation}
Now from (\ref{estimate of (III): how big is radius of w ball}) and (\ref{estimate of (III): how big is radius of x ball}) we get
\begin{equation}\label{estimate (III): how big is the radius of x ball}
2^k y > 2^j+(3\cdot 2^j -4d(x,w))\geq 2^j+2 d(x,w).
\end{equation}
Combining (\ref{estimate (III): triangle inequality for ball subset}) and (\ref{estimate (III): how big is the radius of x ball}) we get
\begin{equation}
2^k y >d(v,w),
\end{equation}
so we proved that 
\begin{equation}\label{estimate of (III): end of intersection proof}
B_d(x,2^j)\subseteq B_d(w,2^k y)
\end{equation}
for every $x,w \in X$, for every $\tilde{j}(x,w)+1\leq j \leq j_{\text{diam}}-1$, for every $k\geq \tilde{k}(j)+2$.\\
Now we use (\ref{estimate of (III): end of intersection proof}) and the definition of $\tilde{k}(j)$ to write
\begin{align}
(III)=\sum_{j=\tilde{j}(x,w)+1}^{j_{\text{diam}}-1}  \bigg[&\sum_{k=0}^{\tilde{k}(j)-1}\frac{m(B_d(w,2^k y))}{2^{[(Q+1)k+(Qs)j]}}+\\
&\sum_{k=\tilde{k}(j)}^{\tilde{k}(j)+1}\frac{m(B_d(x,2^j)\cap B_d(w,2^k y))}{2^{[(Q+1)k+(Qs)j]}}+\nonumber\\
&\sum_{k=\tilde{k}(j)+2}^{+\infty}\frac{m(B_d(x,2^j))}{2^{[(Q+1)k+(Qs)j]}}\bigg]+\nonumber\\
&\sum_{j=j_{\text{diam}}}^{+\infty} \sum_{k=0}^{+\infty}\frac{m(B_d(w,2^k y))}{2^{[(Q+1)k+(Qs)j]}},\nonumber
\end{align}
with the convention $\sum_{k=a}^{b}\varphi(k):=0$ if $b<a$.\\
Let us define
\begin{equation}
k_{\text{diam}}:=\left\lfloor \log_2\left(\frac{\text{diam}(X)}{y}\right) \right \rfloor+1.
\end{equation}
By the same argument used for $j_{\text{diam}}$ we get
\begin{equation}
2^{k_{\text{diam}}} y>\text{diam}(X) \quad \text{and} \quad 2^{k_{\text{diam}}-1} y\leq  \text{diam}(X),
\end{equation}
and we prove
\begin{equation}\label{estimate (III): big k gives whole space}
X=B_d(x,2^k y) \quad \text{for all }k \geq k_{\text{diam}}.
\end{equation}
So we can write 
\begin{align}\label{estimate of (III): equation to estimate}
(III)=&\sum_{j=\tilde{j}(x,w)+1}^{j_{\text{diam}}-1}  \bigg[\sum_{k=0}^{\tilde{k}(j)-1}\frac{m(B_d(w,2^k y))}{2^{[(Q+1)k+(Qs)j]}}+\\
&\quad\quad\quad\quad\quad\sum_{k=\tilde{k}(j)}^{\tilde{k}(j)+1}\frac{m(B_d(x,2^j)\cap B_d(w,2^k y))}{2^{[(Q+1)k+(Qs)j]}}+\nonumber\\
&\quad\quad\quad\quad\;\;\sum_{k=\tilde{k}(j)+2}^{+\infty}\frac{m(B_d(x,2^j))}{2^{[(Q+1)k+(Qs)j]}}\bigg]+\nonumber\\
&\sum_{j=j_{\text{diam}}}^{+\infty} \bigg[\sum_{k=0}^{k_{\text{diam}}-1}\frac{m(B_d(w,2^k y))}{2^{[(Q+1)k+(Qs)j]}}+\nonumber\\
&\quad\quad\quad\quad\sum_{k=k_{\text{diam}}}^{+\infty}\frac{m(X)}{2^{[(Q+1)k+(Qs)j]}}\bigg].\nonumber
\end{align}
Now we are going to prove that $(II)\lesssim_{(Q,s)} (III)$.\\
Suppose  $\tilde{j}(x,w)> j_{\text{diam}}-2$.\\
Then, the expression (\ref{estimate of (III): equation to estimate}) reformulates to
\begin{align}\label{estimate of (III): better equation to estimate}
(III)=&\sum_{j=j_{\text{diam}}}^{+\infty} \bigg[\sum_{k=0}^{k_{\text{diam}}-1}\frac{m(B_d(w,2^k y))}{2^{[(Q+1)k+(Qs)j]}}+\\
&\quad\quad\quad\quad\sum_{k=k_{\text{diam}}}^{+\infty}\frac{m(X)}{2^{[(Q+1)k+(Qs)j]}}\bigg].\nonumber
\end{align}
Let us denote by $\Phi$ the first addend of the first sum in (\ref{estimate of (III): better equation to estimate}) with respect to the index $j$, i.e
\begin{equation}
\Phi= \sum_{k=0}^{k_{\text{diam}}-1}\frac{m(B_d(w,2^k y))}{2^{[(Q+1)k+(Qs)j_{\text{diam}}]}}+\sum_{k=k_{\text{diam}}}^{+\infty}\frac{m(X)}{2^{[(Q+1)k+(Qs)j_{\text{diam}}]}}.
\end{equation}
Now we consider the equation
\[(II)=\sum_{j=\tilde{j}(x,w)-1}^{\tilde{j}(x,w)} \sum_{k=0}^{+\infty}\frac{m(B_d(x,2^j)\cap B_d(w,2^k y))}{2^{[(Q+1)k+(Qs)j]}}\]
and we use monotonicity of the measure $m$ to get
\begin{equation}
(II)\leq \sum_{j=\tilde{j}(x,w)-1}^{\tilde{j}(x,w)} \sum_{k=0}^{+\infty}\frac{ m(B_d(w,2^k y))}{2^{[(Q+1)k+(Qs)j]}},
\end{equation}
which reformulates to
\begin{equation}
(II)\leq \sum_{j=\tilde{j}(x,w)-1}^{\tilde{j}(x,w)}\left[ \sum_{k=0}^{k_{\text{diam}}-1}\frac{ m(B_d(w,2^k y))}{2^{[(Q+1)k+(Qs)j]}}+\sum_{k=k_{\text{diam}}}^{+\infty}\frac{ m(X)}{2^{[(Q+1)k+(Qs)j]}}\right].
\end{equation}
By definition we get $\tilde{j}(x,w)\leq j_{\text{diam}}$, so we have two cases: 
\begin{enumerate}
\item Case $\tilde{j}(x,w)=j_{\text{diam}}-1$: we get
\begin{equation}\label{first estimate for (II) vs phi in first case}
(II)\leq 2^{2Qs} \Phi + 2^{Qs} \Phi.
\end{equation}
\item Case $\tilde{j}(x,w)=j_{\text{diam}}$: we get
\begin{equation}\label{second estimate for (II) vs phi in first case}
(II)\leq 2^{Qs} \Phi+\Phi.
\end{equation}
\end{enumerate}
So from (\ref{first estimate for (II) vs phi in first case}) and (\ref{second estimate for (II) vs phi in first case}) we get
\begin{equation}\label{first part of (II) lesssim than (III)}
(II) \lesssim_{(Q,s)} (III)
\end{equation}
whenever $\tilde{j}(x,w)> j_{\text{diam}}-2$, and the constant associated to equation (\ref{first part of (II) lesssim than (III)}) does not depend on $x$ and $w$.\\

Now suppose $\tilde{j}(x,w)\leq j_{\text{diam}}-2$
Let us denote by $\Phi$ first addend of the first sum in (\ref{estimate of (III): equation to estimate}) with respect to the index $j$, i.e.
\begin{align}
\Phi=&\sum_{k=0}^{\tilde{k}(\tilde{j}(x,w)+1)-1}\frac{m(B_d(w,2^k y))}{2^{[(Q+1)k+(Qs)(\tilde{j}(x,w)+1)]}}+\\
&\sum_{k=\tilde{k}(\tilde{j}(x,w)+1)}^{\tilde{k}(\tilde{j}(x,w)+1)+1}\frac{m(B_d(x,2^{\tilde{j}(x,w)+1})\cap B_d(w,2^k y))}{2^{[(Q+1)k+(Qs)(\tilde{j}(x,w)+1)]}}+\nonumber\\
&\sum_{k=\tilde{k}(\tilde{j}(x,w)+1)+2}^{+\infty}\frac{m(B_d(x,2^{\tilde{j}(x,w)+1}))}{2^{[(Q+1)k+(Qs)(\tilde{j}(x,w)+1)]}}.\nonumber
\end{align}
We may write
\begin{align}
(II)=\sum_{j=\tilde{j}(x,w)-1}^{\tilde{j}(x,w)}& \sum_{k=0}^{+\infty}\frac{m(B_d(x,2^j)\cap B_d(w,2^k y))}{2^{[(Q+1)k+(Qs)j]}}=\\
=\sum_{j=\tilde{j}(x,w)-1}^{\tilde{j}(x,w)} &\bigg[\sum_{k=0}^{\tilde{k}(\tilde{j}(x,w)+1)-1}\frac{m(B_d(x,2^j)\cap B_d(w,2^k y))}{2^{[(Q+1)k+(Qs)j]}}+\\
&\sum_{k=\tilde{k}(\tilde{j}(x,w)+1)}^{\tilde{k}(\tilde{j}(x,w)+1)+1}\frac{m(B_d(x,2^j))\cap B_d(w,2^k y))}{2^{[(Q+1)k+(Qs)j]}}+\nonumber\\
&\sum_{k=\tilde{k}(\tilde{j}(x,w)+1)+2}^{+\infty}\frac{m(B_d(x,2^j)\cap B_d(w,2^k y))}{2^{[(Q+1)k+(Qs)j]}}\bigg].
\end{align}
Now we use the following facts:
\begin{enumerate}
\item $m(B_d(x,2^j)\cap B_d(w,2^k y))\leq m(B_d(x,2^j))$, $m(B_d(x,2^j)\cap B_d(w,2^k y))\leq m(B_d(w,2^k y))$, by the monotonicity of $m$.
\item $m(B_d(x,2^j))\leq m(B_d(x,2^{\tilde{j}(x,w)+1}))$ for $j=\tilde{j}(x,w)-1$ and $j=\tilde{j}(x,w)$, because\\ $B_d(x,2^j)\subseteq B_d(x,2^{\tilde{j}(x,w)+1})$.
\end{enumerate}
So we get
\begin{align}
(II)\leq\sum_{j=\tilde{j}(x,w)-1}^{\tilde{j}(x,w)} &\bigg[\sum_{k=0}^{\tilde{k}(\tilde{j}(x,w)+1)-1}\frac{m(B_d(w,2^k y))}{2^{[(Q+1)k+(Qs)j]}}+\\
&\sum_{k=\tilde{k}(\tilde{j}(x,w)+1)}^{\tilde{k}(\tilde{j}(x,w)+1)+1}\frac{m(B_d(x,2^{\tilde{j}(x,w)+1})\cap B_d(w,2^k y))}{2^{[(Q+1)k+(Qs)j]}}+\nonumber\\
&\sum_{k=\tilde{k}(\tilde{j}(x,w)+1)+2}^{+\infty}\frac{m(B_d(x,2^{\tilde{j}(x,w)+1}))}{2^{[(Q+1)k+(Qs)j]}}\bigg].
\end{align}
So, by computation and using the definition of $\Phi$, we get
\begin{equation}\label{last part of (II) lesssim (III)}
(II)\leq 2^{2Qs} \Phi + 2^{Qs} \Phi.
\end{equation}
Combining (\ref{first part of (II) lesssim than (III)}) and (\ref{last part of (II) lesssim (III)}) we get
\begin{equation}\label{(II) lesssim than (III)}
(II)\lesssim_{(Q,s)} (III),
\end{equation}
for all $x,w \in X$, and the constant associated to (\ref{(II) lesssim than (III)}) does not depend on $x,w\in X$.
Now we go back to estimating equation $(III)$.\\
We consider equation (\ref{estimate of (III): equation to estimate}):
\begin{align*}
(III)=&\sum_{j=\tilde{j}(x,w)+1}^{j_{\text{diam}}-1}  \bigg[\sum_{k=0}^{\tilde{k}(j)-1}\frac{m(B_d(w,2^k y))}{2^{[(Q+1)k+(Qs)j]}}+\\
&\quad\quad\quad\quad\quad\sum_{k=\tilde{k}(j)}^{\tilde{k}(j)+1}\frac{m(B_d(x,2^j)\cap B_d(w,2^k y))}{2^{[(Q+1)k+(Qs)j]}}+\nonumber\\
&\quad\quad\quad\quad\;\;\sum_{k=\tilde{k}(j)+2}^{+\infty}\frac{m(B_d(x,2^j))}{2^{[(Q+1)k+(Qs)j]}}\bigg]+\nonumber\\
&\sum_{j=j_{\text{diam}}}^{+\infty} \bigg[\sum_{k=0}^{k_{\text{diam}}-1}\frac{m(B_d(w,2^k y))}{2^{[(Q+1)k+(Qs)j]}}+\nonumber\\
&\quad\quad\quad\quad\sum_{k=k_{\text{diam}}}^{+\infty}\frac{m(X)}{2^{[(Q+1)k+(Qs)j]}}\bigg],\nonumber
\end{align*}
and we claim that
\begin{align}\label{last claim in estimate of (III)}
\sum_{k=\tilde{k}(j)}^{\tilde{k}(j)+1}\frac{m(B_d(x,2^j)\cap B_d(w,2^k y))}{2^{[(Q+1)k+(Qs)j]}}+\sum_{k=\tilde{k}(j)+2}^{+\infty}\frac{m(B_d(x,2^j))}{2^{[(Q+1)k+(Qs)j]}}\approx_{(Q)} \sum_{k=\tilde{k}(j)+2}^{+\infty}\frac{m(B_d(x,2^j))}{2^{[(Q+1)k+(Qs)j]}}.
\end{align}
Indeed, we can trivially get one part of equation (\ref{last claim in estimate of (III)}) because
\begin{equation}
\sum_{k=\tilde{k}(j)+2}^{+\infty}\frac{m(B_d(x,2^j))}{2^{[(Q+1)k+(Qs)j]}}\leq \sum_{k=\tilde{k}(j)}^{\tilde{k}(j)+1}\frac{m(B_d(x,2^j)\cap B_d(w,2^k y))}{2^{[(Q+1)k+(Qs)j]}}+\sum_{k=\tilde{k}(j)+2}^{+\infty}\frac{m(B_d(x,2^j))}{2^{[(Q+1)k+(Qs)j]}}.
\end{equation}
For the other part we compute
\begin{align}\label{first equation for last comparison in (III)}
\sum_{k=\tilde{k}(j)+2}^{+\infty}\frac{m(B_d(x,2^j))}{2^{[(Q+1)k+(Qs)j]}}=\frac{m(B_d(x,2^j))}{2^{(Qs)j}}\sum_{k=\tilde{k}(j)+2}^{+\infty}\frac{1}{2^{(Q+1)k}}=\frac{m(B_d(x,2^j))}{2^{[(Qs)j+(Q+1)(\tilde{k}(j)+2)]}}\frac{1}{1-2^{-(Q+1)}}
\end{align}
and
\begin{align}\label{second equation for last comparison in (III)}
\sum_{k=\tilde{k}(j)}^{\tilde{k}(j)+1}\frac{m(B_d(x,2^j)\cap B_d(w,2^k y))}{2^{[(Q+1)k+(Qs)j]}}\leq&\sum_{k=\tilde{k}(j)}^{\tilde{k}(j)+1}\frac{m(B_d(x,2^j))}{2^{[(Q+1)k+(Qs)j]}}=\\
&\frac{m(B_d(x,2^j))}{2^{[(Qs)j+(Q+1)(\tilde{k}(j)+2)]}}\left[2^{Q+1}+2^{2(Q+1)} \right].\nonumber
\end{align}
From (\ref{first equation for last comparison in (III)}) and (\ref{second equation for last comparison in (III)}) we get
\begin{equation}
\sum_{k=\tilde{k}(j)}^{\tilde{k}(j)+1}\frac{m(B_d(x,2^j)\cap B_d(w,2^k y))}{2^{[(Q+1)k+(Qs)j]}}\leq\frac{1}{(2^{Q+1}+2^{2(Q+1)} )(1-2^{-(Q+1)})}\sum_{k=\tilde{k}(j)+2}^{+\infty}\frac{m(B_d(x,2^j))}{2^{[(Q+1)k+(Qs)j]}},
\end{equation}
so we proved
\begin{equation}
\sum_{k=\tilde{k}(j)}^{\tilde{k}(j)+1}\frac{m(B_d(x,2^j)\cap B_d(w,2^k y))}{2^{[(Q+1)k+(Qs)j]}}\lesssim_{(Q)}\sum_{k=\tilde{k}(j)+2}^{+\infty}\frac{m(B_d(x,2^j))}{2^{[(Q+1)k+(Qs)j]}},
\end{equation}
and the claim (\ref{last claim in estimate of (III)}) follows.\\

Finally, from (\ref{estimate of (III): equation to estimate}) and  (\ref{last claim in estimate of (III)}) we get
\begin{align*}
(III)\approx_{(Q)}&\sum_{j=\tilde{j}(x,w)+1}^{j_{\text{diam}}-1}  \bigg[\sum_{k=0}^{\tilde{k}(j)-1}\frac{m(B_d(w,2^k y))}{2^{[(Q+1)k+(Qs)j]}}+\\
&\quad\quad\quad\quad\;\;\sum_{k=\tilde{k}(j)+2}^{+\infty}\frac{m(B_d(x,2^j))}{2^{[(Q+1)k+(Qs)j]}}\bigg]+\nonumber\\
&\sum_{j=j_{\text{diam}}}^{+\infty} \bigg[\sum_{k=0}^{k_{\text{diam}}-1}\frac{m(B_d(w,2^k y))}{2^{[(Q+1)k+(Qs)j]}}+\nonumber\\
&\quad\quad\quad\quad\sum_{k=k_{\text{diam}}}^{+\infty}\frac{m(X)}{2^{[(Q+1)k+(Qs)j]}}\bigg].\nonumber
\end{align*}
By construction we have $2^k y \leq \text{diam}(X)$ for all $\tilde{j}(x,w)+1\leq j \leq j_{\text{diam}}-1$, for all  $0\leq k \leq \tilde{k}(j)-1$, $2^j \leq \text{diam}(X)$ for all $\tilde{j}(x,w)+1\leq j \leq j_{\text{diam}}-1$, $2^k y \leq \text{diam}(X)$ for all $k\leq k_{\text{diam}}-1$, so we may apply the Ahlfors-regularity estimate (\ref{third of 3 equations to estimate measure of balls}) to get
\begin{align}\label{final estimate for (III)}
(III)\approx_{(X,Q)}&\sum_{j=\tilde{j}(x,w)+1}^{j_{\text{diam}}-1}  \bigg[\sum_{k=0}^{\tilde{k}(j)-1}\frac{2^{Qk} y^Q}{2^{[(Q+1)k+(Qs)j]}}+\\
&\quad\quad\quad\quad\;\;\sum_{k=\tilde{k}(j)+2}^{+\infty}\frac{2^j}{2^{[(Q+1)k+(Qs)j]}}\bigg]+\nonumber\\
&\sum_{j=j_{\text{diam}}}^{+\infty} \bigg[\sum_{k=0}^{k_{\text{diam}}-1}\frac{2^{Qk} y^Q}{2^{[(Q+1)k+(Qs)j]}}+\nonumber\\
&\quad\quad\quad\quad\sum_{k=k_{\text{diam}}}^{+\infty}\frac{m(X)}{2^{[(Q+1)k+(Qs)j]}}\bigg].\nonumber
\end{align}

Now we recall what we proved up this point and finish the proof of the statement.\\
We proved (\ref{final estimate for (I)}), (\ref{(II) lesssim than (III)}) and (\ref{final estimate for (III)}), i.e.
\begin{equation*}
(I)\approx_{(X,Q,s)}\varphi_1(x,w,y):=\frac{1}{2^{(Q+1)\overline{k}(\tilde{j}(x,w)-2)}}2^{Q(1-s)\tilde{j}(x,w)},
\end{equation*}
\begin{equation*}
0\leq (II) \lesssim_{(Q,s)} (III),
\end{equation*}
\begin{align*}
(III)\approx_{(X,Q)}\varphi_{2}(x,w,y):=&\sum_{j=\tilde{j}(x,w)+1}^{j_{\text{diam}}-1}  \bigg[\sum_{k=0}^{\tilde{k}(j)-1}\frac{2^{Qk} y^Q}{2^{[(Q+1)k+(Qs)j]}}+\\
&\quad\quad\quad\quad\;\;\sum_{k=\tilde{k}(j)+2}^{+\infty}\frac{2^j}{2^{[(Q+1)k+(Qs)j]}}\bigg]+\nonumber\\
&\sum_{j=j_{\text{diam}}}^{+\infty} \bigg[\sum_{k=0}^{k_{\text{diam}}-1}\frac{2^{Qk} y^Q}{2^{[(Q+1)k+(Qs)j]}}+\nonumber\\
&\quad\quad\quad\quad\sum_{k=k_{\text{diam}}}^{+\infty}\frac{m(X)}{2^{[(Q+1)k+(Qs)j]}}\bigg].\nonumber
\end{align*}
We recall that, by definition, we have
\begin{equation}
\sum_{j=-\infty}^{+\infty} \sum_{k=0}^{+\infty}\frac{m(B_d(x,2^j)\cap B_d(w,2^k y))}{2^{[(Q+1)k+(Qs)j]}}=(I)+(II)+(III),
\end{equation}
so we get
\begin{equation}\label{final equation in exchange x and w proof}
\sum_{j=-\infty}^{+\infty} \sum_{k=0}^{+\infty}\frac{m(B_d(x,2^j)\cap B_d(w,2^k y))}{2^{[(Q+1)k+(Qs)j]}}\approx_{(Q,s)} (I)+(III) \approx_{(X,Q,s)} \varphi_{1}(x,w,y)+\varphi_{2}(x,w,y).
\end{equation}
However, $\varphi_{1}(x,w,y)=\varphi_{1}(w,x,y)$, $\varphi_{2}(x,w,y)=\varphi_{2}(w,x,y)$ for all $x,w \in X$, for all $y>0$.\\
Indeed, we have
\begin{equation}
\tilde{j}(x,w)=\left\lceil \log_2(d(x,w)) \right\rceil=\left\lceil \log_2(d(w,x)) \right\rceil=\tilde{j}(w,x),
\end{equation}
\begin{equation}
\overline{k}(j)=\max\left\{\left\lfloor \log_2\left(\frac{d(x,w)-2^j}{y}\right)\right\rfloor+1,0\right\}=\max\left\{\left\lfloor \log_2\left(\frac{d(w,x)-2^j}{y}\right)\right\rfloor+1,0\right\},
\end{equation}
\begin{equation}
\tilde{k}(j)=\max\left\{\left\lfloor \log_2\left(\frac{2^{j}-d(x,w)}{y}\right)\right\rfloor+1,0\right\}=\max\left\{\left\lfloor \log_2\left(\frac{2^{j}-d(w,x)}{y}\right)\right\rfloor+1,0\right\},
\end{equation}
and all the other terms in the definitions of $\varphi_1$ and $\varphi_2$ do not depend on $x$ and $w$, so we get that $\varphi_1, \varphi_2$ are symmetrical with respect to exchanging the roles of $x$ and $w$.\\
We finish the proof by exchanging the roles of $x$ and $w$ in (\ref{final equation in exchange x and w proof}) and we get
\begin{align}
\sum_{j=-\infty}^{+\infty} \sum_{k=0}^{+\infty}\frac{m(B_d(w,2^j)\cap B_d(x,2^k y))}{2^{[(Q+1)k+(Qs)j]}} \approx_{(X,Q,s)}& \varphi_{1}(w,x,y)+\varphi_{2}(w,x,y)=\\
& \varphi_{1}(x,w,y)+\varphi_{2}(x,w,y)\approx_{(X,Q,s)}\nonumber\\
&\sum_{j=-\infty}^{+\infty} \sum_{k=0}^{+\infty}\frac{m(B_d(x,2^j)\cap B_d(w,2^k y))}{2^{[(Q+1)k+(Qs)j]}},\nonumber
\end{align}
which is the statement (\ref{literally final thing to prove}) we needed to prove, so the proof is finished.
\endproof
\subsection{Other properties of the Poisson Integral}
Now we will prove two more properties of the Poisson Integral which are analogous to the properties of the classical Poisson Integral in $\mathbb R^{n+1}$.
\begin{Def}\label{E, E^* and E'}
Let $f$ be a non negative function in $L^p_+(X)$, $\epsilon>0$. We define
\begin{equation}
E(f,\epsilon):=\{(x,y) \in X \times (0,+\infty) \; | \; PI(K_{X,s} * f)(x,y)>\epsilon\},
\end{equation}
\begin{equation}
E^*(f,\epsilon):= \bigcup_{(x,y) \in E(f,\epsilon)} B_d(x,y)\subseteq X,
\end{equation}
\begin{equation}
E'(f,\epsilon):= \bigcup_{(x,y) \in E(f,\epsilon)} B_d(x,y)\times\{y\} \subseteq X \times (0,+\infty).
\end{equation}
\end{Def}
\begin{Lem}\label{Harnack type inequality}[Harnack-type inequality]
We have
\begin{equation}
PI(K_{X,s} * f)(x,y)\gtrsim_{(X,Q)} \epsilon \quad \forall (x,y) \in E'(f,\epsilon).
\end{equation}
\end{Lem}
\proof
Let $(x,y) \in E'(f,\epsilon)$. By definition of $E'$ there exist $(\tilde x,y) \in E(f,\epsilon)$ such that $x \in B_d(\tilde x,y)$. By (\ref{boundedness of constant in Poisson Integral})  the constant $C$ in the Poisson integral
\begin{equation}
PI(K_{X,s} * f)(x,y):= \int_{X} C(x,y) \cdot \frac{1}{y^Q} \sum_{k=0}^{+\infty} \frac{\chi_{B_d(x,2^k y)}(z)}{2^{(Q+1)k}}K_{X,s} * f(z)dm(z)
\end{equation}
 is bounded, so we have $C(x,y) \approx_{(X,Q)} C(\tilde x,y)$, and we get
\begin{equation}
PI(K_{X,s} * f)(x,y)\gtrsim_{(X,Q)} \int_{X} C(\tilde x,y) \cdot \frac{1}{y^Q} \sum_{k=0}^{+\infty} \frac{\chi_{B_d(x,2^k y)}(z)}{2^{(Q+1)k}}K_{X,s} * f(z)dm(z).
\end{equation}
The function $f$ is non negative, so $K_{X,s} * f$ is also non negative, so we have
\begin{equation}
PI(K_{X,s} * f)(x,y)\gtrsim_{(X,Q)} \int_{X} C(\tilde x,y) \cdot \frac{1}{y^Q} \sum_{k=1}^{+\infty}\frac{1}{2^{(Q+1)}} \frac{\chi_{B_d(x,2^k y)}(z)}{2^{(Q+1)(k-1)}}K_{X,s} * f(z)dm(z).
\end{equation}
Now we change the index $k$ of summation and put $1/2^{(Q+1)}$ in the leading constant to get
\begin{equation}\label{second to last step in proof of Harnack inequality property}
PI(K_{X,s} * f)(x,y)\gtrsim_{(X,Q)} \int_{X} C(\tilde x,y) \cdot \frac{1}{y^Q} \sum_{k=0}^{+\infty} \frac{\chi_{B_d(x,2^{k+1} y)}(z)}{2^{(Q+1)k}}K_{X,s} * f(z)dm(z).
\end{equation}
Now let $k \in \{0,1,2,...\}$ be a fixed index. We have $d(x,\tilde x)<y$ by definition of $\tilde x$, so, for all $z \in B_d(\tilde x,2^k y)$, we apply the triangle inequality to get
\begin{equation}
d(x,z)\leq d(\tilde x,x)+d(\tilde x,z)< y+ 2^k y \leq (2^k+1)y \leq 2^{k+1}y.
\end{equation}
So we proved that $B_d(\tilde x,2^k y) \subseteq B_d(x,2^{k+1} y)$ for all $k \in \mathbb N$, hence
\begin{equation}\label{bigger sum of charateristic functions}
 \sum_{k=0}^{+\infty} \frac{\chi_{B_d(x,2^{k+1} y)}(z)}{2^{(Q+1)k}}\geq  \sum_{k=0}^{+\infty} \frac{\chi_{B_d(\tilde x,2^k y)}(z)}{2^{(Q+1)k}}.
\end{equation}
We combine (\ref{second to last step in proof of Harnack inequality property}) and (\ref{bigger sum of charateristic functions}) to get
\begin{equation}
PI(K_{X,s} * f)(x,y)\gtrsim_{(X,Q)} \int_{X} C(\tilde x,y) \cdot \frac{1}{y^Q}\sum_{k=0}^{+\infty} \frac{\chi_{B_d(\tilde x,2^k y)}(z)}{2^{(Q+1)k}}K_{X,s} * f(z)dm(z).
\end{equation}
However, by the definition of Poisson Integral and of $\tilde x$ and of $E(f,\epsilon)$ we get
\begin{equation}
PI(K_{X,s} * f)(x,y)\gtrsim_{(X,Q,M)} PI(K_{X,s} * f)(\tilde x,y)\geq \epsilon,
\end{equation}
which is the required inequality.
\endproof
\begin{Oss}
The previous Harnack-type inequality still holds if we replace $K_{X,s}*f$ in the previous definition and lemma with a generic function $g \in L^1(X)$ such that $g>0$.
\end{Oss}
\begin{Lem}[Uniform continuity at the boundary of the Poisson Integral]\label{Uniform continuity of the Riesz potential}
Let $g\in C(X)$. Then for every $\epsilon>0$ there exists $\delta>0$ such that 
\begin{equation}\label{thesis uniform continuity of potential}
\sup_{x_0 \in X} \left( \sup_{P \in B_{\rho}((x_0,0),\delta)}|PI(g)(P)-g(x_0)|\right)\leq \epsilon
\end{equation}
\end{Lem}
\proof
Let $g \in C_0(X)$. Then $g \in C(x)$. The space $X$ is compact, so by Heine-Cantor theorem we have that for all $\epsilon_1>0$ there exists $\delta_1=\delta_1(\epsilon_1)>0$ such that for every $x_1,x_2 \in X$ such that $d(x_1,x_2)<\delta_1$ we have $|g(x_1)-g(x_2)|<\epsilon_1$.\\
Let $\epsilon_1>0$ be a number to be fixed later. Let $\delta_1(\epsilon_1)>0$ be the number defined by the Heine-Cantor theorem. We claim that for every $0<\epsilon_2<1$ there exist $\tilde y= \tilde y(\epsilon_2,\delta_1)$ such that 
\begin{equation}\label{claim that poisson kernel is concentrated}
I(x,y,\delta_1):=\int_{X\backslash B_d\left(x,\frac{\delta_1}{2}\right)} C(x,y) \cdot \frac{1}{y^Q} \sum_{k=0}^{+\infty} \frac{\chi_{B_d(x,2^k y)}(z)}{2^{(Q+1)k}}dm(z)\leq \epsilon_2
\end{equation}
for all $x \in X$, for all $y< \tilde y$.\\
Indeed, let $\epsilon_2>0$. For all $k\leq \left\lfloor \log_2\left(\frac{\delta_1}{2y}\right)\right\rfloor$ we have $B_d(x,2^k y) \subseteq B_d(x,\frac{\delta_1}{2})$, so we can estimate
\begin{align}
\int_{X\backslash B_d(x,\frac{\delta_1}{2})} C(x,y) \cdot \frac{1}{y^Q} \sum_{k=0}^{+\infty} \frac{\chi_{B_d\left(x,2^k y\right)}(z)}{2^{(Q+1)k}}dm(z)\leq C_2\sum_{k= \left\lfloor \log_2\left(\frac{\delta_1}{2y}\right)\right\rfloor+1}^{+\infty} \frac{m(B_d(x,2^k y))}{y^Q 2^{(Q+1)k}},
\end{align}
here $C_2$ is the uniform upper estimate of $C(x,y)$ given in (\ref{boundedness of constant in Poisson Integral}).\\
By Ahlfors-regularity of $X$ we get the estimate
\begin{equation}
I(x,y,\delta_1)\lesssim_{(X)}  C_2\sum_{k= \left\lfloor \log_2\left(\frac{\delta_1}{2y}\right)\right\rfloor+1}^{+\infty} \frac{y^Q 2^{kQ}}{y^Q 2^{(Q+1)k}}=C_2 \frac{2}{2^{\left(\left\lfloor \log_2\left(\frac{\delta_1}{2y}\right)\right\rfloor+1\right)}}:=\varphi(\delta_1,y).
\end{equation}
We observe that $\varphi(\delta_1,y)\rightarrow 0$ as $y\rightarrow 0$, so there exists $\tilde y=\tilde y(\delta_1,\epsilon_2)$ such that
\begin{equation}
I(x,y,\delta_1)\lesssim_{(X)} \varphi(\delta_1,y) \leq \epsilon_2
\end{equation}
for all $x \in X$, for all $0<y<\tilde y$.
Up to multiplying $\epsilon_2$ by a constant which depends only on $X$ we get
\begin{equation}
I(x,y,\delta_1)\leq \epsilon_2
\end{equation}
for all $x \in X$, for all $0<y<\tilde y$, proving the claim (\ref{claim that poisson kernel is concentrated}).\\
Now we are going to prove the statement. Let $\epsilon>0$. Fix
\begin{equation}
\epsilon_1=\epsilon_1(\epsilon):=\frac{\epsilon}{2}.
\end{equation}
Let $\delta_1=\delta_1(\epsilon_1)$ be the number defined by Heine-Cantor theorem. Let
\begin{equation}
\text{osc}(g):=\sup(g)-\inf(g)
\end{equation}
be the oscillation of the function $g$. If $\text{osc}(g)=0$ then the functions $g$ and $PI(g)$ are constant, so the claim is trivial. Suppose $\text{osc}(g)>0$. \\
The function $g$ is continuous over a compact set $X$, so $\text{osc}(g)<+\infty$. Fix
\begin{equation}
\epsilon_2=\epsilon_2(\epsilon):=\frac{\epsilon}{2 \cdot \text{osc}(g)}.
\end{equation}
Let $\tilde y=\tilde y(\epsilon_2,\delta_1)$ be the number previously defined.\\
Let $x_0 \in X$, let $P=(x,y) \in X \times (0,+\infty)$ such that $y<\tilde y$ and $d(x,x_0)<\frac{\delta_1}{2}$. We compute
\begin{align}
PI(g)(P)-g(x_0)=&\int_{X}  \frac{C(x,y)}{y^Q} \sum_{k=0}^{+\infty} \frac{\chi_{B_d(x,2^k y)}(z)}{2^{(Q+1)k}}g(z)dm(z)-g(x_0)=\\
&\int_{X}  \frac{C(x,y)}{y^Q} \sum_{k=0}^{+\infty} \frac{\chi_{B_d(x,2^k y)}(z)}{2^{(Q+1)k}}\left(g(z)- g(x_0)\right)dm(z).\nonumber
\end{align}
Now we write
\begin{align}
|PI(g)(P)-g(x_0)|\leq&\underset{B_d\left(x,\frac{\delta_1}{2}\right)}{\int}  \frac{C(x,y)}{y^Q} \sum_{k=0}^{+\infty} \frac{\chi_{B_d(x,2^k y)}(z)}{2^{(Q+1)k}}\left|g(z)- g(x_0)\right|dm(z)+\\
&\underset{X\backslash B_d\left(x,\frac{\delta_1}{2}\right)}{\int}  \frac{C(x,y)}{y^Q} \sum_{k=0}^{+\infty} \frac{\chi_{B_d(x,2^k y)}(z)}{2^{(Q+1)k}}\left|g(z)- g(x_0)\right|dm(z).\nonumber
\end{align}
If $z \in B_d\left(x,\frac{\delta_1}{2}\right)$ then, by triangle inequality, we have 
\begin{equation}
d(x_0,z)\leq d(x_0,x)+d(x,z)\leq \frac{\delta_1}{2}+\frac{\delta_1}{2}=\delta_1,
\end{equation}
so by the definition of $\delta_1$ we have $|g(z)-g(x_0)|<\epsilon_1$, so we get the estimate
\begin{align}
|PI(g)(P)-g(x_0)|\leq&\underset{B_d\left(x,\frac{\delta_1}{2}\right)}{\int} C(x,y) \cdot \frac{1}{y^Q} \sum_{k=0}^{+\infty} \frac{\chi_{B_d(x,2^k y)}(z)}{2^{(Q+1)k}}\cdot \epsilon_1 dm(z)+\\
&\underset{X\backslash B_d\left(x,\frac{\delta_1}{2}\right)}{\int} C(x,y) \cdot \frac{1}{y^Q} \sum_{k=0}^{+\infty} \frac{\chi_{B_d(x,2^k y)}(z)}{2^{(Q+1)k}}\text{osc}(g)dm(z).\nonumber
\end{align}
However, $y<\tilde y$, so we get
\begin{align}
|PI(g)(P)-g(x_0)|\leq 1 \cdot \epsilon_1 +\text{osc}(g)\cdot \epsilon_2.
\end{align}
So we substitute the values of $\epsilon_1$ and $\epsilon_2$ and we get
\begin{equation}
|PI(g)(P)-g(x_0)|\leq \epsilon
\end{equation}
for all $P=(x,y) \in X \times (0,+\infty)$ such that $y<\tilde y$, $d(x,x_0)<\frac{\delta_1}{2}$. By defining
\begin{equation}
\delta=\delta(\epsilon):=\min\left\{\frac{\delta_1}{2},\tilde y\right\}
\end{equation}
we have that if $\rho((x,y),(x_0,0))<\delta$ then
\begin{equation}
d(x,x_0)<\delta\leq \frac{\delta_1}{2}, \quad \text{and } \quad |y|<\delta\leq \tilde y.
\end{equation}
So we proved that for all $\epsilon>0$ there exists $\delta=\delta(\epsilon)$ such that
\begin{equation}
|PI(g)(P)-g(x_0)|\leq \epsilon
\end{equation}
for all $x_0 \in X$, for all $P \in B_{\rho}((x_0,0),\delta)$.\\
Taking the supremum over all $x_0$ and $P$ gives us (\ref{thesis uniform continuity of potential}), ending the proof.
\endproof

\section{Convergence at the boundary}
In this section we prove several technical lemmas and propositions, and then we prove the two main results of this work: the non tangential convergence at the boundary of the harmonic extension of a Riesz potential up to an exceptional set of zero capacity and the tangential convergence at the boundary of the harmonic extension of a Riesz potential up to an exceptional set of null measure.\\\\
\subsection{$C_{K_{X,s},p}$-thinness at the boundary}
Let $(X,d,m)$ be an Ahlfors $Q$-regular space. The following definitions generalize the concept of zero capacity to the space $X \times (0,+\infty)$, and allow us to formulate the main result of this work.
\begin{Def}\label{E E_t and E^{M,*} definition}
Let $E \subseteq X \times (0,+\infty)$.  We define
\begin{equation}
E_t:=\{(x,y) \in E \; | \; 0<y<t\},
\end{equation}
\begin{equation}\label{definition of E^{M,*}}
E^*:= \bigcup_{(x,y) \in E} B_d(x,y) \subseteq X,
\end{equation}
\begin{equation}
E^*_t := \bigcup_{\underset{0<y<t}{(x,y) \in E}}B_d(x,y) \subseteq X.
\end{equation}
\end{Def}
\begin{Def}\label{C-tinness definition}
Let $E \subseteq X \times(0,+\infty)$. $E$ is $C_{K_{X,s},p}$-thin at $X \times\{0\}$ if
\begin{equation}
\lim_{t \rightarrow 0} C_{K_{X,s},p}(E^*_t)=0.
\end{equation}
\end{Def}
\begin{Rem}
If $E$ is $C_{K_{X,s},p}$-thin at  $X \times\{0\}$ then the essential projection of E
\begin{equation*}
\{x \in X \; | \; \forall t >0 \; \exists y<t \; \text{ such that } (x,y) \in E\}
\end{equation*}
is of $C_{K_{X,s},p}$-capacity 0, and hence of measure 0.
\end{Rem}

For every $x \in X$ and $r>0$ let us define (when it exists) the radius
\begin{equation}
\eta_{X,p}(x,r):= \inf \left\{ R>0 \; | \; m(B_d(x,R))\geq C_{K_{X,s},p}(B_d(x,r))\right\},
\end{equation}
and let us define
\begin{equation}
\eta^*_{X,p}(x,r):=\max \{ r,\eta_{X,p}(x,r)\}.
\end{equation}
Let $C\geq 1$. Let $E \subseteq X$. We define
\begin{equation}
\tilde{E}_{K,p,C}:=\bigcup_{x \in E} B_d(x,C \cdot \eta^*_{X,p}(x,\delta_E(x))),
\end{equation}
where
\begin{equation}
\delta_E(x):=d(x,X \backslash E).
\end{equation}
The following lemmas and propositions will be used to prove the main results of this work.
\begin{Prop}\label{Theorem 6}
Let $E\subseteq X$ be a Borel set. Under the previous notations we have
\begin{equation}
m(\tilde{E}_{K,p,C}) \lesssim_{(X,Q,s,p,C)}  C_{K_{X,s},p}(E),
\end{equation}
for all constants $C\geq\Omega$, where $\Omega$ is the constant defined by Theorem \ref{Quasi additivity Ahlfors space theorem}.
\end{Prop}
\proof
Let $\Psi \geq 1$ be the constant defined in the proof of Theorem \ref{Quasi additivity Ahlfors space theorem}. Let $C\geq \Psi$. Let $F$ be an arbitrary compact subset of $\tilde{E}_{K,p,C}$. We claim that we can find a finite family of points $x_j \in E$ such that
\begin{align*}
&F \subset \bigcup_{j} B_d(x_j,5C \cdot \eta^*_{X,p}(x_j, r_j)),\\
&\{B(x_j, \Psi\cdot  \eta^*_{X,p}(x_j, r_j))\} \; \text{ is disjoint,}\\
&r_j=\delta_E(x_j).
\end{align*}
Indeed, we consider the open covering
\begin{equation}
 \bigcup_{x \in E} B_d(x,C \cdot \eta^*_{X,p}(x, \delta_E(x)))=\tilde{E}_{K,p,C} \supseteq F,
\end{equation}
and by compactness of $F$ we find a finite covering
\begin{equation}\label{finite covering in covering lemma of F}
\bigcup_{j=1}^{N} B_d(x_j,C \cdot \eta^*_{X,p}(x,r_j))\supseteq F.
\end{equation}
It is not restrictive to assume that
\begin{equation}
r_1\geq r_2 \geq \dots \geq r_N.
\end{equation}
We can find a finite covering 
\begin{equation}
\bigcup_{k=1}^{\tilde N} B_d(x_{j_{k}},5C \cdot \eta^*_{X,p}(x,r_{j_{k}}))\supseteq F,
\end{equation}
such that
\begin{equation}
\left\{ B_d(x_{j_{k}},C \cdot \eta^*_{X,p}(x, r_{j_{k}})) \right\}_{k=1}^{\tilde N} \quad \text{is disjoint.}
\end{equation}
Indeed, consider $r_1$, which is the greatest radius $r_j$ for $j=1,2,\dots,N$. Suppose there are exactly $N(1)$ indices $j_{k_1}, \dots, j_{k_{N(1)}} \neq 1$ such that
\begin{equation}
 B_d(x_{j_{k}},C \cdot \eta^*_{X,p}(x, r_{j_{k}})) \cap  B_d(x_{1},C \cdot \eta^*_{X,p}(x, r_{1}))\neq \emptyset,
\end{equation}
for some $0\leq N(1)\leq N-1$. Then we have $r_{j_{h}}\leq r_{1}$ for $h=1,2,\dots,N(1)$. By triangle inequality we get
\begin{equation}\label{intersetction between big B_1 and other balls in covering lemma}
B_d(x_{j_h},C \cdot \eta^*_{X,p}(x, r_{j_h}))\subseteq B_d(x_{1},5C \cdot \eta^*_{X,p}(x, r_{1}))
\end{equation}
for all $h=1,2,\dots,N(1)$. So, from (\ref{finite covering in covering lemma of F}) and (\ref{intersetction between big B_1 and other balls in covering lemma}), we get
\begin{equation}
F \subseteq  B_d(x_{1},5C \cdot \eta^*_{X,p}(x, r_{1}))\cup \bigcup_{\underset{j \not \in \left\{j_{k_1},\dots,j_{k_{N(1)}}\right\}}{j=2,\dots,N}} B_d(x_{j},C \cdot \eta^*_{X,p}(x, r_{j})),
\end{equation}
and we have
\begin{equation}
B_d(x_{1},C \cdot \eta^*_{X,p}(x,r_{1}))\cap B_d(x_{j},C \cdot \eta^*_{X,p}(x,r_{j}))=\emptyset
\end{equation}
for all $j \not \in \{1,j_{k_1},j_{k_2},\dots, j_{k_{(N_1)}}\}$.\\
So we iterate this procedure a finite amout of times, considering each time $r_j$ the greatest radius in the family $\{r_{j_1}, \dots, r_{j_{M}}\}$, and we prove that there exists a family of indices $\{\tilde j_1, \dots, \tilde j_{\tilde N}\}$ such that
\begin{equation}
F \subseteq \bigcup_{k=1,\dots,\tilde N} B_d(x_{\tilde j_{k}},5C \cdot \eta^*_{X,p}(x,r_{\tilde j_{k}})),
\end{equation}
and
\begin{equation}
\left\{B_d(x_{\tilde j_{k}},C \cdot \eta^*_{X,p}(x, r_{j}))\right\}_{k=1,\dots,\tilde N} \quad \text{is disjoint.}
\end{equation}
The claim follows because $C \geq \Psi$.\\
Let $E'=\bigcup_j B_d(x_j,r_j)$. By definition of $\delta_E$ this is a subset of $E$. We apply Theorem \ref{Quasi additivity Ahlfors space theorem} for $B(x_j,r_j)$ and $E'$ and we get
\begin{equation}
\sum_{j} C_{K_{X,s},p}(B_d(x_j,r_j))  \lesssim_{(X,Q,s,p)} C_{K_{X,s},p}(E') \leq C_{K_{X,s},p}(E).
\end{equation}
We observe that we apply Theorem $\ref{Quasi additivity Ahlfors space theorem}$ instead of the finite quasi-additivity formula because the number of sets in the family $\{B_d(x_j,r_j)\}_{j}$ depends on the choice of the sets $F$ and $E$, and it can be arbitrarily large.\\
Now we observe that
\begin{equation}
m(F) \leq \sum_j m(B_d(x_j,5C \cdot\eta^*_{X,p}(x_j,r_j)))\approx_{(X,Q,s,p,C)} \sum_j  m(B_d(x_j,\eta^*_{X,p}(x_j,r_j))).
\end{equation}
By definition of $\eta^*_{X,p}$, using properties of the Riesz capacity (see Proposition \ref{capacity of ball general formula proposition}) and the compactness of $X$, it follows that
\begin{equation}
m( B_d(x_j, \eta^*_{X,p}(x_j,r_j)))\approx_{(X,Q,s)} m( B_d(x_j, \eta_{X,p}(x_j,r_j))),
\end{equation}
and by definition of $\eta_{X,p}$ and Ahlfors-regularity we have
\begin{equation}
m( B_d(x_j, \eta_{X,p}(x_j,r_j)))\approx_{(X,Q,s)} C_{K_{X,s},p}(B_d(x_j,r_j)).
\end{equation}
Hence we get
\begin{equation}
m(F)\lesssim_{(X,Q,s,p,C)} C_{K_{X,s},p}(E).
\end{equation}
A measure $m$ on an Ahlfors-regular space $(X,d,m)$ is regular, so, since $F$ is an arbitrary compact subset of $\tilde{E}_{K,p,C}$, the required inequality follows and the proposition is proved.
\endproof

\begin{Lem}\label{Lemma 3}
Let $f \in L^p(X)$, $f\geq 0$. Let $\epsilon >0$. Consider
\begin{equation}
E(f,\epsilon)=\{(x,y) \in X \times (0,+\infty) \; | \; PI(K_{X,s}*f)(x,y)>\epsilon\},
\end{equation}
\begin{equation}
E^*(f,\epsilon)=\bigcup_{(x,y) \in E(f,\epsilon)} B_d(x,y).
\end{equation}
Then
\begin{equation}
C_{K_{X,s},p}(E^*)(f,\epsilon))\lesssim_{(X,Q,s,p)} \left(\frac{\| f \|_{L^p(X)}}{\epsilon}\right)^p.
\end{equation}
\end{Lem}
\proof
Let
\begin{equation}
E'(f,\epsilon)=\bigcup_{(x,y) \in E(f,\epsilon)} B_d(x,y)\times\{y\} \subseteq X \times(0,+\infty).
\end{equation}
Since $f\geq0$ we have $PI(K_{X,s}*f)\geq0$, so we may apply Lemma \ref{Harnack type inequality} to get
\begin{equation}\label{Harnack consequence}
PI(K_{X,s}*f)(x,y)\gtrsim_{(X,Q)} \epsilon \quad \forall (x,y) \in E'(f,\epsilon).
\end{equation}
Let us consider the maximal function
\begin{equation}
F(x):= \sup_{y>0} PI(f)(x,y).
\end{equation}
By the maximal inequality (see \cite[Theorem 3.7]{SW}) we get
\begin{equation}\label{maximal inequality consequence}
\|F\|_{L^p(X)} \lesssim_{(X,p)} \|f\|_{L^p(X)},
\end{equation}
so $F \in L^p(X)$.\\
Now we compute the potential of $F$ and we get that, for every $y>0$, we have
\begin{align}
K_{X,s}*F(x)=&\int_X K_{X,s}(x,z)F(z)dm(z)=\\
&\int_X  K_{X,s}(x,z) \sup_{\tilde y>0}PI(f)(z,\tilde y)dm(z)\geq \nonumber\\
&\int_X K_{X,s}(x,z) PI(f)(z,y)dm(z)= \nonumber\\
& K_{X,s}* (PI(f)(\cdot,y))(x).\nonumber
\end{align}
Now we apply Lemma \ref{Lemma exchange PI and Riesz potential} and we get
\begin{align}
K_{X,s}*F(x)\gtrsim_{(X,Q,s)}&  PI(K_{X,s}*f)(x,y),
\end{align}
for every $y>0$.\\
By construction $E^*(f,\epsilon)$ is the projection on the space $X$ of the set $E'(f,\epsilon) \subseteq X \times (0,+\infty)$, so $\forall x \in E^*(f,\epsilon)$ $\exists y(x)>0$ such that $(x,y(x))\in E'(f,\epsilon)$.\\
We use (\ref{Harnack consequence}) to get
\begin{equation}
K_{X,s}*F(x)\gtrsim_{(X,Q,s)}  PI(K_{X,s}*f)(x,y(x)) \gtrsim_{(X,Q)}\epsilon, \quad \forall x \in E^*(f,\epsilon).
\end{equation}
However, $F\in L^p(X)$, so from the definition of capacity we get that
\begin{equation}
C_{K_{X,s},p}(E^*(f,\epsilon))\lesssim_{(X,Q,s,p)} \left(\frac{\| F \|_{L^p(X)}}{\epsilon}\right)^p.
\end{equation}
Using the maximal inequality (\ref{maximal inequality consequence}) we get
\begin{equation}
C_{K_{X,s},p}(E^*(f,\epsilon))\lesssim_{(X,Q,s,p)} \left(\frac{\| f \|_{L^p(X)}}{\epsilon}\right)^p,
\end{equation}
ending the proof.
\endproof
\begin{Lem}\label{Lemma 2}
Let $E \subseteq X \times (0,+\infty)$. Let 
\begin{eqnarray}
f:X\times[0,+\infty) \longrightarrow [0,+\infty)
\end{eqnarray}
 such that 
\begin{enumerate}
\item[1)] \[f(x,y_1)\leq f(x,y_2) \quad \text{for all }y_1\leq y_2,\]
\item[2)] There exists a constant $\alpha\geq 1$ such that for all $x_1,x_2 \in X$, for all $y\geq 0$ we have
\[f(x_1,y)\leq \alpha f(x_2,y).\]
\end{enumerate}
Let $\Omega_{f,x_0}:=\{(x,y) \; | \; d(x,x_0)\leq f(x_0,y)\}$. Then
\begin{equation}
\{x \in X \; | \: \Omega_{f,x} \cap E \neq \emptyset\} \subseteq \bigcup_{x \in E^*} \Omega_{\alpha f,x}(\delta_{E^*}(x)),
\end{equation}
where $\Omega_{f,x_0} (y):=\{x \in X \; | \; (x,y) \in \Omega_{f,x_0} \}$, and $\delta_{E^*}(x)=d\left(x,X \backslash{E^*}\right)$.
\end{Lem}
\proof 
Let $x_0 \in \{x \in X \; | \: \Omega_{f,x} \cap E \neq \emptyset\}$. Since $\Omega_{f,x_0}  \cap E \neq \emptyset$ there exist $(\tilde x, \tilde y) \in \Omega_{f,x_0}  \cap E$. We have
\begin{itemize}
\item $\tilde x \in \Omega_{f,x_0} (\tilde y)$, hence property 2) entails $x_0 \in \Omega_{\alpha f,x} (\tilde y)$.
\item $(\tilde x,\tilde y) \in E$, hence $B_d(\tilde x,\tilde y) \subseteq E^*$ by definition of $E^*$, so $\delta_{E^*}(\tilde x)\geq \tilde y$.
\end{itemize}
The monotonicity of $f$ entails the monotonicity of the regions $\Omega_{\alpha f,x}$, so we have $\Omega_{\alpha f,x}(r_1)\subseteq \Omega_{\alpha f,x}(r_2)$ $ \forall r_1<r_2$, for all $x \in X$, so we get
\begin{equation} 
x_0 \in \Omega_{\alpha f,\tilde x}( \tilde y) \subseteq  \Omega_{\alpha f,\tilde x}( \delta_{E^*}(\tilde x)),
\end{equation}
which entails the required inequality.
\endproof
\begin{Prop}\label{eta star has the properties needed proposition}
Let $C\geq 1$. Let $p>1$, let $\frac{1}{p'}\leq s<1$. Consider the function
\begin{equation}\label{def function in stupid geometric lemma for alpha}
f(x,y):=C \cdot \eta^*_{X,p}(x,y).
\end{equation}
It can be proved that the function $f$ satisfies the following two properties:
\begin{enumerate}
\item[1)]  \[f(x,y_1)\leq f(x,y_2) \quad \text{for all }y_1\leq y_2,\]
\item[2)]  There exists a constant $\tilde\alpha=\tilde \alpha(X,s,p)> 1$ such that for all $x_1,x_2 \in X$, for all $y\geq 0$ we have
\[f(x_1,y)\leq \tilde \alpha  f(x_2,y).\]
\end{enumerate}
\end{Prop}
The proof of the previous proposition follows from Proposition \ref{capacity of ball general formula proposition}.
\begin{Def}
Let $\tilde\alpha=\tilde \alpha(X,s,p)$ denote the constant defined by property 2) in the previous proposition.
\end{Def}
The following proposition will be used to prove the main result of this work.
\begin{Prop}\label{Theorem 9}
Let $p>1$, $\frac{1}{p'}\leq s<1$. Let $\Psi\geq 1$ be the constant defined by Theorem \ref{Quasi additivity Ahlfors space theorem}. If $E \subseteq X \times (0,+\infty)$ is $C_{K_{X,s},p}$-thin at $X \times\{0\}$ then, given
\begin{equation}
\Omega_{x_0,K_{X,s},p,\Psi}:=\left\{ (x,y) \; | \; x \in B_d(x_0,  \Psi \cdot \eta^*_{X,p}(x_0,y))\right\},
\end{equation}
for $x_0 \in X$, we have
\begin{equation}
m\left( \bigcap_{t>0} \left\{ x \in X \; | \; \Omega_{x_0,K_{X,s},p,\Psi} \cap E_t \neq 0\right\}\right) =0.
\end{equation}
\end{Prop}
\proof
By Proposition \ref{eta star has the properties needed proposition} the function
\begin{equation}
f(x,y):=\Psi \cdot \eta^*_{X,p}(x,y).
\end{equation}
satisfies the hypotheses of Lemma \ref{Lemma 2}. So we apply Lemma \ref{Lemma 2} to the region 
\begin{equation}
\Omega_{x_0,K_{X,s},p,\Psi}:=\Omega_{f,x_0}=\{(x,y) \; | \; d(x,x_0)\leq \Psi \cdot \eta^*_{X,p}(x_0,y)\}
\end{equation}
 and we get 
\begin{eqnarray}
\{x \in X \; | \; \Omega_{x,K_{X,s},p,\Psi} \cap E = \emptyset\} &\subseteq&
\bigcup_{x \in E^*}\Omega_{\tilde \alpha f,x}(  \delta_{E^*}(x))=\\
& &\bigcup_{x \in E^*} \Omega_{x,K_{X,s},p,\tilde\alpha\Psi} (\delta_{E^*}(x)).\nonumber
\end{eqnarray}
So by Proposition \ref{Theorem 6}, with $C=\tilde\alpha\Psi$, we get
\begin{equation}\label{step in theorem 9 proof}
m( \{x \in X \; | \; \Omega_{x,K_{X,s},p,\Psi} \cap E \neq \emptyset\} ) \lesssim_{(X,Q,s,p,\tilde\alpha\Psi)}  C_{K_{X,s},p}(E^*).
\end{equation}
Now we apply equation (\ref{step in theorem 9 proof}) to $E=E_t$ and we get
\begin{equation}
m( \{x \in X \; | \; \Omega_{x,K_{X,s},p,\Psi} \cap E_t \neq \emptyset\})  \lesssim_{(X,Q,s,p,\tilde\alpha\Psi)} C_{K_{X,s},p}(E_t^*)\rightarrow 0 \quad \text{as } t \rightarrow 0
\end{equation}
because $E$ is $C_{K_{X,s},p}$-thin at $X \times\{0\}$, so the theorem follows.
\endproof
\subsection{Convergence at the boundary}
We will now prove one more lemma and then we will prove the main results of this work.
\begin{Lem}\label{Lemma 4}
Let $f \in L^p(X)$. Let $\delta>0$. Then there exist $E \subseteq X \times (0,+\infty)$ and $F\subseteq X$ such that
\begin{enumerate}
\item $C_{K_{X,s},p}(E^*)<\delta \quad \text{and } \quad C_{K_{X,s},p}(F)<\delta$.
\item $\forall \epsilon>0$ there exists $r>0$ such that
\begin{equation}
\sup_{x \in X\backslash F} \left( \sup_{P \in B_{\rho}((x,0),r)\backslash E}| PI(K_{X,s}*f)(P)-K_{X,s}*f(x)| \right)<\epsilon.
\end{equation}
\end{enumerate}
\end{Lem}
\proof
Let $\delta>0$. Let $\epsilon>0$. Let $A=A(X,Q,s,p)>0$ be the constant defined by Lemma \ref{Lemma 3} such that
\begin{equation}
C_{K_{X,s},p}(E^*(f,\epsilon))\leq A \cdot \left(\frac{\| g \|_{L^p(X)}}{\epsilon}\right)^p
\end{equation}
for all $g\in L^p_+(X)$, for all $\epsilon>0$. Let $f \in L^p(x)$. Consider
\begin{equation}
f^+:=\max\{f,0\}, \quad f^-:=\max\{-f,0\}.
\end{equation}
Let $\epsilon_1=\epsilon_1(j)>0$ arbitrary to be fixed later. By Lusin's theorem and Urysohn's lemma for all $\epsilon_1>0$, for all $j \in \mathbb N$ there exist $g^+_j,g^-_j \in L^p(X)\cap C_0(X)$ and there exist sets $S^+_j,S^-_j \subseteq X$  such that $m(S^+_j)< \epsilon_1$, $m(S^-_j)< \epsilon_1$, such that $g^+_j \equiv f^+$ in $X\backslash S^+_j$, $g^-_j \equiv f^-$ in $X\backslash S^-_j$ and such that $0\leq g^+_j \leq f^+$ and $0\leq g^-_j \leq f^-$.\\
Indeed, let $j \in \mathbb N$. We apply Lusin's theorem to the function $f^+$ and we get that there exist an open set $S^+_{1,j}$ such that $m(S^+_{1,j})< \frac{\epsilon_1}{2}$ and $f^+$ is continuous in $X\backslash S^+_{1,j}$. Let $S^+_{2,j}$ be an arbitrary closed set such that $S^+_{2,j} \subseteq X\backslash \overline{S^+_{1,j}}$ and $m\left(\left(X\backslash\overline{ S^+_{1,j}}\right)\backslash S^+_{2,j}\right)< \frac{\epsilon_1}{2}$. Such set exists because $X$ is an Ahlfors-regular space, so the measure $m$ is regular. We apply Urysohn's lemma and we get that there exists a function
\begin{equation}
h_j:X\longrightarrow[0,1]
\end{equation}
such that $h_j$ is continuous, $h_j \equiv 1$ on $S^+_{2,j}$, $h_j \equiv 0$ on $\overline{S^+_{1,j}}$. We define
\begin{equation}
g^+_j:=f^+ \cdot h_j,
\end{equation}
and $S^+_j:= X \backslash S^+_{2,j}$.\\
By construction $g^+_j \in L^p(X)\cap C_0(X)$ (because $X$ is compact), and we have $g^+_j \equiv f^+$ in $X\backslash S^+_j$, and $0\leq g^+_j \leq f^+$. Moreover, by construction we have
\begin{equation}
m(S^+_j)=m(X \backslash S^+_{2,j})\leq m\left(\left(X\backslash\overline{ S^+_{1,j}}\right)\backslash S^+_{2,j}\right)+m(S^+_{1,j})\leq \epsilon_1,
\end{equation}
proving the claim for the function $f^+$. We repeat the same argument for $f^-$ and the claim is proved.\\
We have
\begin{equation}
\|f^+-g^+_j\|_{L^p(X)}=\left(\int_{ S^+_j} |f^+-g^+_j|^p dm\right)^{\frac 1 p} \longrightarrow 0 \quad \text{as } \epsilon_1 \rightarrow 0.
\end{equation}
Now we repeat the same argument for $f^-$ and choose $\epsilon_1=\epsilon_1(\delta, A, j,p)$ small enough such that
\begin{equation}
\|f^+-g^+_j\|_{L^p(X)} \leq 2^{-j} \left( \frac{2^{-j} \delta}{2A}\right)^{\frac 1 p},
\end{equation}
\begin{equation}
\|f^--g^-_j\|_{L^p(X)} \leq 2^{-j} \left( \frac{2^{-j} \delta}{2A}\right)^{\frac 1 p}.
\end{equation}
We define
\begin{equation}
E_{+,j}:=E(f^+-g^+_j,2^{-j})=\{(x,y) \in X \times (0,+\infty) \; | \; PI(K_{X,s} * (f^+-g^+_j))(x,y)>2^{-j}\},
\end{equation}
\begin{equation}
E_{-,j}:=E(f^--g^-_j,2^{-j})=\{(x,y) \in X \times (0,+\infty) \; | \; PI(K_{X,s} * (f^--g^-_j))(x,y)>2^{-j}\}.
\end{equation}
Following definition \ref{E E_t and E^{M,*} definition} we consider
\begin{equation}
 E_{+,j}^*:= \bigcup_{(x,y) \in E_{+,j}} B_d(x,y)\subseteq X,
\end{equation}
\begin{equation}
 E_{-,j}^*:= \bigcup_{(x,y) \in E_{-,j}} B_d(x,y)\subseteq X.
\end{equation}
By construction we may apply Lemma \ref{Lemma 3} and we get
\begin{equation}\label{how big is capacity of E_j plus *}
C_{K_{X,s},p}( E_{+,j}^*)\leq A \left( \frac{\| f^+-g^+_j\|_{L^p(X)}}{2^{-j}}\right)^p\leq A \left( \frac{2^{-j}}{2^{-j}} \left( \frac{2^{-j} \delta}{2A}\right)^{\frac 1 p}\right)^p=2^{-j}\frac{\delta}{2},
\end{equation}
\begin{equation}\label{how big is capacity of E_j minus *}
C_{K_{X,s},p}( E_{-,j}^*)\leq A \left( \frac{\| f^--g^-_j\|_{L^p(X)}}{2^{-j}}\right)^p\leq A \left( \frac{2^{-j}}{2^{-j}} \left( \frac{2^{-j} \delta}{2A}\right)^{\frac 1 p}\right)^p=2^{-j}\frac{\delta}{2}.
\end{equation}
Let us define
\begin{equation}
E:=\bigcup_{j=1}^{+\infty} E_{-,j} \cup \bigcup_{j=1}^{+\infty}E_{+,j}.
\end{equation}
By construction
\begin{equation}
E^*=  \bigcup_{(x,y) \in E} B_d(x,y)=\bigcup_{j=1}^{+\infty} E^*_{-,j} \cup \bigcup_{j=1}^{+\infty}E^*_{+,j}.
\end{equation}
By the subadditivity of the capacity and by equations (\ref{how big is capacity of E_j plus *}) and (\ref{how big is capacity of E_j minus *}) we get
\begin{equation}
C_{K_{X,s},p}(E^*)\leq\sum_{j=1}^{+\infty} 2^{-j}\frac{\delta}{2}+\sum_{j=1}^{+\infty} 2^{-j}\frac{\delta}{2}=\delta.
\end{equation}
Now we define
\begin{equation}
F_{+,j}:=\{x \in X \; | \; K_{X,s}*(f^+-g^+_j) \geq 2^{-j}\},
\end{equation}
\begin{equation}
F_{-,j}:=\{x \in X \; | \; K_{X,s}*(f^--g^+_j) \geq 2^{-j}\}.
\end{equation}
For all $x \in F_{+,j}$ we have
\begin{equation}
 K_{X,s}*\left(\frac{f^+-g^+_j }{2^{-j}}\right)(x)\geq 1,
\end{equation}
so by the definition of capacity we get
\begin{equation}
C_{K_{X,s},p}(F_{+,j})\leq \| 2^j(f^+-g^+_j)\|^p_{L^p_{(X)}}\leq \frac{2^{-j} \delta}{2A}.
\end{equation}
Up to multiplying $\delta$ by the constant $A$ we can reformulate the last equation to get
\begin{equation}\label{how big is capacity of F_j plus}
C_{K_{X,s},p}(F_{+,j})\leq 2^{-j}\frac{ \delta}{2}.
\end{equation}
By the same argument we also get
\begin{equation}\label{how big is capacity of F_j minus}
C_{K_{X,s},p}(F_{-,j})\leq  2^{-j}\frac{ \delta}{2}.
\end{equation}
Let us define
\begin{equation}
F:=\bigcup_{j=1}^{+\infty} F_{-,j} \cup \bigcup_{j=1}^{+\infty}F_{+,j}.
\end{equation}
By the subadditivity of the capacity and by equations (\ref{how big is capacity of F_j plus}) and (\ref{how big is capacity of F_j minus}) we get
\begin{equation}
C_{K_{X,s},p}(F)\leq\sum_{j=1}^{+\infty} 2^{-j}\frac{\delta}{2}+\sum_{j=1}^{+\infty} 2^{-j}\frac{\delta}{2}=\delta.
\end{equation}
We define
\begin{equation}
g_j:=g^+_j-g^-_j.
\end{equation}
By linearity we get
\begin{align}
PI(K_{X,s} * f)(x,y)-&PI(K_{X,s} * g_j)(x,y)=PI(K_{X,s} * (f-g_j))(x,y)=\\
 & PI(K_{X,s} * (f^+-f^- -g_j^++g_j^-))(x,y)=\nonumber\\
& PI(K_{X,s} * (f^+-g_j^+))(x,y)-PI(K_{X,s} * (f^- -g_j^-))(x,y),\nonumber
\end{align}
so by triangle inequality we get
\begin{align}\label{triangle inequality for PI(K*(f-g_j))}
|PI(K_{X,s} * f)(x,y)-PI(K_{X,s} * g_j)(x,y)|\leq& |PI(K_{X,s} * (f^+-g_j^+))(x,y)|+\\
&|PI(K_{X,s} * (f^- -g_j^-))(x,y)|.\nonumber
\end{align}
Let $j\geq 1$. Suppose $(x,y)\in X \times(0,+\infty) \backslash E$. Then $(x,y) \not \in E_{+,j} \cup E_{-,j}$, so by definition of $E_{+,j}$ and $E_{-,j}$ we get
\begin{equation}\label{how big is the difference of PI(K) of f and g_j}
|PI(K_{X,s} * (f^+-g_j^+))(x,y)|<2^{-j}, \quad \text{and} \quad |PI(K_{X,s} * (f^- -g_j^-))(x,y)|<2^{-j}.
\end{equation}
From equations (\ref{triangle inequality for PI(K*(f-g_j))}) and (\ref{how big is the difference of PI(K) of f and g_j}) we get
\begin{equation}
|PI(K_{X,s} * f)(x,y)-PI(K_{X,s} * g_j)(x,y)|\leq 2^{-j+1}
\end{equation}
for all $(x,y) \in X \times(0,+\infty) \backslash E$. So we proved that
\begin{equation}\label{uniform convergence of poisson integral of potential}
PI(K_{X,s} * g_j)\overset{j \rightarrow +\infty}{\longrightarrow} PI(K_{X,s} * f) \quad \text{uniformly on } X \times(0,+\infty) \backslash E.
\end{equation}
Suppose $x \in X \backslash F$.  By additivity of the potential and by triangle inequality we get that
\begin{equation}
|K_{X,s} * f(x)-K_{X,s} *  g_j(x)| \leq |K_{X,s} * f^+(x)-K_{X,s} *  g^+_j(x)|+ |K_{X,s} * f^-(x)-K_{X,s} *  g^-_j(x)|.
\end{equation}
However, by the definition of $F_{+,j}$ and $F_{-,j}$, it follows that
\begin{equation}
|K_{X,s} * f(x)-K_{X,s} *  g_j(x)| \leq 2^{-j+1}
\end{equation}
for all $x \in X \backslash F$. So we proved that
\begin{equation}\label{uniform convergence of potential}
K_{X,s} * g_j\overset{j \rightarrow +\infty}{\longrightarrow} K_{X,s} * f \quad \text{uniformly on } X  \backslash F.
\end{equation}
So from (\ref{uniform convergence of poisson integral of potential}) and (\ref{uniform convergence of potential}) we get that there exists $j_0 \in \mathbb N$ such that
\begin{equation}\label{1 last inequality of lemma 4}
\sup_{x \in X \backslash F} |K_{X,s} * f(x)-K_{X,s} *  g_{j_0}(x)| \leq \frac{\epsilon}{3},
\end{equation}
\begin{equation}\label{2 last inequality of lemma 4}
\sup_{P \in X\times(0,+\infty) \backslash E} |PI(K_{X,s} * f)(P)-PI(K_{X,s} *  g_{j_0})(P)| \leq \frac{\epsilon}{3}.
\end{equation}
By construction $g_{j_0} \in C_0(X)$, so $K_{X,s} *  g_{j_0} \in C(X)$ so we apply Lemma \ref{Uniform continuity of the Riesz potential} to the function $K_{X,s} *  g_{j_0}$ and we get that there exists $r>0$ such that
\begin{equation}\label{3 last inequality of lemma 4}
\sup_{x \in X}\left(\sup_{P \in B_{\rho}((x,0),r)}|PI(K_{X,s} * g_{j_0})(P)-K_{X,s} *  g_{j_0}(x)| \right)\leq \frac{\epsilon}{3}.
\end{equation}
The statement follows from (\ref{1 last inequality of lemma 4}), (\ref{2 last inequality of lemma 4}) and (\ref{3 last inequality of lemma 4}) by triangle inequality.
\endproof
\  \\
We are now going to prove the main results of this work.
\begin{Teo}[Non tangential convergence for the Riesz potential]\label{Theorem 10}
Let $(X,d,m)$ be a compact Ahlfors-regular space. Let $f \in L^p(X)$. Then $\exists E \subseteq X \times (0,+\infty)$ such that $E$ is  $C_{K_{X,s},p}$-thin at  $X \times\{0\}$
and
\begin{equation}\label{theorem 10 thesis}
\lim_{\underset{\underset{(x,y)\not \in E}{x \in B_d(x_0, y)}}{(x,y)=P \rightarrow (x_0,0)}} PI(K_{X,s} * f)(P)=K_{X,s} * f(x_0)
\end{equation}
for $C_{K_{X,s},p}$-almost everywhere $x_0 \in X$, i.e. $\exists F \subset X$ such that $C_{K_{X,s},p}(F)=0$ and (\ref{theorem 10 thesis}) holds $\forall x_0 \in X\backslash F$.
\end{Teo}
\proof
Let $f \in L^p(X)$. Let $\epsilon_j>0$ be a sequence such that $\epsilon_j \downarrow 0$ as $j \rightarrow +\infty$. By Lemma \ref{Lemma 4} there exist $E_j \subseteq X \times(0,+\infty)$, $F_j \subseteq X$ and $r_j \downarrow 0$ such that
\begin{equation}\label{first consequence from lemma 4}
\sum_j C_{K_{X,s},p}(E_j^*) < +\infty, \quad \text{and} \quad \sum_j C_{K_{X,s},p}(F_j) <+\infty,
\end{equation}
\begin{equation}\label{consequence of lemma 4}
\sup_{x \in X\backslash F_j} \left( \sup_{P \in B_{\rho}((x,0),r_j)\backslash E_j}| PI(K_{X,s}*f)(P)-K_{X,s}*f(x)| \right)<\epsilon_j.
\end{equation}
Let us choose $t_j \downarrow 0$ such that
\begin{equation}
t_j < r_{j+1},
\end{equation}
\begin{equation}\label{intersection property with balls and half planes in theorem 10}
\{(x,y)\; | \; x \in B_{d}(x_0,y)\} \cap B_{\rho}((x_0,0),r_i) \subseteq \bigcup_{j=i}^{+\infty} B_{\rho}((x_0,0),r_j) \cap \{(x,y) \; | \; y \geq t_j\}.
\end{equation}
A sequence $t_j$ with such properties exists thanks to the definition of the distance $\rho$.\\
Let us define
\begin{equation}\label{definition of E'_j}
E'_j = E_j \cap \{ (x,y) \; | \; y \geq t_j\},
\end{equation}
and 
\begin{equation}
E=\bigcup_{j} E'_j.
\end{equation}
Let $i \in \mathbb N$ be a fixed index. By construction $t_i>t_j$ for all $j>i$, and
\begin{equation}
E_{t_i}=\{(x,y) \in E \; | \; y<t_i\} = \bigg\{(x,y) \in \bigcup_{j=1}^{+\infty}E'_j \; \bigg | \; y<t_i\bigg\} = \bigcup_{j=1}^{+\infty}\bigg\{(x,y) \in E'_j \; \bigg | \; y<t_i\bigg\}.
\end{equation}
However from (\ref{definition of E'_j}) we have
\begin{equation}
\bigg\{(x,y) \in E'_j \; \bigg | \; y<t_i\bigg\}=\emptyset \quad \forall j \leq i,
\end{equation}
so we get
\begin{equation}\label{E_{t_i} is in E_j for j>i}
E_{t_i}=\bigcup_{j=i}^{+\infty}\bigg\{(x,y) \in E'_j \; \bigg | \; y<t_i\bigg\} \subseteq \bigcup_{j=i}^{+\infty} E_j.
\end{equation}
From (\ref{definition of E^{M,*}}), (\ref{E_{t_i} is in E_j for j>i}) and by subadditivity of the capacity we get
\begin{equation}
C_{K_{X,s},p}(E^*_{t_i})\leq C_{K_{X,s},p}\left( \bigcup_{j=i}^{+\infty} E^*_j\right)\leq \sum_{j=i}^{+\infty}C_{K_{X,s},p}\left(  E^*_j\right)
\end{equation}
Equation (\ref{first consequence from lemma 4}) entails
\begin{equation}
C_{K_{X,s},p}(E_{t_i}^*)\leq \sum_{j=i}^{+\infty}C_{K_{X,s},p}\left(E_j^*\right) \overset{i \rightarrow +\infty}{\longrightarrow} 0,
\end{equation}
so, using the monotonocity of the capacity, we get
\begin{equation}
C_{K_{X,s},p}(E_{t}^*) \overset{t \rightarrow 0}{\longrightarrow} 0,
\end{equation}
i.e. $E$ is a $C_{K_{X,s},p}$-thin set at $X\times\{0\}$.\\
Now we define
\begin{equation}
F:=\bigcap_{i=1}^{+\infty} \bigcup_{j=i}^{+\infty} F_j.
\end{equation}
By (\ref{first consequence from lemma 4}) we have
\begin{equation}
C_{K_{X,s},p}(F)\leq C_{K_{X,s},p}\left(\bigcup_{j=i}^{+\infty} F_j\right)\leq \sum_{j=i}^{+\infty} C_{K_{X,s},p}(F_j)\longrightarrow 0\quad \text{as }j \rightarrow +\infty,
\end{equation}
so we have $C_{K_{X,s},p}(F)=0$.\\
Let $x_0 \in X \backslash F$. By definition $\exists j_0=j_0(x_0)$ such that $x_0 \in F_j$, $\forall j \geq j_0$.
From (\ref{consequence of lemma 4}) we get 
\begin{equation}
\sup_{P \in B_{\rho}((x_0,0),r_j)\backslash E_j}| PI(K_{X,s}*f)(P)-K_{X,s}*f(x_0)|\leq \epsilon_j \quad \forall j \geq j_0.
\end{equation}
Let $i \geq j_0$. Using (\ref{intersection property with balls and half planes in theorem 10}) and the definition of $E$ we get
\begin{align}
&\sup_{  (\{P=( x,y)\; | \; x \in B_{d}(x_0,y)\}\cap B_{\rho}((x_0,0),r_i))\backslash E}| PI(K_{X,s}*f)(P)-K_{X,s}*f(x_0)|\leq\\
&\sup_{j\geq i}\bigg(\sup_{P \in (B_{\rho}((x_0,0),r_j)\cap\{(x,y)\; | \; y\geq t_j\})\backslash E}| PI(K_{X,s}*f)(P)-K_{X,s}*f(x_0)|\bigg)\leq\nonumber\\
&\sup_{j\geq i}\bigg(\sup_{P \in (B_{\rho}((x_0,0),r_j)\cap\{(x,y)\; | \; y\geq t_j\})\backslash E_j}| PI(K_{X,s}*f)(P)-K_{X,s}*f(x_0)|\bigg)\leq\nonumber\\
&\sup_{j\geq i} \epsilon_j = \epsilon_i \longrightarrow 0 \quad \text{as } i\rightarrow +\infty, \nonumber
\end{align}
which entails the thesis, finishing the proof.
\endproof
\begin{Teo}[Tangential convergence for the Riesz potential]\label{Theorem 11} 
Let $(X,d,m)$ be a compact Ahlfors-regular space. Let $p>1$, let $\frac{1}{p'}\leq s<1$. Let $\Psi\geq1$ be the constant defined by Theorem \ref{Quasi additivity Ahlfors space theorem}. Consider the region
\begin{equation}
\Omega_{x_0,K_{X,s},p,\Psi}:=\left\{ (x,y) \; | \; x \in B_d(x_0, \Psi \cdot \eta^*_{X,p}(x_0, y))\right\}.
\end{equation}
Let $f \in L^p(X)$. Then
\begin{equation}
\lim_{\underset{P \in \Omega_{x_0,K_{X,s},p,\Psi}}{(x,y)=P \rightarrow (x_0,0)}}PI(K_{X,s} * f)(P)=K_{X,s} * f(x_0)
\end{equation}
for $m$-almost all $x_0 \in X$.
\end{Teo}
\proof
The proof is similar to the one of Theorem \ref{Theorem 10}.
Let $f \in L^p(X)$. Let $\epsilon_j>0$ be a sequence such that $\epsilon_j \downarrow 0$ as $j \rightarrow +\infty$. By Lemma \ref{Lemma 4} there exist $E_j \subseteq X \times(0,+\infty)$, $F_j \subseteq X$ and $r_j \downarrow 0$ such that
\begin{equation}\label{first consequence from lemma 4 theorem 11}
\sum_j C_{K_{X,s},p}(E_j^*) < +\infty, \quad \text{and} \quad \sum_j C_{K_{X,s},p}(F_j) <+\infty,
\end{equation}
\begin{equation}\label{consequence of lemma 4 theorem 11}
\sup_{x \in X\backslash F_j} \left( \sup_{P \in B_{\rho}((x,0),r_j)\backslash E_j}| PI(K_{X,s}*f)(P)-K_{X,s}*f(x)| \right)<\epsilon_j.
\end{equation}
Let us choose $t_j \downarrow 0$ such that
\begin{equation}
t_j < r_{j+1},
\end{equation}
\begin{equation}\label{intersection property with balls and half planes in theorem 11}
\{(x,y)\; | \; x \in B_{d}(x_0, \Psi \cdot \eta^*_{X,p}(x_0, y)))\} \cap B_{\rho}((x_0,0),r_i) \subseteq \bigcup_{j=i}^{+\infty} B_{\rho}((x_0,0),r_j) \cap \{(x,y) \; | \; y \geq t_j\}.
\end{equation}
A sequence $t_j$ with such properties exists thanks to the definition of the distance $\rho$.\\
Let us define
\begin{equation}
E'_j = E_j \cap \{ (x,y) \; | \; y \geq t_j\},
\end{equation}
and 
\begin{equation}
E=\bigcup_{j} E'_j.
\end{equation}
By the same argument used in the proof of Theorem \ref{Theorem 10} we get
\begin{equation}
C_{K_{X,s},p}(E_{t}^*) \overset{t \rightarrow 0}{\longrightarrow} 0,
\end{equation}
i.e. $E$ is a $C_{K_{X,s},p}$-thin set at $X\times\{0\}$.\\
Now we define
\begin{equation}
F:=\bigcap_{i=1}^{+\infty} \bigcup_{j=i}^{+\infty} F_j,
\end{equation}
and
\begin{equation}
S:= F \cup \left( \bigcap_{t>0} \left\{ x \in X \; | \; \Omega_{x,K_{X,s},p,\Psi} \cap E_t \neq \emptyset\right\}\right).
\end{equation}
By the same argument in the proof of Theorem \ref{Theorem 10} we have $C_{K_{X,s},p}(F)=0$, and hence $m(F)=0$.
By Proposition \ref{Theorem 9} we have
\begin{equation}
m\left( \bigcap_{t>0} \left\{ x \in X \; | \; \Omega_{x,K_{X,s},p,\Psi} \cap E_t \neq \emptyset\right\}\right)=0,
\end{equation}
so we proved
\begin{equation}
m(S)=0.
\end{equation}
Let $x_0 \in X \backslash S$. By definition of $S$ there exists $j_0=j_0(x_0)$ such that $x_0 \in F_j$, $\forall j \geq j_0$, and such that
\begin{equation}\label{no intersection tangential region in theorem 11}
\Omega_{x_0,K_{X,s},p,\Psi} \cap E_{t}=\emptyset \quad \text{for all }t<t_{j_0}.
\end{equation}
From (\ref{consequence of lemma 4 theorem 11}) we get 
\begin{equation}
\sup_{P \in B_{\rho}((x_0,0),r_j)\backslash E_j}| PI(K_{X,s}*f)(P)-K_{X,s}*f(x_0)|\leq \epsilon_j \quad \forall j \geq j_0.
\end{equation}
Let $i \geq j_0$. Using (\ref{intersection property with balls and half planes in theorem 11}), (\ref{no intersection tangential region in theorem 11}) and the definition of $E$ we get
\begin{align}
&\sup_{ \underset{y<t_{j_0}}{P=(x,y) \in  \Omega_{x_0,K_{X,s},p,\Psi}}}| PI(K_{X,s}*f)(P)-K_{X,s}*f(x_0)|\leq\\
&\sup_{ \underset{y<t_{j_0}} {P=(x,y) \in  \Omega_{x_0,K_{X,s},p,\Psi}\backslash  E} }| PI(K_{X,s}*f)(P)-K_{X,s}*f(x_0)|\leq\\
&\sup_{j\geq i}\bigg(\sup_{P \in (B_{\rho}((x_0,0),r_j)\cap\{(x,y)\; | \; y\geq t_j\})\backslash E}| PI(K_{X,s}*f)(P)-K_{X,s}*f(x_0)|\bigg)\leq\nonumber\\
&\sup_{j\geq i}\bigg(\sup_{P \in (B_{\rho}((x_0,0),r_j)\cap\{(x,y)\; | \; y\geq t_j\})\backslash E_j}| PI(K_{X,s}*f)(P)-K_{X,s}*f(x_0)|\bigg)\leq\nonumber\\
&\sup_{j\geq i} \epsilon_j = \epsilon_i \longrightarrow 0 \quad \text{as } i\rightarrow +\infty, \nonumber
\end{align}
which entails the thesis, finishing the proof.
\endproof
\begin{Oss}\label{final observation}
The region $\Omega_{x_0,K_{X,s},p,\Psi}$ has the following properties:
\begin{itemize}
\item Case $\frac{1}{p'}<s<1$. There exist constants $0<\tilde C_1< \tilde C_2<+\infty$, $y_0>0$ that depend only on $X$, $s$ and $p$ such that
\begin{align*}
\Omega_{x_0,K_{X,s},p,\Psi}\cap\{(x,y) \; | \; y<y_0\}\supseteq&\bigg\{(x,y)   \mid  y<y_0, \; x \in B_d\left(x_0,\tilde C_1 \cdot y^{p\left(s-\frac{1}{p'}\right)}\right)\bigg\}, \\
 \Omega_{x_0,K_{X,s},p,\Psi}\cap\{(x,y) \; | \; y<y_0\} \subseteq& \bigg\{(x,y)\mid y<y_0, \; x \in B_d\left(x_0,\tilde C_2 \cdot y^{p\left(s-\frac{1}{p'}\right)}\right)\bigg\},
\end{align*}
\item Case $s=\frac{1}{p'}$. There exist constants   $0<\tilde K_1 <\tilde K_2<+\infty$, $y_0>0$ that depend only on $X$, $s$ and $p$  such that
\begin{align*}
\Omega_{x_0,K_{X,s},p,\Psi}\cap\{(x,y) \; | \; y<y_0\}\supseteq&\bigg\{(x,y)  \; | \;y<y_0, \; x \in B_d\left(x_0,\tilde K_1 \cdot \log\left(\frac{1}{y}\right)^{-Q}\right)\bigg\}, \\
 \Omega_{x_0,K_{X,s},p,\Psi}\cap\{(x,y) \; | \; y<y_0\} \subseteq&\bigg \{(x,y)  \; | \;y<y_0, \; x \in B_d\left(x_0,\tilde K_2 \cdot \log\left(\frac{1}{y}\right)^{-Q}\right)\bigg\},
\end{align*}
here $Q$ is the dimension of the Ahlfors $Q$-regular space $X$.
\end{itemize}
\end{Oss}
The proof of this observation follows from Proposition \ref{capacity of ball general formula proposition}.
\section*{Acknowledgments}
\addcontentsline{toc}{section}{Acknowledgements}
We would like to thank Andrea Carbonaro and Nikolaos Chalmoukis for our useful talks about the topics in this article. The topics in this article constitute a part of the autor's PhD thesis ``Potential theory on metric spaces", which was written under the supervision of Nicola Arcozzi.\\


\begin{thebibliography}{13}
\addcontentsline{toc}{section}{References}
\bibitem{AH} D. R. Adams; L. I. Hedberg. \textit{Function Spaces and Potential Theory}, Springer, 1996.
\bibitem{A1} H. Aikawa. \textit{Quasiadditivity of capacity and minimal thinness}. (English summary) Ann. Acad. Sci. Fenn. Ser. A I Math. 18 (1993), no. 1, 65--75. 31B15 (31B25)
\bibitem{AB} H. Aikawa; A.A. Borichev. \textit{Quasiadditivity and measure property of capacity and the tangential boundary behavior of harmonic functions.} Potential Theory - ICPT 94: Proceedings of the International Conference on Potential Theory held in Kouty, Czech Republic, August 13-20, 1994, edited by Josef Kral, Jaroslav Lukes, Ivan Netuka and Jiri Vesely, Berlin, Boston: De Gruyter, 2011, pp. 219--228. https://doi.org/10.1515/9783110818574.219
\bibitem{AE} H. Aikawa; M. Essen. \textit{Potential Theory - Selected Topics}, Lecture Notes in Mathematics, 1633. Springer, 1996.
\bibitem{ARSW1} N. Arcozzi; R. Rochberg; E. Sawyer; B. D. Wick. \textit{Potential Theory on Trees, Graphs and Ahlfors-Regular Metric Spaces.} Potential Analysis 41, 317--366 (2014). https://doi.org/10.1007/s11118-013-9371-8
\bibitem{Cartan} H. Cartan. \textit{Theorie generale du balayage en potentiel newtonien.} Ann. Univ. Grenoble, Math. Phys. 22 (1946) 221--280.
\bibitem{Choquet} G. Choquet. \textit{La naissance de la th\'eorie des capacit\'es: r\'eflexion sur une exp\'erience personnelle.} Comptes rendus de l'Académie des sciences. Série générale, La Vie des sciences, vol. 3, n. 4, 1986, pp. 385--397.
\bibitem{Frostman} O. Frostman. \textit{Potentiel d'equilibre et capacite des ensembles avec quelques applications il la theorie des fonctions}. Medd. Lunds Univ. Mat. Sem. 3 (1935) 1--118. 
\bibitem{Gauss} C. F. Gauss. \textit{Allgemeine Lehrs\"atze in Beziehung auf die im verkehrten Verhaltnisse des Quadrats der Entfernung wirkenden Anziehungs und Abstossungs Kr\"afte.} Resultate aus den Beobachtungen des magnetischen Vereins im Jahre 1839, Leipzig, 1840.
\bibitem{HK} J. Heinonen; P. Koskela. \textit{Quasiconformal maps in metric spaces with controlled geometry}. Acta Mathematica, Acta Math. 181(1), 1--61, 1998.
\bibitem{Mazya} V. Maz'ya. \textit{Sobolev spaces with applications to elliptic partial differential equations.} Second, revised and augmented edition. Grundlehren der Mathematischen Wissenschaften [Fundamental Principles of Mathematical Sciences], 342. Springer, Heidelberg, 2011.
\bibitem{NRS} A. Nagel; W. Rudin, J.H. Shapiro. \textit{Tangential Boundary Behavior of Function in Dirichlet-Type Spaces.} Annals of Mathematics , Sep., 1982, Second Series, Vol. 116, No. 2 (Sep., 1982), pp. 331--360
\bibitem{SW} E.M. Stein; G. Weiss. \textit{Introduction to Fourier Analysis on Euclidean Spaces.} (PMS-32). Princeton University Press, 1971. http://www.jstor.org/stable/j.ctt1bpm9w6.
\end{thebibliography}
\end{document}